\documentclass[11pt]{article}

\usepackage{latexsym,amsmath,amsfonts,amssymb,amstext,theorem,euscript,array,enumerate,mathrsfs}
\usepackage{color}
\usepackage{comment}

\newtheorem{Theorem}{Theorem}[section]
\newtheorem{Definition}{Definition}[section]
\newtheorem{Proposition}{Proposition}[section]

\newtheorem{Lemma}{Lemma}[section]

\newtheorem{Remark}{Remark}[section]
\newtheorem{Example}{Example}[section]
\numberwithin{equation}{section}

\def\esssup_#1{\underset{#1}{\mathrm{ess\,sup\, }}}
\def\essinf_#1{\underset{#1}{\mathrm{ess\,inf\, }}}

\def\qed{{\hfill\hbox{\enspace${ \square}$}} \smallskip}
\def\sqr#1#2{{\vcenter{\vbox{\hrule height .#2pt \hbox{\vrule
 width .#2pt height#1pt \kern#1pt \vrule
width .#2pt} \hrule height .#2pt}}}}
\def\square{\mathchoice\sqr54\sqr54\sqr{4.1}3\sqr{3.5}3}

\def\ds{\begin{displaystyle}}
\def\eds{\end{displaystyle}}
\def\dis{\displaystyle }
\def\<{\langle }
\def\>{\rangle }

\def \N{\mathbb{N}}
\def \R{\mathbb{R}}

\def \E{\mathbb{E}}
\def \F{\mathbb{F}}
\def \G{\mathbb{G}}
\def \P{\mathbb{P}}
\def \Q{\mathbb{Q}}

\def \H{\mathbb{H}}

\def \Ac{{\cal A}}
\def \Bc{{\cal B}}

\def \Ec{{\cal E}}
\def \Fc{{\cal F}}

\def \Ic{{\cal I}}
\def \Kc{{\cal K}}

\def \Pc{{\cal P}}

\def \Nc{{\cal N}}
\def \Rc{{\cal R}}
\def \Sc{{\cal S}}
\def \Tc{{\cal T}}
\def \Uc{{\cal U}}
\def \Vc{{\cal V}}

\def \Yc{{\cal Y}}
\def \Vc{{\cal V}}

\def\calb{{\cal B}}
\def\calc{{\cal C}}

\def\calf{{\cal F}}
\def\calg{{\cal G}}
\def\calh{{\cal H}}

\def\calm{{\cal M}}
\def\caln{{\cal N}}

\def\calp{{\cal P}}

\def\calv{{\cal V}}
\def\call{{\cal L}}

\def\bfB{{\bf B}}
\def\bfC{{\bf C}}
\def\bfM{{\bf M}}

\def \eps{\varepsilon}

\def \ep{\hbox{ }\hfill$\Box$}

\def\reff#1{{\rm(\ref{#1})}}

\def\beqs{\begin{eqnarray*}}
\def\enqs{\end{eqnarray*}}
\def\beq{\begin{eqnarray}}
\def\enq{\end{eqnarray}}

\allowdisplaybreaks

\begin{document}

\title{Randomization method and backward SDEs
  for optimal control\\
 of partially observed
 path-dependent  stochastic systems}

\author{Elena BANDINI\thanks{Politecnico di Milano, Dipartimento di Matematica, via Bonardi 9, 20133 Milano, Italy; ENSTA ParisTech, Unit\'e de Math\'ematiques appliqu\'ees, 828, boulevard des Mar\'echaux, F-91120 Palaiseau, France; e-mail: \texttt{elena.bandini@polimi.it}}
\and
Andrea COSSO\thanks{Laboratoire de Probabilit\'es et Mod\`eles Al\'eatoires, CNRS, UMR 7599, Universit\'e Paris Diderot; Politecnico di Milano, Dipartimento di Matematica, via Bonardi 9, 20133 Milano, Italy; e-mail: \texttt{cosso@math.univ-paris-diderot.fr}}
\and
Marco FUHRMAN\thanks{Politecnico di Milano, Dipartimento di Matematica, via Bonardi 9, 20133 Milano, Italy; e-mail: \texttt{marco.fuhrman@polimi.it}}
\and
Huy\^{e}n PHAM\thanks{Laboratoire de Probabilit\'es et Mod\`eles Al\'eatoires, CNRS, UMR 7599, Universit\'e Paris Diderot, and CREST-ENSAE; e-mail: \texttt{pham@math.univ-paris-diderot.fr}}
}

\maketitle

\begin{abstract}
We introduce a suitable backward stochastic differential equation (BSDE) to
represent the value of an optimal control problem with
partial observation for a controlled stochastic equation driven
by Brownian motion. Our model is general enough to
include cases with latent factors
in Mathematical Finance. By a standard reformulation
based on the reference probability method, it also includes the
classical model where the observation process is affected
by a Brownian motion (even in presence of a correlated noise),
a case where a BSDE representation of the value was not available so far.
This approach based on BSDEs  allows for greater generality
beyond the Markovian case, in particular our model
may include  path-dependence in the coefficients (both with respect to the state
and  the control), and does not require any non-degeneracy
condition on the controlled equation. We also discuss the issue of
numerical treatment of the proposed BSDE.

We use a randomization method, previously adopted only for cases
of full observation, and consisting, in a first step,  in replacing the control by an exogenous process
independent of the driving noise and in formulating
an auxiliary (``randomized'')
control problem where optimization is performed over changes of equivalent probability measures affecting the characteristics of the exogenous process.
Our  first  main result is to prove the equivalence between the
original partially observed control problem and
the randomized problem. In a second step we prove that
the latter can be associated by duality to a BSDE, which
then characterizes  the value of the original problem as well.

\end{abstract}

\vspace{5mm}


\vspace{5mm}


\maketitle

\date{}


\section{Introduction}

The main aim of this paper is to prove a representation formula for the value
of a general class of stochastic optimal control problems with
partial observation by means of
an appropriate backward stochastic differential equation (backward SDE, or BSDE).

To motivate our results, let us start with a classical optimal control problems with
partial observation, where we consider an $\R^n$-valued controlled process $X$ solution
to an equation of the form
$$
dX_t = b(X_t, \alpha_t )\,dt + \sigma^1(X_t, \alpha_t )\,dV^1_t
+ \sigma^2(X_t, \alpha_t )\,dV^2_t,
$$
with  initial condition $X_0=x_0$, possibly random. The equation
is driven by two processes $V^1$, $V^2$ which are
independent Wiener processes under some probability $\bar \P$,
and the coefficients depend on a control process $\alpha$. The aim is to maximize a reward
functional of the form
$$
J( \alpha)   =   \bar\E\bigg[\int_t^T f(X_s ,\alpha_s)\,ds + g(X_T)\bigg],
$$
where $\bar\E$ denotes the expectation under $\bar\P$.
In the partial observation problem  the control $\alpha$
is constrained to being adapted to the filtration $\F^W=(\calf^W_t)_{t\ge 0}$
generated by another process $W$, called the observation process.
A standard model, widely used in applications, consists in assuming
that $W$ is defined by the formula
$$
dW_t = h(X_t, \alpha_t )\,dt + dV^2_t,\qquad  W_0=0.
$$
In this problem $b,\sigma^1,\sigma^2,f,g,h$ are given data
satisfying appropriate assumptions. We also introduce the
value
$$
{\text{\Large$\upsilon$}}_0=\sup_{\alpha}J(\alpha),
$$
where $\alpha$ ranges in the class of admissible control processes,
i.e. $\F^W$-progressive processes
with values in some set $A$ of control actions.
A very effective approach to this problem is the so-called  reference probability
method, which consists in introducing, by means of
a Girsanov transformation, a probability $\P$  under which
 $V^1$ and $W$   are
independent Wiener processes. Explicitly, one defines
 $d\P    =  Z_T^{-1} \,d\bar\P$ where the density process $Z$ satisfies
 $$ dZ_t = Z_t\,h(X_t,  \alpha_t )\,dW_t,\qquad Z_0=1. $$
Next, one introduces the process of unnormalized conditional distributions
defined for every test function $\phi:\R^n\to \R$ by the formula
$$
\rho_t(\phi)=\E\,[\phi(X_t)Z_t\mid \calf^W_t]
$$
and proves that $\rho$ is a solution to the so-called
 controlled Zakai  equation:
\begin{equation}\label{Zakai_controlled_intro}
    d\rho_t(\phi)=\rho_t(\call^{ \alpha_t } \phi)\,dt
+ \rho_t(h(\cdot,  \alpha_t )
\phi+\calm^{ \alpha_t } \phi)\,dW_t,
\end{equation}
where
$\call^a \phi=\frac12\,Tr(\sigma\sigma^T(\cdot,a)D^2\phi)+D\phi\, b(\cdot,a)$,
$\calm^a \phi=D\phi \,\sigma^2(\cdot,a)$, $\sigma=(\sigma^1, \sigma^2)$,
and initial condition $\rho_0$ equal to the law of $x_0$.
The reward functional to be
maximized can be re-written as
\begin{equation}\label{Zakai_reward_intro}
J( \alpha )= \E\left[
\int_0^T\rho_t\,(f(\cdot, \alpha_t ))\,dt+
\rho_T\,(g(\cdot))\right].
\end{equation}
Under appropriate assumptions, for every
admissible control process, the equation \eqref{Zakai_controlled_intro}
 has a unique  solution $\rho$ in some class  of $\F^W$-progressive
 processes with values in the set of nonnegative Borel measures on $\R^n$,
and thus
\eqref{Zakai_controlled_intro}-\eqref{Zakai_reward_intro}
can be seen as an optimal control problem with full observation, called
the separated problem,  having
the same value ${\text{\Large$\upsilon$}}_0$ as the original one
(properly reformulated).  Often, conditions are given so that $\rho_t(dx)$
 admits a density $\eta_t(x)$ with respect to the Lebesgue measure on $\R^n$
 and the controlled Zakai equation is then written as an equation for
 $\eta$, considered as a process with values in some Hilbert space, for instance
 the space $L^2(\R^n)$
 or some weighted $L^2$-space. We also mention that,
 as an alternative to the Zakai equation,
 one could use the (controlled version of the) Kushner-Stratonovich equation
 and repeat similar considerations.

However, the new controlled  state   $\rho$, or equivalently $\eta$, is now an
 infinite-dimensional process, which makes the separated control problem
rather challenging, and the subject of intensive study.
It is not possible to give here a satisfactory description of the obtained
results and we will limit ourselves
to a brief sketch of the possible approaches and refer the reader to the treatises
\cite{par88}
and \cite{Be} for more complete results and references.
A first approach consists in the construction of the so-called
Nisio semigroup, strictly related to the dynamic programming principle
(see for instance \cite{DavisVaraiya73} for the case of partial
observation)
which aims at describing the time evolution of the value function
in this Markovian case. Another method is the stochastic maximum principle,
in the sense of Pontryagin, which provides necessary conditions
for the optimality of a control process in terms of an adjoint equation
and can be used to solve successfully the problem in a number of cases:
see  \cite{Be} and \cite{Tang98}.
The Hamilton-Jacobi-Bellman (HJB) equation corresponding to
\eqref{Zakai_controlled_intro}-\eqref{Zakai_reward_intro}
has been studied by viscosity methods in \cite{lio89},
 later significantly extended in
\cite{gozzi_swiech}.
The issue of existence of an optimal control is addressed for instance in
\cite{elketal88}.

It is a remarkable fact that the use of BSDEs in this context is limited
to the adjoint equations in the stochastic maximum principle cited above,
in spite of the fact that BSDEs are used extensively and successfully
in many areas of stochastic optimization to represent directly
the value function. The reason for this can be explained
by looking at the simple case when
$\sigma^1=I$, $\sigma^2=0$, $h(x,a)=h(x)$, i.e. we have the
partially observed controlled system
$$
dX_t =b(X_t, \alpha_t )\,dt + dV^1_t,
\quad
dW_t=
h(X_t)\,dt +dV^2_t.
$$
and the corresponding controlled Zakai equation for the density process $\eta$
is the following stochastic PDE  in $\R^n$
\begin{equation}\label{Zakai_example_controlled_density}
    d\eta_t(x) = \frac12\,\Delta \eta_t(x)\,dt - div\,  (
b(x,  \alpha_t ) \,\eta_t(x)
)\,dt
+ \eta_t(x)\,  h(x)
\,dW_t.
\end{equation}
The standard method to represent the value based on BSDEs fails for such
a control problem, since the diffusion coefficient is degenerate
and the drift lacks the required structural condition,
see for instance \cite{gozswi15} for the infinite-dimensional case:
roughly speaking, an associated BSDE could immediately be written
for equations having the form
$$
d\eta_t(x)=\frac12\,\Delta \eta_t(x)\,dt
+\eta_t(x)\,  h(x)\,\Big[
r(\eta_t(\cdot), t,x,  \alpha_t )
\,dt+dW_t\Big]
$$
for some coefficient $r$. This  often implies, by an application of the
 Girsanov theorem, that the laws of the controlled processes $\eta$
 (depending on the control process $\alpha$)
are all absolutely continuous with respect to the law of the uncontrolled
process corresponding to $r=0$, but this property fails in general
for the solutions to \eqref{Zakai_example_controlled_density}.
The same difficulty can also be seen at the level of the corresponding
HJB equation: in the simple case we are addressing
this equation is of semilinear type, but it does not fall into the class
of PDEs whose solution can be represented by means of an associated BSDE,
see for instance
\cite{PardouxPeng92}
or \cite{gozswi15}. The problem of representing the value function
of the classical partially observed control problem
by means of a suitable BSDE becomes even more difficult
in the general case when the HJB equation
is fully non-linear, and has remained unsolved so far.

It is the purpose of this paper to fill this gap in the existing literature, by introducing
a suitable BSDE whose solution provides such a representation formula.
As a motivation, we note that
methods based on BSDEs have the advantage that they easily generalize
beyond the Markovian framework. As a matter of fact we will be able
to treat control problems where the coefficients in the state equation
and in the reward functional exhibit memory effects both with respect to
state and to control, i.e. their value at some time
$t$ may depend in a rather general way
on the whole trajectory of the state and control processes
on the time interval $[0,t]$. Various models with delay effects, or hereditary
systems, are thus included in our treatment. For these models
there is no direct application of methods which exploit Markovianity,
such as the Nisio semigroup or the HJB equation
(although the Markovian character can be retrieved in a number
of cases after an appropriate
reformulation, which requires however non-trivial efforts
and often introduces additional assumptions).

Another motivation for introducing BSDEs
is the fact that their solutions can be approximated
numerically. This way our result opens perspectives to finding
an effective way to approximate
the value of a partially observed problem, which is a difficult task
due to the infinite-dimensional character of the Zakai
equation and its corresponding HJB equation.

To perform our program we first formulate a general
control problem of the form
\beq \label{dynX1intro}
dX_t^\alpha & = & b_t( X^\alpha, \alpha)\,dt +
\sigma_t( X^\alpha, \alpha)\,dB_t,
\qquad X^\alpha_0=x_0,
\enq
for $t\in [0,T]$, with  reward functional and value defined by
\beq\label{gainvalueintro}
J(\alpha) & = & \E\Big[\int_0^Tf_t(X^\alpha,\alpha)\,dt+g(X^\alpha)\Big],
\qquad {\text{\Large$\upsilon$}}_0=\sup_\alpha J(\alpha)
\enq
where the coefficients
$b,\sigma,f,g$ depend on the whole trajectories of $X$ and $\alpha$
in an non-anticipative way.  The partial observation character is modeled as follows:
in the  Wiener process $B$ we distinguish  two components (possibly multidimensional) and write it in the form
 $B=(V,W)$. We call $W$ the observation process and we require
 the control process $\alpha$   to be adapted to the filtration
 generated by $W$ and taking values in some set $A$,
 so that  the supremum in \eqref{gainvalueintro}
 is taken over such controls. In  section \ref{SubS:ClassicalPr}
 we prove that this model includes
 the classical partial observation problem described
 above as a special case and, moreover, that
this formulation  is general enough to
 include large classes of optimization models
 with latent factors of interest in mathematical finance: see
 section \ref{SubS:Finance} below.

To tackle this problem we will use a \emph{randomization method}, introduced  in  \cite{KP12} for  classical Markovian models, but earlier considered in \cite{KMPZ10} in connection with impulse control and in \cite{bou09}, \cite{EKa}, \cite{EKb} on optimal switching problems. The idea of using the randomization method was inspired by the fact that
it allows to represent (or construct) viscosity solutions to some classes
of fully non-linear PDEs.  Other methods yield similar results, for instance those based on
 the notion of second order BSDEs \cite{STZ} or the theory of $G$-expectations \cite{Peng07}.  It is likely that they might also be successfully applied
to optimal control problems with partial observation. We also note that the randomization method has  already been applied
to a variety of situations, see \cite{FP15} (compare also Remark \ref{comparison_with_FP15} below),
\cite{ChoukrounCosso16}, \cite{CossoFuhrmanPham16}, \cite{FuhrmanPhamZeni16}, \cite{Bandini15} \cite{BandiniFuhrman15} in addition to the references given above.
In order to present this method applied to
the problem \eqref{dynX1intro}-\eqref{gainvalueintro},
we assume for simplicity that $A$ is a subset of a Euclidean space and we take
a finite measure $\lambda$ on $A$ with full support.
Then,   enlarging the original
 probability space if needed, we introduce  a Poisson random measure $\mu$ on $\R_+\times A$ with
intensity $\lambda(da)$ and  independent of the Brownian motion $B$.
Then we consider the stepwise process
 $I$ associated with $\mu$ and replace the control process
 $\alpha$ by $I$, thus arriving at the following dynamics:
\[
\begin{cases}
\vspace{2mm}
\dis
d X_t =  b_t(  X,I)\,dt + \sigma_t(  X,  I)\,dB_t,
 \\
\dis   I_t \ = \ a_0 + \int_0^t\int_A(a-  I_{s-})\, \mu(ds\, da).
\end{cases}
\]
Next we consider an auxiliary optimization problem, called
\emph{randomized} or \emph{dual} problem
(in contrast to the starting optimal control problem with partial observation
which we refer to also as \emph{primal problem}), which consists in
optimizing among equivalent changes of probability measures which only affect the intensity measure of $\mu$ but not the law of $W$.
In the randomized problem, an admissible control is a bounded
 positive   map  $\nu$ defined on $\Omega\times\R_+\times A$,
 which is predictable with respect to the filtration $\F^{W,\mu}$ generated
 by $W$ and $\mu$.
 Given $\nu$,
by means of an absolutely continuous  change of measure of
 Girsanov type
we construct
a probability measure $\P^\nu$ such that
the compensator of $\mu$  is given by $\nu_t(a)\lambda(da)dt$
and
 $B$ remains a Brownian motion under $\P^\nu$.
 Then we introduce the reward and the value as
\beqs
J^\Rc( \nu)  =    \E^{ \nu}
\left[\int_0^Tf_t(  X,  I)\,dt+g(  X)\right],
\qquad
{\text{\Large$\upsilon$}}_0^\Rc =    \sup_{  \nu } J^\Rc(  \nu),
\enqs
where $\E^\nu$ denotes the expectation under $\P^\nu$.
One of our main results  states that the
 two control problems presented above are equivalent, in the sense that they share the same value:
 \begin{equation}\label{equiv_Intro}
    {\text{\Large$\upsilon$}}_0= {\text{\Large$\upsilon$}}_0^\Rc.
 \end{equation}
See Theorems \ref{MainThm} and \ref{MainThm_bis}, where some additional comments can be found.

The reason for this construction is that, as shown
in section \ref{Sec:separandom},  the randomized control problem is associated
  to the following   BSDE  with a sign constraint,
  which then also characterizes  the value function of the initial control problem \reff{gainvalueintro}.
For any bounded measurable functional $\varphi$ on the  space of continuous paths
with values in $\R^n$
define
\beqs
\rho_t (\varphi) &=&  \E\big[ \varphi(X_{\cdot\wedge t}) \mid \Fc_t^{W,\mu} \big],
\enqs
and consider the BSDE
\begin{equation}\label{BSDEconstrained_Intro}
\begin{cases}
\vspace{2mm} \dis Y_t \ = \ \rho_T(g)  + \int_t^T\rho_s(f_s(\cdot,I)) \,ds
+ K_T - K_t - \int_t^TZ_s\,dW_s - \int_t^T\!\int_A U_s(a)\,\mu(ds\,da), \\
\dis U_t(a) \ \le \ 0,
\end{cases}
\end{equation}
In Theorem \ref{Thm:RandomizedFormula}, which is another of our main results,
we prove that there exists a unique minimal solution $(Y,Z,U,K)$ to \eqref{BSDEconstrained_Intro} (i.e. among all solutions we take the minimal one in terms of the $Y$-component)
in a suitable space of stochastic processes adapted to the filtration
$\F^{W,\mu}$, and moreover
\beq
\label{RandomizedFormula_Intro}
 Y_0={\text{\Large$\upsilon$}}^\Rc_0={\text{\Large$\upsilon$}}_0,     \quad
\text{and more generally}\quad
Y_t =  \esssup_{\nu} \E^\nu\Big[ \int_t^T\rho_s (f_s(\cdot,I))\,ds + \rho_T(g) \,\big|\, \calf_t^{W,\mu}\Big].
\enq
The BSDE \eqref{BSDEconstrained_Intro}
  is called  the
{\it randomized  equation}, and corresponds  to the HJB  equation
of the classical Markovian framework.
Note that the introduction of the measure-valued process $\rho$ and its occurrence
in the generator and the terminal condition of the BSDE is reminiscent
of the  separated problem in classical optimal control
with partial observation. We  study in a companion paper \cite{BCFP15} how one can also derive such kind of HJB equation in the context of partially observed Markovian control problems.

We note that
 probabilistic numerical methods have already been designed for BSDEs
 with constraints similar to
  \eqref{BSDEconstrained_Intro} in
  \cite{KLP14} and \cite{KLP15}.   We also propose, in the Markovian case,
an approximation scheme for  \eqref{BSDEconstrained_Intro}
and hence for the value of the partially
observed control problem, leaving however aside a detailed
analysis which is left for future work.

We would like to point out that in our approach
the original partially observed optimal control problem is formulated
in the strong form, i.e. with a fixed probability space.
This is probably a more natural setting, especially in connection
with modeling applications, and
it is  customary
for  the  stochastic maximum principle and for other classes of optimization
problems like optimal stopping and switching.
However,  almost all applications
of BSDE techniques to the search of an optimal continuous control process
are set in the weak formulation, since this avoids some difficulties
(one exception may be found for instance in \cite{FT04}). In spite of that, in the present paper we have chosen to adopt
the strong formulation, at the expense of additional technical difficulties.
Moreover,  we note that
our main results
are stated in a fairly general framework,
allowing for locally Lipschitz coefficients with linear growth
and without any
non-degeneracy condition imposed
on the diffusion coefficient $\sigma$.  In particular, when $\sigma=0$, this includes the case of deterministic control problem with a path-dependent state dynamics and delay on control. Finally, when the diffusion coefficient of the Brownian motion $V$ is zero, meaning that the dynamics of $X$ is driven only by $W$, we are reduced to the case of full observation control problem. Therefore, we have provided  a general equivalence and representation result in a unifying framework embedding several classical cases in optimal control theory and the proofs
 we present are almost entirely self-contained.


The rest of the paper is organized as follows. In Section \ref{S:Formulation} we formulate the general optimal control problem
\eqref{dynX1intro}-\eqref{gainvalueintro}
(the  primal problem) with partial observation and
path-dependence in the state and the control.  We then present two motivating particular cases:
a general optimization model with latent factors and uncontrolled observation process, which finds usual applications in mathematical finance, and
the classical optimal control problem with partial observation
discussed above
(but including also path-dependence). Then, in Section \ref{Randomized} we implement the randomization method
and formulate the randomized optimal control problem associated with the primal problem. We state in  Theorem \ref{MainThm}  the basic equivalence result between the primal and
the randomized pro\-blem.
Section \ref{proofthm}  is entirely devoted to the proof of  Theorem \ref{MainThm}, which requires for both inequalities sharp approximation results and suitable constructions with marked point processes.
We also
extend this result to the case of locally Lipschitz coefficients with linear growth, which is essential in order
to cover the case of the classical partially observed optimal control.
 In Section \ref{Sec:separandom} we show a separation principle for the randomized control problem using nonlinear filtering arguments, and then relate  by duality the separated randomized problem  to a constrained BSDE, which may be viewed consequently as the randomized equation for the primal control problem. Finally in Section \ref{Sec:numerics} we comment
on numerical issues  and propose, in the Markovian case,
an approximation scheme for the value of the partially
observed control problem.

\section{General formulation and applications}
\label{S:Formulation}

\subsection{Basic notation and assumptions}
\label{SubS:Notation}

In the following we will consider controlled stochastic equations  of the form
\beq \label{dynX1}
dX_t^\alpha & = & b_t( X^\alpha, \alpha)\,dt +
\sigma_t( X^\alpha, \alpha)\,dB_t,
\enq
for $t\in [0,T]$, where $T>0$ is a fixed deterministic and finite
terminal time, and gain functionals
\beqs
J(\alpha) & = & \E\Big[\int_0^Tf_t(X^\alpha,\alpha)\,dt+g(X^\alpha)\Big].
\enqs
The initial condition in \eqref{dynX1} is $X_0^\alpha=x_0$, a given
random variable with law denoted $\rho_0$.
Before formulating precise assumptions let us explain informally the meaning of several terms occurring in these
expressions. The controlled process $X^\alpha$ takes values in $\R^n$ while $B$ is a  Wiener process in $\R^{m+d}$. We write $B=(V,W)$ when we need to distinguish  the first $m$ components of $B$ from the other
$d$ components.
The control process, denoted by $\alpha$, takes values in a set $A$ of control actions. The partial observation available to the controller will be described
by imposing that the control process should be adapted to the filtration generated by the process $W$ alone. Our formulation includes path-dependent (or hereditary) systems, i.e. it allows for the presence of memory
effects both on the state and the control. Indeed, the coefficients
$b,\sigma,f,g$ depend on the whole trajectory of $X^\alpha$ and $\alpha$. The dependence
will be non-anticipative, in the sense that their values at time $t$ depend
on the values $X_s^\alpha$ and $\alpha_s$ for $s\in [0,t]$: this is expressed below in a standard way
by requiring that they should be progressive with respect to
the canonical filtration on the space of  paths.

Now let us come to precise assumptions and notations.
Let us denote by $\bfC_n$   the  space  of continuous paths from $[0,T]$ to $\R^n$, equipped  with the usual supremum norm
$\|x\|_{_\infty}$ $=$ $x^*_T$, where we set  $x^*_t$ $:=$ $\sup_{s\in [0,t]}|x(s)|$, for $t$ $\in$ $[0,T]$ and $x$ $\in$ $\bfC_n $. We define the filtration
$(\calc_t^n)_{t\in[ 0,T]}$, where $\calc_t^n$ is the $\sigma$-algebra generated by the canonical coordinate maps $\bfC_n \to \R^n$,
$x(\cdot)\mapsto x(s)$ up to time $t$, namely
\beqs
\calc_t^n  &:=  &  \sigma \{ x(\cdot)\mapsto x(s)\; :\, s\in [0,t] \},
\enqs
and we denote $Prog(\bfC_n )$  the progressive $\sigma$-algebra on $[0,T]\times\bfC_n $ with respect to $(\calc_t^n)$.

We will require that the space of control actions $A$
is a Borel space. We  recall that a  Borel  space   $A$ is a topological space  homeomorphic to a  Borel subset of a Polish space
(some authors use the terminology  Lusin space). When needed, $A$ will be  endowed with its Borel $\sigma$-algebra  $\calb(A)$.
We denote by $\bfM_A$   the  space  of Borel measurable
paths  $a:[0,T]\to  A$, we introduce the filtration
$(\calm_t^A)_{t\in[ 0,T]}$, where $\calm_t^A$ is the $\sigma$-algebra
\beqs
\calm_t^A  &:=  &  \sigma \{ a(\cdot)\mapsto a(s)\; :\, s\in [0,t] \}
\enqs
and
we denote $Prog(\bfC_n \times \bfM_A)$
the progressive $\sigma$-algebra on $[0,T]\times\bfC_n $ with respect to
the filtration $(\calc_t^n\otimes \calm_t^A)_{t\in [0,T]}$.

\vspace{3mm}

\noindent {\bf (A1)}
\begin{itemize}
\item [(i)]  $A$ is a  Borel   space.
\item [(ii)]
The functions $b,\sigma,f$ are defined on $[0,T]\times \bfC_n\times  {\bf M}_A$ with
values in $\R^n$, $\R^{n\times (m+d) }$ and $\R$ respectively,  are assumed
to be $Prog(\bfC_n \times \bfM_A)$-measurable
(see also Remark \ref{A1_misurabilita} below).
\item [(iii)]
 The function $g$ is continuous on $\bfC_n$, with respect to the supremum norm.
The functions $b$,  $\sigma$ and $f$ are assumed to satisfy the following
sequential continuity
condition: whenever $x_k,x\in\bfC_n$, $\alpha_k,\alpha\in
\bfM_A$, $\|x_k-x\|_\infty\to 0$, $\alpha_k(t)\to\alpha(t)$ for $dt$-a.e.
  $t$ $\in$ $[0,T]$ as $k\to\infty$ we have
  $$
b_t(x_k,a_k) \to b_t(x ,a),\;\;  \sigma_t(x_k,a_k)\to \sigma_t(x,a),
\;\; f_t(x_k,a_k) \to f_t(x ,a)
\quad
\text{ for } dt\text{-a.e. }
 t \in [0,T].
 $$
\item [(iv)]
 There exist nonnegative constants  $L$ and $r$ such that
\beq
|b_t(x,a) - b_t(x',a)| + |\sigma_t(x,a)-\sigma_t(x',a)|
& \leq & L   (x-x')^*_t, \label{lipbsig} \\
|b_t(0,a)| + |\sigma_t(0,a)| & \leq & L,  \label{borbsig}
\\
|f_t(x,a)| + |g(x)| & \leq & L \big(1 + \|x\|_{_\infty}^r \big), \label{PolGrowth_f_g}
\enq
for all $(t,x,x',a)$ $\in$ $[0,T]\times\bfC_n \times\bfC_n \times \bfM_A$.
\item [(v)]  $\rho_0$ is a  probability measure  on the Borel subsets of $\R^n$
satisfying $\int_{\R^n}|x|^p\rho_0(dx)<\infty$
for some $p\ge \max(2,2r)$.
\end{itemize}

\vspace{2mm}

\begin{Remark} \label{A1_misurabilita}
{\rm
The measurability condition
{\bf (A1)}-(ii) is assumed because it guarantees  the following property,
which is easily deduced:
\begin{itemize}

\item [(ii)'] Whenever
$(\Omega,\calf,\P)$ is a probability space with a  filtration $\F$, and $\alpha$ and $X^\alpha$  are $\F$-progressive processes
 with values in $A$ and $\R^n$ respectively, then the process
 $(b_t( X^\alpha , \alpha)$,
$\sigma_t( X^\alpha, \alpha )$,
$f_t(X^\alpha,\alpha))_{t\in [0,T]}$ is also
   $\F$-progressive.
\end{itemize}
All the results in this paper still hold, with the same proofs,
if property (ii)' is assumed to hold instead of (ii). There are cases when  (ii)' is easy to be checked directly.

We shall also discuss in paragraph \ref{Sec:locallylip} how the global Lipschitz condition in {\bf (A1)}-(iv) can be weakened to local Lipschitz condition.

We finally note that the function $g$, being continuous,  is also $\calc_T^n$-measurable.
\ep
}
\end{Remark}

\begin{Remark} \label{remDelta}
{\rm
Assumption {\bf (A1)}
allows us to model various memory effects of the control on the state process,
including important and usual cases of delay in the control. For instance
suppose that
$A$ is a bounded  Borel subset  of a Banach space and $\bar b:A\to \R^n$
is Lipschitz continuous.
Then we may  consider a weighted combination of pure delays:
$$
b_t(x,a) = \bar b\left( \sum_{i=1}^q \pi_i(t) a({t-\delta_i })\right),
$$
where $0<\delta_1<\ldots<\delta_q<T$,
$\pi_i$ are bounded measurable real-valued functions
and we use the convention that $\alpha_t=\bar \alpha$
(a fixed element of $A$) if $t<0$.  We may also allow
the delays $\delta_i$ to depend on $t$ in an appropriate way.
Alternatively, we may have
$$
b_t(x,\alpha) = \bar b\left( \int_0^t \pi(t,s)\, a(s)\, ds\right),
$$
with $\pi$ bounded measurable and real-valued.
Note that in the latter case the measurability condition
{\bf (A1)}-(i) fails in general, since the
$\sigma$-algebras $\calm_t^A $ are
determined by a countable number of times, but the
property (i)' in the previous remark is easy to verify.

Clearly, we may address more complicated situations
which are combinations of the two previous
cases and may also include a dependence on the path  $x$.
\ep
}
\end{Remark}

\begin{Remark}
{\rm  We mention that no non-degeneracy assumption on the diffusion coefficient $\sigma$ is imposed, and  in particular, some lines or columns of $\sigma$ may be equal to zero.  We can then consider a priori more general model than \reff{dynX1} by adding dependence of the coefficients $b$, $\sigma$ on another diffusion process $M$, for example an unobserved and uncontrolled  factor (see Application  in subsection \ref{SubS:Finance}).
This generality is only apparent since it can be embedded in a standard way in our framework by considering the enlarged state process $(X,M)$.
\ep
}
\end{Remark}

\begin{Remark}
{\rm
The requirement that $p\ge \max(2,2r)$ in {\bf (A1)}-(v) can be weakened for specific results
in the sequel. For instance, Theorem \ref{MainThm} below
still holds provided we only require $p\ge \max(2,r)$ and Theorem \ref{MainThm_bis}
remains valid provided we have $p\ge  2$, $p>r$. The corresponding proofs remain
unchanged.
}
\end{Remark}

\subsection{Formulation of the partially observed control problem}
\label{Primal}

We assume that
$A,b,\sigma,f,g,\rho_0$
 are given and
satisfy the assumptions {\bf (A1)}. We formulate a control problem fixing a setting
$(\Omega,\calf,\P, \F, V,W, x_0)$, where
 $(\Omega,\calf,\P)$ is a complete probability space with
 a right-continuous and $\P$-complete filtration $\F=(\calf_t)_{t\ge 0}$,
$V$ and $W$ are processes with values in $\R^{m}$ and $\R^{d}$ respectively,
such that $B=(V,W)$ is an $\R^{m+d}$-valued
standard Wiener process with respect to $\F$ and $\P$, and
$x_0$ is an  $\R^n$-valued
random variable, with law $\rho_0$ under $\P$, which is assumed to be
$\calf_0$-measurable.
Note that $V$ and $W$ are also
standard Wiener processes and that
$V$, $W$, $x_0$ are all
independent.

Let us denote $\F^W=(\calf^W_t)_{t\ge 0}$ the right-continuous and $\P$-complete filtration generated by $W$.
An admissible control process is any $\F^W$-progressive process $\alpha$ with values in $A$.
The set of admissible control processes is denoted by $\Ac^W$. The controlled equation has the form
\begin{equation}\label{stateeq}
    dX_t^\alpha \ = \ b_t( X^\alpha, \alpha)\,dt +
\sigma_t( X^\alpha, \alpha)\,dB_t
\end{equation}
on the interval $[0,T]$ with initial condition $X_0^\alpha=x_0$,
and the gain functional is
\begin{equation}\label{gaineq}
J(\alpha) \ = \ \E\left[\int_0^Tf_t(X^\alpha,\alpha)\,dt+g(X^\alpha)\right].
\end{equation}

Since we assume that {\bf (A1)} holds,  by standard results (see e.g. \cite{rowi} Thm V. 11.2,
 or   \cite{jacod_book}  Theorem 14.23),
there exists a unique $\F$-adapted strong solution $X^\alpha$ $=$ $(X_t^\alpha)_{0\leq t \leq T}$   to \eqref{stateeq}
with continuous trajectories and such that (with the same $p$ for which $\E|x_0|^p<\infty$)
$$
    \E\,\Big[\sup_{t\in [0,T]}|X_t^\alpha|^p\Big]   <  \infty.
$$
The stochastic optimal control problem  under partial observation consists in maximizing $J(\alpha)$ over all $\alpha\in\Ac^W$:
\begin{equation}\label{primalvalue}
{\text{\Large$\upsilon$}}_0 = \max_{\alpha\in\Ac^W} J(\alpha).
\end{equation}

\begin{Remark}\label{R:equivprimal}
\begin{em}
Let  $\F^B=(\calf^B_t)_{t\ge 0}$ be the right-continuous and
$\P$-complete filtration generated by $B$. Then $B$ is clearly
an $\F^B$-Brownian motion, the processes
$\alpha$ and $ X^\alpha$ are $\F^B$-progressive and the filtration
$\F$ does not play any role in determining $J(\alpha)$
and {\Large$\upsilon$}$_0$. So we might assume from the beginning that $\F=\F^B$
and even that $\calf=\calf^B_\infty$ whenever convenient, but
in the sequel we keep
the previous framework unless explicitly mentioned.
 \end{em}
\ep
\end{Remark}

\subsection{Two basic applications}
\label{S:Applications}

In this paragraph,  we address   two classical optimal
control problems with partial observation, and we show that they
can be recast in the form outlined in the previous subsection.

\subsubsection{Model with latent factors and uncontrolled observation process}
\label{SubS:Finance}

Let $(\Omega,\calf,\bar\P)$ be a complete probability space with a
right-continuous and $\bar\P$-complete filtration
$\F=(\calf_t)_{t\ge 0}$. Let $V,\bar W$ be  independent standard Wiener processes with respect to $\F$, with values
in $\R^{m}$ and $\R^{d}$ respectively. We assume that a controller, for instance a financial agent, wants to optimize her/his position, described by an $\bar n$-dimensional
stochastic process $\bar X^\alpha$ solution on the interval $[0,T]$ to an equation of the form
\begin{equation}\label{stateeqfinance}
d\bar X_t^\alpha \; = \;  \bar b_t( \bar X^\alpha, M, O, \alpha)\,dt +
\bar \sigma^1_t(\bar X^\alpha, M, O, \alpha) \,d V_t
+ \bar\sigma^2_t( \bar X^\alpha, M,O, \alpha)\,d\bar W_t
\end{equation}
with coefficients $\bar b, \bar\sigma^1, \bar\sigma^2$ defined on $[0,T]\times \bfC_{\bar n + \bar m+d}\times  {\bfM_A }$ valued
in $\R^{\bar n}$, $\R^{\bar n\times m}$, $\R^{\bar n\times d}$ respectively,  and $Prog(\bfC_{\bar n+\bar m+d} \times\bfM_A) $-measurable.
Here the process $M$, valued in $\R^{\bar m}$,
 represents a  latent factor that can influence  the dynamics of
 $\bar X^\alpha$ and is governed by a dynamics of the form:
\beq \label{dynM}
dM_t &=& \bar\beta_t(M) dt + \gamma_t^1(M) dV_t + \gamma_t^2(M) d\bar W_t,
\enq
for some  coefficients $\bar\beta$, $\gamma^1$, $\gamma^2$ defined on $[0,T]\times \bfC_{\bar m}$ valued in
$\R^{\bar m}$, $\R^{\bar m\times m}$, $\R^{\bar m\times d}$ respectively, and  $Prog(\bfC_{\bar m} )$-measurable.  The process $M$ is not directly observed, and actually the agent takes her/his decisions based on a noisy observation represented by a process $O$ in $\R^d$ solution to an equation of the form
\begin{equation}\label{observation}
dO_t =h_t(M,O)\,dt + k_t(O)\,d\bar W_t, \qquad t\in [0,T],
\end{equation}
for some coefficients $h$ and $k$  defined on $[0,T]\times \bfC_{\bar m+d}$ and $[0,T]\times \bfC_{d}$,  $Prog(\bfC_{\bar m+d} )$-measurable and  $Prog(\bfC_{d} )$-measurable,
valued in $\R^d$ and $\R^{d\times d}$  respectively.  For instance, $O_t$ may be related to
the market price of financial risky assets at time $t$.
 We denote $\F^O=(\calf^O_t)_{t\ge 0}$
the right-continuous and $\bar\P$-complete filtration generated by $O$.
An admissible control process, representing
for instance the agent's investment strategy,
is any $\F^O$-progressive process $\alpha$
with values in the Borel space $A$.

The agent wishes to maximize, among all admissible control processes, a  gain functional of the form
\beqs
J(\alpha) &=& \bar\E\Big[\int_0^T\bar f_t(\bar X^\alpha,M,O, \alpha)\,dt+
\bar g(\bar X^\alpha, M,O)\Big],
\enqs
where $\bar\E$ denotes expectation with respect to $\bar\P$, for real-valued coefficients $\bar f,\bar g$ defined on $[0,T]\times \bfC_{\bar n+\bar m+d}\times  {\bfM_A}$ and
$\bfC_{\bar n+\bar m+d}$, $Prog(\bfC_{\bar n+\bar m+d}\times \bfM_A) $-measurable and  $\calc_T^{\bar n+\bar m+d}$-measurable, respectively.

In order to put this problem in the form addressed in the previous subsection we make a change of probability measure and pass from the ``physical'' probability $\bar \P$
to a ``reference'' probability $\P$. Assuming that $k_t(y)$ is invertible for all $t$ $\in$ $[0,T]$ and $y$ $\in$ $\bfC_d$, and that the process $\{k_t^{-1}(O)h_t(M,O),0\leq t\leq T\}$ is bounded, we
define  a process $Z$ setting
\beqs
Z_t^{-1} &=& \exp\left(- \int_0^t k_s(O)^{-1}h_s(M,O)\,d\bar W_s
-\frac12 \int_0^t |k_s(O)^{-1}h_s(M,O)|^2\,ds\right), \;\;\;  t\in [0,T].
\enqs
The process $Z^{-1}$ is a martingale under $\bar \P$, and by the Girsanov theorem, under the probability $\P(d\omega)=Z_T(\omega)^{-1}\bar\P(d\omega)$ the pair $(V,W)$ is a standard
Wiener process in $\R^{d+m}$ with respect to $\F$, where  $W_t$ $=$ $\bar W_t+ \int_0^t k_s(O) ^{-1}\,h_s(M,O)\,ds$,  $t\in [0,T]$.
We  denote by  $\F^W=(\calf^W_t)_{t\in[ 0,T]}$ the right-continuous and $\P$-complete filtration generated by $W$, and see that  the observation process $O$ is a solution under $\P$ to the equation:
\begin{equation}\label{statouno}
dO_t \; = \;  k_t(O)\,d W_t.
\end{equation}
Assuming a Lipschitz condition on $k$, i.e. there exists a constant $K$ such that
\beqs
|k_t(y) - k_t(y^1)| &\leq&  K  (y-y^1)^*_t,
\enqs
for all $(t,y,y^1 )$ $\in$ $[0,T]\times\bfC_d \times\bfC_d$, we deduce that $O$ must be $\F^W$-adapted and therefore that $\calf^O_t\subset\calf^W_t$ for $t\in[ 0,T]$. On the other hand, since
$W_t$ $=$ $\int_0^t k_s(O)^{-1}\,  dO_s$, the opposite inclusion also holds and we conclude that $\F^O =\F^W$.   Moreover, it is easily checked that $Z$ is
a $\P$-martingale satisfying the equation
\begin{equation}\label{statodue}
    dZ_t =Z_t
k_t(O)^{-1}h_t(M,O)\,d W_t,
\end{equation}
and that the equation \eqref{stateeqfinance}-\reff{dynM} for $(\bar X^{\alpha},M)$  can be re-written under $\P$ as
\beq\nonumber
    d\bar X_t^\alpha &=&\big[\bar b_t( \bar X^\alpha, M,O, \alpha)
    -  \bar\sigma^2_t(\bar X^\alpha, M,O, \alpha)k_t(O)^{-1}h_t(M,O)\big]\,dt
    \\\label{statotre}  && \; + \;  \bar\sigma^1_t(\bar X^\alpha,M,O, \alpha)\,d V_t
+ \bar\sigma^2_t( \bar X^\alpha,M,O,\alpha)\,dW_t, \\
dM_t &=&  \big[ \bar\beta_t(M) - \gamma_t^2(M) k_t(O)^{-1}h_t(M,O) \big] dt + \gamma_t^1(M) dV_t + \gamma_t^2(M) dW_t, \label{dynMP}
\enq
while the gain functional is re-written  as an expectation  under $\P$  from the Bayes formula:
\beq
J(\alpha)&=&\E\Big[\int_0^TZ_t \bar f_t(\bar X^\alpha,M,O, \alpha)\,dt+
Z_T\bar g(\bar X^\alpha,M,O)\Big].\label{gainrecast}
\enq

Now let us define the four-component process  $X^\alpha$ $=$ $(\bar X^\alpha,M,O,Z)$ and note that
the equations \eqref{statouno}-\eqref{statodue}-\eqref{statotre}-\eqref{dynMP}  specify a controlled stochastic equation for $X^\alpha$
of the form \eqref{stateeq} (with the obvious choice of $b$ and $\sigma$ in that equation). Similarly,  the gain functional \eqref{gainrecast}
can be put in the form  \eqref{gaineq} (with the obvious choice of $f$ and $g$).
We will see later, in Section \ref{Sec:locallylip}, how our general
results can be applied to the optimization problem
formulated in this way.


\begin{Example}\label{R:Finance}\begin{em}

As an  example of financial application, let us mention the case of a risky asset whose price $S_t$ satisfies
\beqs
dS_t &=& S_t( \rho(M_t) \,dt + \sigma_t(S)\,d\bar W_t)
\enqs
for a scalar Brownian motion $\bar W$, a volatility which is a functional   of the past values of $S$, and an unobserved return process $M$ governed by \reff{dynM}.
We assume that $\rho$, $\sigma_t(.)$ and $\sigma_t^{-1}(.)$ are bounded functions.
The wealth $\bar X_t^\alpha$ of an investor that invests a fraction $\alpha_t$ of her/his wealth in this asset (and the rest in a risk-free asset with interest rate $r$)
evolves according to the self-financing equation:
\beq
d\bar X_t^\alpha  &=&    \alpha_t\, \bar X_t^\alpha\,\frac{dS_t}{S_t} + (1-\alpha_t)\,\bar X_t^\alpha \,r\,dt  \label{autofinance} \\
&=&\bar X_t^\alpha [  r +  \alpha_t (\rho(M_t) -r) ]dt   +  \bar X_t^\alpha \alpha_t \sigma_t(S) d\bar W_t \nonumber
\enq
The investor typically observes the risky price process or equivalently the log-price process  $O_t:=\log S_t$  that solves the equation
\beqs
dO_t &=&  \big(\rho(M_t)\,-\frac{\sigma_t(S)^2}{2}\big) dt  + \sigma_t(S)\,d\bar W_t,
\enqs
which can be put in the form \eqref{observation} setting $k_t(y)=\sigma_t(\exp(y))$ and  $h_t(z,y)$ $=$ $\rho(z)-k_t(y)^2/2$.
Notice that the wealth process is  $\F^O$-adapted, since it is solution to equation \reff{autofinance}.  Therefore, when choosing the investment
strategy $\alpha$ the investor gains no additional information by observing the wealth process, and so it is reasonable to impose the condition
that $\alpha$ should be adapted to the filtration $\F^O$ alone, rather than to the one generated by $O$ and $\bar X^\alpha$.
\ep
\end{em}
\end{Example}

\subsubsection{A classical partially observed control problem}
\label{SubS:ClassicalPr}

In the previous example the observed process $O$ was not affected by the choice of the control. We next  remove this restriction, adopting
a classical approach which consists in starting with the ``reference'' probability $\P$ and introducing the  ``physical'' probability   later, as presented   e.g. in the book \cite{Be}.

Let $(\Omega,\calf, \P)$ be a complete probability space with a right-continuous and $ \P$-complete filtration $\F=(\calf_t)_{t\ge 0}$. Let $V,W$ be  independent standard Wiener
processes with respect to $\F$, with values in $\R^{m}$ and $\R^{d}$ respectively, and consider  the observation process solution to the equation in $\R^d$
\beq  \label{statounobis}
 dO_t &=& k_t(O)\,d W_t,
\enq
where $k_t(y)$  is defined on  $[0,T]\times \bfC_{d}$,  $Prog(\bfC_{d} )$-measurable,  Lipschitz  in $y$, and invertible with bounded inverse.
Similarly as in the previous paragraph, we see that $\F^W$ $=$ $\F^O$, and  an admissible control process is any $\F^W$-progressive process $\alpha$ with values in a Borel space $A$.

We are given coefficients $\bar b$, $h$, $\bar\sigma^1$, $\bar\sigma^2$  defined on $[0,T]\times\bfC_{\bar n+d}\times{\bfM_A}$, valued in
in $\R^{\bar n}$, $\R^d$, $\R^{\bar n\times m}$, $\R^{\bar n\times d}$ respectively,  and $Prog(\bfC_{\bar n+d} \times \bfM_A)$-measurable.   Then, for any admissible control process $\alpha$,
let the process $\bar X^\alpha$ be defined as the solution to the equation in $\R^{\bar n}$:
\beq
d\bar X_t^\alpha  &=&  [\bar b_t(\bar X^\alpha, O, \alpha)    - \bar\sigma^2_t(\bar X^\alpha, O, \alpha)k_t(O)^{-1}h_t(\bar X^\alpha,O,\alpha)]\,dt     \label{statotrebis} \\
& & \;\; + \; \bar\sigma^1_t(\bar X^\alpha,O,\alpha)\,d V_t \; + \; \bar\sigma^2_t(\bar X^\alpha,O,\alpha)\,dW_t. \nonumber
\enq

We introduce the gain functional $J(\alpha)$ associated to a control $\alpha$  by means of a change of probability in the following way.
 Assuming that the function $k^{-1}h$ is bounded, let us define for any admissible control process $\alpha$, the $\P$-martingale:
\beqs
Z_t^\alpha &=&  \exp\Big(  \int_0^t k_s(O)^{-1}h_s(\bar X^\alpha,O,\alpha)\,d  W_s
-\frac12 \int_0^t |k_s(O)^{-1}h_s(\bar X^\alpha,O,\alpha)|^2\,ds\Big),
\enqs
solution to the equation
\beq \label{statoduebis}
    dZ_t ^\alpha & = &  Z_t^\alpha
k_t(O)^{-1}h_t(\bar X^\alpha, O,\alpha)\,d W_t,
\enq
and introduce the ``physical'' probability $\P^\alpha$ setting   $\P^\alpha(d\omega)=Z_T^\alpha(\omega) \P(d\omega)$.
Given  real-valued coefficients $\bar f,\bar g$ defined on $[0,T]\times \bfC_{\bar n+d}\times  {\bfM_A}$ and
$\bfC_{\bar n+d}$, $Prog(\bfC_{\bar n+d} \times \bfM_A)$-measurable and  $\calc_T^{\bar n+d}$-measurable, respectively, the gain functional  is then defined as
\beqs
J(\alpha) &=&  \E^\alpha\Big[\int_0^T \bar f_t( \bar X^\alpha,O, \alpha)\,dt+  \bar g( \bar X^\alpha , O)\Big].
\enqs
The interpretation of this formulation is the following. By defining the process $W^\alpha$ as
\beqs
W_t^\alpha &=&  W_t- \int_0^s k_s(O) ^{-1} \, h_s(\bar X^\alpha,O,\alpha)\,ds, \qquad t\in [0,T],
\enqs
for any admissible control process $\alpha$, we see, by the Girsanov theorem, that the pair $(V,W^\alpha)$ is a standard
Wiener process in $\R^{m+d}$ under the probability  $\P^\alpha$ and with respect to $\F$. Moreover, the dynamics of $(\bar X^\alpha,O)$ is written under $\P^\alpha$ as:
\beqs  \label{physicalstate}
d\bar X_t^\alpha &=& \bar b_t( \bar X^\alpha, O, \alpha)\,dt+
\bar \sigma^1_t(\bar X^\alpha, O , \alpha)\,d V_t
+ \bar \sigma^2_t(\bar X^\alpha, O , \alpha)\,dW_t^\alpha,
\\\label{physicalobservation}
dO_t &=&h_t(\bar X^\alpha,O,\alpha)\,dt + k_t(O)\,d W_t^\alpha.
\enqs
We then obtain a classical  controlled state equation, and an observation process perturbed by noise
and also   affected by the choice of the control.

Finally, we notice that this problem is recast in the framework of  subsection \ref{Primal} by rewriting from Bayes formula and the $\P$-martingale property of $Z^\alpha$, the gain functional as
an expectation under $\P$:
\beq \label{gainrecastbis}
J(\alpha) &=& \E\Big[\int_0^TZ_t^\alpha \bar f_t(\bar X^\alpha,O, \alpha)\,dt+ Z_T^\alpha \bar g(\bar X^\alpha,O)\Big].
\enq
Thus, by defining the three-component process  $X^\alpha$ $=$ $(\bar X^\alpha,Z^\alpha,O)$, we see that the equations  \eqref{statounobis}-\reff{statotrebis}-\reff{statoduebis}
specify a controlled stochastic equation for $X^\alpha$  of the form \eqref{stateeq}, and the gain functional \eqref{gainrecastbis}  can be put in the form  \eqref{gaineq}.

\section{The randomized stochastic optimal control problem}
\label{Randomized}

We still assume that
$A,b,\sigma,f,g,\rho_0$
 are given and
satisfy the assumptions {\bf (A1)}.
We implement the randomization method and formulate the randomized stochastic optimal control problem associated with the control problem of subsection
\ref{Primal}. To this end we
suppose we are also given $\lambda, a_0$ satisfying
the following conditions, which are assumed to hold from now on:

\vspace{3mm}

\noindent {\bf (A2)}
\begin{itemize}
\item [(i)]  $\lambda$ is a finite positive measure on  $(A,\calb(A))$
with full topological support.
\item [(ii)] $a_0$ is a fixed, deterministic point in $A$.
\end{itemize}

\vspace{3mm}

We anticipate that  $\lambda$ will play the role of an intensity measure and
$a_0$ will be the starting point of some auxiliary
process introduced later.
Notice that the initial problem \reff{primalvalue}  does not depend on
$\lambda,a_0$, which
only appear in order to give a randomized representation of the partially observed control problem.
In this sense, {\bf (A2)} is not a restriction imposed on the
original problem
and we have the choice to fix $a_0\in A$ and
an intensity measure $\lambda$ satisfying this condition.

\subsection {Formulation of the randomized control problem}
\label{randomizedformulation}

The randomized control problem is formulated fixing a setting
$(\hat \Omega, \hat \calf,\hat \P, \hat V,\hat W, \hat \mu, \hat x_0)$, where
$(\hat \Omega, \hat \calf,\hat \P)$
is an arbitrary complete probability space
with independent random elements $\hat V$, $\hat W$, $\hat \mu$, $\hat x_0$.
The random variable
$\hat x_0$ is   $\R^n$-valued, with law $\rho_0$ under $\hat \P$.
The  process
$\hat B:=(\hat V,\hat W)$
is a standard Wiener process in $\R^{m+d}$ under $\hat \P$.
$\hat \mu$ is a Poisson random measure on  $A$ with intensity $\lambda(da)$
  under $\hat \P$;
thus, $\hat \mu$ is a sum
 of Dirac measures of the form $\hat\mu=\sum_{n\ge 1}\delta_{(\hat S_n,\hat \eta_n)}$,
 where $(\hat \eta_n)_{n\ge 1}$ is a sequence of $A$-valued random variables
 and
 $(\hat S_n)_{n\ge 1}$ is a strictly increasing sequence of random variables
with values in $(0,\infty)$, and for any $C\in\calb(A)$ the process
$\hat \mu((0,t]\times C)-t\lambda(C)$, $t\ge 0$, is a $\hat \P$-martingale.
We also define the $A$-valued process
\begin{equation}
\label{I}
\hat I_t \ = \ \sum_{n\ge 0}\hat \eta_n\,1_{[\hat S_n,\hat S_{n+1})}(t), \qquad t\ge 0,
\end{equation}
where we use the convention that $\hat S_0=0$ and $\hat I_0=a_0$, the
  point in  assumption {\bf (A2)}-(ii).
Notice that the formal sum in \eqref{I} makes sense
even if there is no addition operation defined in $A$ and that,
when $A$ is a subset of a linear space, formula \eqref{I} can be written as
\[
\hat I_t \ = \ a_0 + \int_0^t\int_A(a-\hat I_{s-})\,\hat \mu(ds\, da), \qquad t\ge 0.
\]
Let $\hat X$ be the solution to the equation
\beq \label{dynXrandom}
d\hat X_t &=&  b_t( \hat X,\hat I)\,dt + \sigma_t(\hat X,\hat I)\,dB_t,
\enq
for $t\in [0,T]$, starting from $\hat X_0$ $=$ $\hat x_0$. We define two filtrations
$\F^{\hat W,\hat \mu}=(\calf^{\hat W,\hat \mu}_t)_{t\ge 0}$
and
$\F^{\hat x_0,\hat B,\hat \mu}=(\calf^{\hat x_0,\hat B,\hat \mu}_t)_{t\ge 0}$
setting
\beq\nonumber
\calf^{\hat W,\hat \mu}_t&=&\sigma (\hat W_s,
\hat \mu((0,s]\times C)\,:\, s\in [0,t],\, C\in\calb(A))
\vee \caln,
\\
\calf^{\hat x_0,\hat B,\hat \mu}_t&=&\sigma (\hat x_0,\hat B_s,
\hat \mu((0,s]\times C)\,:\, s\in [0,t],\, C\in\calb(A))
\vee \caln, \label{expandedfiltration}
\enq
where $\caln$ denotes the family of $\hat \P$-null sets of $\hat \calf$.
We denote $\calp(\F^{\hat W,\hat \mu})$, $\calp(\F^{\hat x_0,\hat B,\hat \mu})$
 the corresponding predictable
$\sigma$-algebras.

Under {\bf (A1)} it is well-known (see e.g. Theorem 14.23 in \cite{jacod_book}) that there exists a unique $\F^{\hat x_0,\hat B,\hat \mu}$-adapted strong solution
$\hat X$ $=$ $(\hat X_t)_{0\leq t \leq T}$   to \eqref{dynXrandom}, satisfying
$\hat X_0=\hat x_0$, with continuous trajectories and such that (with the same $p$ for which
$\hat \E|\hat x_0|^p<\infty$)
\beq\label{EstimateX}
    \hat \E\,\Big[\sup_{t\in [0,T]}|\hat X_t|^p\Big]  & < & \infty.
\enq

We can now define the randomized optimal control problem as follows: the set $\hat \calv$
of admissible controls consists of all $\hat\nu=\hat \nu_t(\hat \omega,a):
\hat \Omega\times \R_+\times A\to (0,\infty)$,
which are $\calp(\F^{\hat W,\hat \mu})\otimes \calb(A)$-measurable and bounded.
Then the Dol\'eans exponential process
\beq \nonumber
\kappa_t^{\hat\nu} \ &=& \
\Ec_t\bigg(\int_0^\cdot\int_A (\hat \nu_s(a) - 1)\,(\mu(ds\, da)- \lambda(da)\,ds)\bigg) \\
&=& \ \exp\left(\int_0^t\int_A (1 - \hat \nu_s(a))\lambda(da)\,ds
\right)\prod_{0<\hat S_n\le t}\nu_{\hat S_n}(\hat \eta_n),\qquad t\ge 0,\label{doleans}
\enq
is a martingale with respect to $\hat\P$ and $\F^{\hat W,\hat \mu}$,
and we can define a new probability setting
$\hat\P^{\hat\nu}(d\hat\omega)=\kappa_T^{\hat\nu}(\hat\omega)\,\hat\P(d\hat\omega)$. From the
Girsanov theorem for multivariate point processes (\cite{ja}) it follows that under $\hat \P^{\hat\nu}$
the $\F^{\hat W,\hat \mu}$-compensator of $\hat\mu$ on the set
$[0,T]\times A$ is the random measure $\hat\nu_t(a)\lambda(da)dt$. Notice that $\hat B$ remains a Brownian motion under $\hat\P^{\hat\nu}$, and using \eqref{lipbsig}-\eqref{borbsig} we can generalize estimate \eqref{EstimateX} as follows
\beq\label{EstimateX_nu}
    \sup_{\hat\nu\in\Vc}\,\hat\E^{\hat\nu}\,\Big[\sup_{t\in [0,T]}|\hat X_t|^p\Big]  & < & \infty,
\enq
where $\hat\E^{\hat\nu}$ denotes the expectation with respect to $\hat\P^{\hat\nu}$.
We finally introduce the gain functional of the randomized control problem
\beq \label{defJrandomized}
J^\Rc(\hat\nu) &=&  \hat\E^{\hat\nu}
\left[\int_0^Tf_t(\hat X,\hat I)\,dt+g(\hat X)\right].
\enq
The randomized stochastic optimal control problem consists in maximizing
$ J^\Rc(\hat\nu)$ over all $\hat\nu\in\hat\calv$. Its value
is defined as
\begin{equation}\label{dualvalue}
{\text{\Large$\upsilon$}}_0^\Rc \;=\;   \sup_{\hat \nu\in\hat \calv} J^\Rc(\hat \nu).
\end{equation}

\begin{Remark}\label{infzero}\emph{
Let us define
$\hat\Vc_{\inf\,>\,0}=\{\hat\nu\in \hat\Vc\,:\, \inf_{\hat \Omega\times[0,T]\times A}\hat\nu>0\}$.
Then
\begin{equation}\label{eqinfzero}
{\text{\Large$\upsilon$}}_0^\Rc \;= \;
\sup_{\hat\nu\in\hat\Vc_{\inf\,>\,0}}J^\Rc(\hat\nu).
\end{equation}
Indeed, given $\hat\nu\in\hat\Vc$ and $\epsilon>0$, define $\hat\nu^\epsilon=
\hat\nu\vee \epsilon\in \hat\Vc_{\inf\,>\,0}$
and write the gain \eqref{defJrandomized} in the form
 $$
J^\Rc(\hat\nu^\epsilon) =  \hat\E
\left[\kappa_T^{\hat\nu^\epsilon}\left( \int_0^Tf_t(\hat X,\hat I)\,dt+g(\hat X)\right)\right].
 $$
 It is easy to see that $J^\Rc(\hat\nu^\epsilon) \to J^\Rc(\hat\nu) $ as $\epsilon\to 0$,
 which implies
 $${\text{\Large$\upsilon$}}_0^\Rc=\sup_{\hat\nu\in\hat\Vc}J^\Rc(\hat\nu)\le \sup_{\hat\nu\in\hat\Vc_{\inf\,>\,0}}J^\Rc(\hat\nu).
 $$
 The other inequality being obvious, we obtain  \eqref{eqinfzero}.
}
\ep
\end{Remark}

\vspace{2mm}

\begin{Remark} \label{remchoice}
{\rm
We end this section noting that a randomized control problem can be constructed starting from the initial control problem
with partial observation. Indeed, let $(\Omega,\calf,\P, \F, V,W, x_0)$
be the setting for the  stochastic optimal control problem formulated in subsection \ref{Primal}.
Suppose that  $(\Omega',\calf',\P')$ is another  probability space where a  Poisson random
 measure $\mu$ with intensity $\lambda$ is defined (for instance by a classical result,
 see \cite{Zabczyk96} Theorem 2.3.1, we may take $\Omega'=[0,1]$,
 $\calf'$ the corresponding Borel sets and $\P'$ the Lebesgue measure).
Then we define $\bar \Omega=\Omega\times \Omega'$, we denote by $\bar \calf$ the completion of
$\calf\otimes \calf'$ with respect to $\P\otimes \P'$ and by
$\bar \P$ the extension of $\P\otimes \P'$ to $\bar \calf$.
The random elements $V,W,x_0$ in $\Omega$ and the random measure
$\mu$ in $\Omega'$ have obvious extensions to $\bar \Omega$, that will be denoted by
the same symbols. Clearly, $(\bar \Omega, \bar \calf,\bar \P,   V,  W,   \mu, x_0)$
is a setting for a  randomized control problem as formulated before, that we call
\emph{product extension} of the  setting $(\Omega,\Fc,\P,V,W,x_0)$ for the initial control problem \reff{primalvalue}.

We note that the initial formulation of a randomized
setting $(\hat \Omega, \hat \calf,\hat \P, \hat V,\hat W, \hat \mu, \hat x_0)$ was more general, since  it was not required that
$\hat \Omega$ should be a product space $\Omega\times \Omega'$ and, even if it were the case, it was not required that the process $\hat B=(\hat V,\hat W)$
should depend only on $\omega\in\Omega$ while the random measure $\hat\mu$ should depend only on  $\omega'\in\Omega'$.
}
\ep
\end{Remark}

\subsection {The value of the randomized control problem}

In this section it is our purpose to show that the value ${\text{\Large$\upsilon$}}_0^\Rc$
of the randomized control problem
defined in \eqref{dualvalue}
does not depend on the specific
setting $(\hat \Omega, \hat \calf,\hat \P, \hat V,\hat W, \hat \mu, \hat x_0)$,
so that
 it is just a functional of the (deterministic) elements
$A,b,\sigma,f,g,\rho_0,\lambda, a_0$.
 Later on, in
 Theorem
\ref{MainThm},  we will prove
that in fact ${\text{\Large$\upsilon$}}_0^\Rc$
does not depend on the choice of $\lambda$ and $a_0$ either.

So let now $(\tilde \Omega, \tilde \calf,\tilde \P, \tilde V,\tilde W, \tilde \mu, \tilde x_0)$
be another setting for the randomized control problem,
as in Section \ref{randomizedformulation}, and let $\F^{\tilde W,\tilde \mu}$,
$\F^{\tilde x_0,\tilde B,\tilde \mu}$,   $\tilde X$, $\tilde I$, $\tilde \calv$ be defined in analogy
 with what was done before. So, for any
admissible control  $\tilde\nu\in \tilde \calv$, we also define $\kappa^{\tilde\nu}$ and the probability
$d\tilde\P^{\tilde\nu} =\kappa_T^{\tilde\nu} \,d\tilde\P$ as well as the gain and the value
\beqs
\tilde J^\Rc(\tilde\nu) \; = \;   \tilde\E^{\tilde\nu}
\left[\int_0^Tf_t(\tilde X,\tilde I)\,dt+g(\tilde X)\right],
\qquad
\tilde{\text{\Large$\upsilon$}}_0^\Rc \; = \;    \sup_{\tilde \nu\in\tilde \calv} \tilde J^\Rc(\tilde\nu).
\enqs
We recall that the  gain functional and value for the setting
$(\hat \Omega, \hat \calf,\hat \P, \hat V,\hat W, \hat \mu, \hat x_0)$
was defined in \reff{defJrandomized} and \eqref{dualvalue} and
denoted   by $J^\Rc$ and ${\text{\Large$\upsilon$}}_0^\Rc$ rather than $\hat J^\Rc$ and $\hat{\text{\Large$\upsilon$}}_0^\Rc$, to simplify
the notation in the following sections.

\begin{Proposition} \label{indepofthesetting}
With the previous notation, we have ${\text{\Large$\upsilon$}}_0^\Rc$ $=$ $\tilde{\text{\Large$\upsilon$}}_0^\Rc$.
In other words, ${\text{\Large$\upsilon$}}_0^\Rc$ only depends on the objects $A,b,\sigma,f,g,\rho_0,\lambda, a_0$ appearing
in the assumptions {\bf (A1)} and {\bf (A2)}.
\end{Proposition}

 \noindent {\bf Proof.} It is enough to prove that ${\text{\Large$\upsilon$}}_0^\Rc\le \tilde{\text{\Large$\upsilon$}}_0^\Rc$,
 since the opposite inequality is established by the same arguments. Writing the gain $J^\Rc(\hat\nu)$ defined in
 \eqref{defJrandomized}
 in the form
 $$
J^\Rc(\hat\nu) =  \hat\E
\left[\kappa_T^{\hat\nu}\left( \int_0^Tf_t(\hat X,\hat I)\,dt+g(\hat X)\right)\right],
 $$
 recalling the definition \eqref{doleans} of the process $\kappa^{\hat\nu}$
 and noting that the
process $\hat I$ is completely determined by $\hat\mu$,
we see that  $J^\Rc(\hat\nu)$ only depends on the
(joint) law of $(\hat X,\hat \mu,\hat \nu)$ under $\hat\P$.
Since, however,  $\hat X$ is the solution to equation \eqref{dynXrandom} with initial
condition  $\hat X_0$ $=$ $\hat x_0$,
 it is easy to check that under  our assumptions
the  law of $(\hat X,\hat \mu,\hat \nu)$  only depends on the law of
$(\hat x_0, \hat V,\hat W,\hat \mu, \hat\nu)$.
Since $\hat x_0$, $\hat V$ and $(\hat W,\hat \mu, \hat\nu)$
are all independent, and the laws of
$\hat x_0$ and $\hat V$ are
  fixed (since
$\hat V$ is a standard Wiener process and $ \hat x_0$
has law $\rho_0$) we conclude that $J^\Rc(\hat\nu)$ only depends on the law of
$(\hat W,\hat \mu, \hat\nu)$ under $\hat\P$.  Similarly, $\tilde J^\Rc(\tilde\nu)$ only depends on the law of $(\tilde W,\tilde \mu,\tilde \nu)$
 under $\tilde\P$.

 Next we claim that, given $\hat\nu\in\hat \calv$ there exists
 $\tilde\nu\in\tilde \calv$ such that the law of
 $(\hat W,\hat \mu,\hat \nu)$ under $\hat \P$ is the same as
 the law of
 $(\tilde W,\tilde\mu,\tilde \nu)$ under $\tilde\P$.
 Assuming the claim for a moment, it follows from the previous
 discussion that for this choice of $\tilde\nu$ we have
 \beqs
 J^\Rc(\hat\nu) \; = \; \tilde J^\Rc(\tilde\nu) &\le&  \tilde{\text{\Large$\upsilon$}}_0^\Rc,
 \enqs
 and taking the supremum over $\hat\nu\in\hat \calv$ we deduce that ${\text{\Large$\upsilon$}}_0^\Rc$ $\le$ $\tilde{\text{\Large$\upsilon$}}_0^\Rc$, which proves the result.

It only remains to prove the claim. By a monotone class argument
we may suppose that $\hat\nu_t(a)=k(a)\,\phi_t\,\psi_t$,
where $k$ is a $\calb(A)$-measurable, $\phi$ is $\F^{\hat W}$-predictable
and $\psi$ is $\F^{\hat \mu}$-predictable (where these filtrations
are the ones generated by $\hat W$ and $\hat \mu$ respectively).
We may further suppose that
$$
\phi_t=1_{(t_0,t_1]}(t)\phi_0(\hat W_{s_1},\ldots,\hat W_{s_h})
$$
for an integer $h$ and deterministic times $0\le s_1\le\ldots s_h\le t_0 <t_1$
and a Borel function $\phi_0$ on $\R^h$,
since this class of processes generates the predictable $\sigma$-algebra
of $\F^{\hat W}$, and that
$$
\psi_t=1_{(\hat S_n,\hat S_{n+1}]}(t)\psi_0(\hat S_1,\ldots,
\hat S_n, \hat \eta_1,\ldots, \hat \eta_n,t)
$$
for an integer $n\ge 1$ and
 a Borel function $\psi_0$ on $\R^{2n+1}$,
since this class of processes generates the predictable $\sigma$-algebra
of $\F^{\hat \mu}$ (see \cite{ja}, Lemma (3.3)). It is immediate to verify that
the required process $\tilde\nu$ can be defined setting
$$
\tilde\nu_t(a)=k(a)\,
1_{(t_0,t_1]}(t)\phi_0(\tilde W_{s_1},\ldots, \tilde W_{s_h}) \,
1_{(\tilde S_n,\tilde S_{n+1}]}(t)\psi_0(\tilde S_1,\ldots,
\tilde S_n, \tilde \eta_1,\ldots, \tilde \eta_n,t),
$$
where $(\tilde S_n,\tilde \eta_n)_{n\ge 1}$
are associated to the measure $\tilde \mu$, i.e.
$\tilde\mu=\sum_{n\ge 1}\delta_{(\tilde S_n,\tilde \eta_n)}$.
 \qed

\subsection{Equivalence of the partially observed and the randomized control pro\-blem}
We can now state one of the main results of the paper.

\begin{Theorem}
\label{MainThm}
Assume that {\bf (A1)} and {\bf (A2)} are satisfied.
Then the values of the partially observed control problem
and of the randomized control problem are equal:
\beq \label{equivrandom}
{\text{\Large$\upsilon$}}_0 &=& {\text{\Large$\upsilon$}}_0^\Rc,
\enq
where ${\text{\Large$\upsilon$}}_0$ and ${\text{\Large$\upsilon$}}_0^\Rc$ are defined by  \eqref{primalvalue} and \eqref{dualvalue} respectively.
This common value only depends on the objects $A,b,\sigma,f,g,\rho_0$ appearing in assumption {\bf (A1)}.
\end{Theorem}

The last sentence follows immediately from
Proposition
\ref{indepofthesetting}, from the equality ${\text{\Large$\upsilon$}}_0$ $=$  ${\text{\Large$\upsilon$}}_0^\Rc$
and from the obvious fact that ${\text{\Large$\upsilon$}}_0$ cannot depend on $\lambda, a_0$ introduced in assumption {\bf (A2)}.
The proof of the equality is contained in the next section.

Before giving the proof of Theorem \ref{MainThm}, let us discuss the significance of this equivalence result. The randomized control problem involves an uncontrolled state process $(X,I)$ solution to  \reff{I}-\reff{dynXrandom}, and the optimization is done over a set of equivalent probability measures whose effect is to change the characteristics (the intensity) of the auxiliary randomized process $I$
without impacting on the Brownian motion $B$ driving $X$. Therefore,  the equivalence result  \reff{equivrandom} means that by performing such optimization in the randomized problem, we
 achieve the same value as in the original control problem where controls affect directly the drift and diffusion of the state process.  As explained in the Introduction,
such equivalence result has important implications that will be addressed in Section \ref{Sec:separandom} where it is shown that the randomized control problem is associated by duality to a backward stochastic differential equation (with nonpositive jumps), called  the
{\it randomized  equation}, which then also
characterizes  the value function of the initial control problem \reff{primalvalue}.

\begin{Remark} \label{comparison_with_FP15}
{\rm
We mention that in the article \cite{FP15}   an equivalence result
similar to Theorem \ref{MainThm} was proved.
However, in \cite{FP15} only the case of full observation was addressed and there was no memory effect with respect to the control, whereas path-dependence in the state variable was allowed.
But
the main difference with respect to our setting
 is that in  \cite{FP15} the primal problem was formulated in a weak form,
 i.e. taking the supremum of the gain functional
 \eqref{gainvalueintro}
 also over all possible choices of the probability space
 $(\Omega,\calf,\P)$.
This simplifies many arguments, and in particular makes the inequality
${\text{\Large$\upsilon$}}_0\ge {\text{\Large$\upsilon$}}_0^\Rc$
 trivial. Finally,  in Theorem
\ref{MainThm_bis} below we will present an extension to the case
of locally Lipschitz coefficients with linear growth,
a more general situation that was not
addressed  in \cite{FP15}.
}
\ep
\end{Remark}

\section{Proof of Theorem \ref{MainThm}}
\label{proofthm}

The proof  is split into two parts, corresponding to the inequalities ${\text{\Large$\upsilon$}}_0^\Rc$ $\le$ ${\text{\Large$\upsilon$}}_0$ and ${\text{\Large$\upsilon$}}_0$ $\le$ ${\text{\Large$\upsilon$}}_0^\Rc$.
In the sequel, {\bf (A1)} and {\bf (A2)} are always assumed to hold.
However, instead of the inequality $p\ge \max(2,2r)$, in {\bf (A1)}-(v)
 it is enough to suppose that  $p\ge \max(2,r)$.

Before starting with the rigorous proof, let us have a look at the main points.
\begin{itemize}
\item ${\text{\Large$\upsilon$}}_0^\Rc$ $\le$ ${\text{\Large$\upsilon$}}_0$. First, we prove that the value of the primal problem ${\text{\Large$\upsilon$}}_0$ does not change if we reformulate it on the enlarged probability space where the randomized problem lives, taking the supremum over $\Ac^{W,\mu'}$, which is the set of controls $\bar\alpha$ progressively measurable with respect to the filtration generated by $W$ and the Poisson random measure $\mu'$ (Lemma \ref{L:AWmu}; actually, we take $\bar\alpha$ progressively measurable with respect to an even larger filtration, denoted $\F^{W,\mu'_\infty}$). Second, we prove that for every $\nu\in\Vc_{\inf>0}$ there exists $\bar\alpha^\nu\in\Ac^{W,\mu'}$ such that $\mathscr L_{\P^\nu}(x_0,B,I)=\mathscr L_{\bar\P}(x_0,B,\bar\alpha^\nu)$ (Proposition \ref{P:Ineq_I}). This result is a direct consequence of the key Lemma \ref{L:Replication_I}. From $\mathscr L_{\P^\nu}(x_0,B,I)=\mathscr L_{\bar\P}(x_0,B,\bar\alpha^\nu)$ we obtain that $J^\Rc(\nu)=\bar J(\bar\alpha^\nu)$, namely
\[
{\text{\Large$\upsilon$}}_0^\Rc \ := \ \sup_{\nu\in\calv_{\inf>0}}J^\Rc(\nu) \ = \ \sup_{\substack{\bar\alpha^\nu\\ \nu\in\Vc_{\inf>0}}}\bar J(\bar\alpha^\nu).
\]
Since ${\text{\Large$\upsilon$}}_0$ $=$ $\sup_{\bar\alpha\in\Ac^{W,\mu'}} \bar J(\bar\alpha)$ by Lemma \ref{L:AWmu}, and every $\bar\alpha^\nu$ belongs to $\Ac^{W,\mu'}$, we easily obtain the inequality ${\text{\Large$\upsilon$}}_0^\Rc$ $\le$ ${\text{\Large$\upsilon$}}_0$.
\item ${\text{\Large$\upsilon$}}_0$ $\le$ ${\text{\Large$\upsilon$}}_0^\Rc$. The proof of this inequality is based on a ``density'' result (which corresponds to the key Proposition \ref{extensionapproximation}) in the spirit of Lemma 3.2.6 in \cite{80Krylov}. Roughly speaking, we prove that the class $\{\bar\alpha^\nu\colon\nu\in\Vc_{\inf>0}\}$ is dense in $\Ac^{W,\mu'}$, with respect to the metric $\tilde\rho$ defined in \eqref{MetricKrylov} (the same metric used in Lemma 3.2.6 in \cite{80Krylov}). Then, the inequality ${\text{\Large$\upsilon$}}_0$ $\le$ ${\text{\Large$\upsilon$}}_0^\Rc$ follows from the stability Lemma \ref{contrhotilde}, which states that, under Assumption {\bf (A1)}, the gain functional is continuous with respect to the metric $\tilde\rho$.
\end{itemize}

\subsection{Proof of the inequality ${\text{\Large$\upsilon$}}_0^\Rc$ $\le$ ${\text{\Large$\upsilon$}}_0$}

We note at the outset that the requirement that $\lambda$ has full support will not be used in the proof of the inequality ${\text{\Large$\upsilon$}}_0^\Rc$ $\le$
${\text{\Large$\upsilon$}}_0$.

 Let $(\Omega,\calf,\P, \F, V,W, x_0)$
be a setting for the stochastic optimal control problem with partial observation formulated in subsection \ref{Primal}. We construct a setting for a randomized control problem in the form of a product
 extension as described at the end of Section \ref{randomizedformulation}.

Let $\lambda$ be a Borel measure on $A$ satisfying {\bf (A2)}. As a first step, we need to
construct a suitable surjective measurable map $\pi:\R\to A$ and to
introduce a properly chosen measure $\lambda'$ on the Borel subsets of the real line
such that in particular $\lambda$ $=$ $\lambda' \circ \pi^{-1}$.
We also recall that the space of control actions $A$ is assumed to be a Borel space and it is  known  that  any such  space is either finite or countable (with the discrete topology) or isomorphic, as a measurable space, to the real line (or equivalently
to the half line $(0,\infty)$): see e.g. \cite{BertsekasShreve78}, Corollary 7.16.1.

Let us denote by $A_c$ the subset of $A$ consisting of all points
$a\in A$ such that $\lambda(\{a\})>0$, and let
$A_{nc}=A\backslash A_c$. Since $\lambda$ is finite, the set
$A_c$ is either empty or countable, and it follows in particular
that both $A_c$ and $A_{nc}$ are also Borel spaces.

 In the construction of
$\lambda'$ we distinguish the following three cases.

\begin{enumerate}
  \item $A_c=\emptyset$, so that $A=A_{nc}$ is uncountable. Then,
  as recalled above, there exists a
  bijection $\pi:\R\to A$ such that $\pi$ and its inverse are both
  Borel measurable. We define a measure $\lambda'$ on  $(\R,\calb(\R))$   setting $\lambda'(B)=\lambda (\pi(B))$ for $B\in\calb(\R)$.
Even if we cannot guarantee that $\lambda'$ has full support, it clearly holds that  $\lambda'(\{r\})=0$ for every $r\in\R$.
Basically, in this case we are identifying $A$ with $\R$ and $\lambda$ with its image measure $\lambda'$.
\item $A_{nc}=\emptyset$, so that $A=A_c$ is  countable, with the discrete topology. For every
$j\in A$ choose a (nontrivial) interval   $\Ic_j\subset \R$ in such a way that $\{\Ic_j,\,j\in A\}$ is a partition of $\R$.
Choose an arbitrary nonatomic finite measure on $(\R,\calb(\R))$  with full support  (say, the standard Gaussian measure, denoted by $\gamma$)
 and denote by  $\lambda'$  the unique positive measure  on $(\R,\Bc(\R))$ such that
\[
\lambda'(B) \ = \ \lambda(\{j\})\gamma(B)/\gamma(\Ic_j),
\qquad \text{for every }B\subset \Ic_j,\,  B\in\calb(\R),\,j\in A.
\]
Notice that $\lambda'$ is a finite measure ($\lambda'(\R)=\lambda(A)$),  satisfying $\lambda'(\Ic_j)=\lambda(\{j\})$ for every $j\in A$
and $\lambda'(\{r\})=0$ for every $r\in\R$.
We also define the projection $\pi\colon\R\rightarrow A$ given by
\begin{equation}\label{defproiez}
    \pi(r)=j,
\qquad \text{if }
r\in \Ic_j  \text{ for some } j\in A.
\end{equation}
Clearly, $\lambda$ $=$ $\lambda' \circ \pi^{-1}$.
\item $A_{c}\neq\emptyset$ and $A_{nc}\neq\emptyset$.
For every
$j\in A_c$ choose a (nontrivial) interval   $\Ic_j\subset (-\infty,0]$
in such a way that $\{\Ic_j,\,j\in A_c\}$ is a partition of $(-\infty,0]$.
Moreover, there exists a
  bijection $\pi_1:(0,\infty)\to A_{nc}$ such that $\pi_1$ and its inverse are both
  Borel measurable.
Denote by  $\lambda'$  the unique positive measure  on $(\R,\Bc(\R))$ such that
\begin{eqnarray*}
\lambda'(B) \ = \ \lambda(\{j\})\gamma(B)/\gamma(\Ic_j),
& \text{for every }B\subset \Ic_j,\,  B\in\calb(\R),\,j\in A_c,
\nonumber\\
\lambda'(B)=\lambda (\pi_1(B))&
 \text{for every }B\subset (0,\infty),\, B\in\calb(\R).
\end{eqnarray*}
Again, $\lambda'$ is a finite measure
 satisfying $\lambda'(\Ic_j)=\lambda(\{j\})$ for every $j\in A_c$
and $\lambda'(\{r\})=0$ for every $r\in\R$.
We also define the projection $\pi\colon\R\rightarrow A$ given by
\begin{equation}\label{defproiezdue}
\pi(r) \ = \
\begin{cases}
j, \qquad & \text{if }
r\in \Ic_j  \text{ for some } j\in A_c, \\
\pi_1(r),
& \text{if }
r\in (0,\infty),
\end{cases}
\end{equation}
so that in particular $\lambda$ $=$ $\lambda' \circ \pi^{-1}$.
\end{enumerate}

Now  let  $(\Omega',\calf',\P')$ denote the canonical probability space of a
 non-explosive Poisson point process on $\R_+\times \R$ with intensity $\lambda'$.
 Thus, $\Omega'$ is the set of sequences $\omega'=(t_n,r_n)_{n\geq1}\subset(0,\infty)\times \R$ with $t_n<t_{n+1}\nearrow\infty$,
 $(T_n,R_n)_{n\geq1}$ is the canonical marked point process  (i.e.  $T_n(\omega')=t_n$, $R_n(\omega')=r_n$),  and
 $\mu'$ $=$ $\sum_{n\ge 1}\delta_{(T_n,R_n)}$ is the corresponding random measure.
 Let $\Fc'$ denote the   smallest $\sigma$-algebra such that all the maps $T_n,R_n$ are measurable, and   $\P'$  the unique probability on
$\Fc'$ such that  $\mu'$ is a Poisson random measure with intensity $\lambda'$ (since $\lambda'$ is a finite measure,  this probability actually exists).
 We will also use the completion of the space  $(\Omega',\calf',\P')$, still denoted by the same symbol by abuse of notation. In all the cases considered
above,  setting
$$
A_n=\pi(R_n), \qquad
\mu=\sum_{n\ge 1}\delta_{(T_n,A_n)},
$$
it is easy to verify that $\mu$ is a Poisson random measure on $(0,\infty)\times A$ with intensity $\lambda$, defined in $(\Omega',\calf',\P')$. Then, following \reff{I}, we associate to this Poisson random measure on $(0,\infty)\times A$, the $A$-valued process
\beqs
I_t &=& \sum_{n\ge 0}  A_n\,1_{[T_n,T_{n+1})}(t), \qquad t\ge 0,
\enqs
where we use the convention that $T_0=0$ and $I_0=a_0$ the  point in  assumption {\bf (A2)}-(ii).
In $(\Omega',\Fc')$ we define the natural filtrations $\F^\mu=(\calf^\mu_t)_{t\ge 0}$,
$\F^{\mu'}=(\calf^{\mu'}_t)_{t\ge 0}$
given by
\begin{eqnarray}
\calf^\mu_t&=&\sigma\big(\mu((0,s]\times C)\,:\, s\in [0,t],\, C\in\calb(A)\big)
\vee \caln',\nonumber\\
\calf^{\mu'}_t&=&\sigma\big(\mu'((0,s]\times B)\,:\, s\in [0,t],\, B\in\calb(\R)\big)
\vee \caln',\nonumber
\end{eqnarray}
where $\caln'$ denotes the family of $\P'$-null sets of $\Fc'$. We denote by $\Pc(\F^\mu)$,
 $\Pc(\F^{\mu'})$the corresponding predictable $\sigma$-algebras.
 Note that $\calf^{\mu}_t\subset \calf^{\mu'}_t$ and
$ \calf^{\mu'}_\infty=\calf'$.

Then we define $\bar \Omega=\Omega\times \Omega'$, we denote by $\bar \calf$ the completion of
$\calf\otimes \calf'$ with respect to
$\P\otimes \P'$ and by
$\bar \P$ the extension of $\P\otimes \P'$ to $\bar \calf$.
The random elements $V,W,x_0$ in $\Omega$ and the random measures
$\mu,\mu'$ in $\Omega'$ have obvious extensions to $\bar \Omega$, that will be denoted by
the same symbols. Then $(\bar \Omega, \bar \calf,\bar \P,   V,  W,   \mu,  x_0)$
is a setting for a  randomized control problem as formulated  in section \ref{randomizedformulation}.
Recall that $\F^{ W}$ denotes the $\P$-completed filtration in $( \Omega,   \calf)$ generated by the Wiener
process $W$. All filtrations $\F^{ W}$,  $\F^\mu$,   $\F^{\mu'}$  can also be lifted to  filtrations in $(\bar \Omega, \bar \calf)$, and $\bar\P$-completed.
In the sequel it should be clear from the context whether they are considered as filtrations in $(\bar \Omega, \bar \calf)$ or in their original spaces.
As in Section \ref{randomizedformulation} we   define the filtration
$\F^{ W,  \mu}=(\calf^{ W,  \mu}_t)_{t\ge 0}$ in $(\bar \Omega, \bar \calf)$  by
$$
\calf^{  W,  \mu}_t  =
\Fc_t^W\vee\calf^\mu_t\vee\caln,
$$
 ($\caln$ denotes the family of $\bar \P$-null sets of $\bar \calf$),   we   introduce the classes $  \calv, \Vc_{\inf\,>\,0}$ and, for any
 admissible control $ \nu \in \calv$, the corresponding martingale $\kappa^{\nu}$, the probability
 $\P^{ \nu}(d\omega\, d\omega' )=\kappa_T^{\nu}(\omega,\omega')\,\bar\P(d \omega\,d\omega')$ and the gain $J^\Rc(\nu)$.
 For technical purposes, we need to introduce the set  $\calv'$ of  elements  $\nu'$ $=$
 $\nu'_t(\omega',a): \Omega'\times \R_+\times A\to (0,\infty)$, which are $\calp(\F^{\mu})\otimes \calb(A)$-measurable and bounded.
We also define  another filtration
$\F^{W,\mu'_\infty}=(\calf^{W,\mu'_\infty}_t)_{t\ge 0}$
in $(\bar \Omega, \bar \calf)$ setting
$$
\calf^{W,\mu'_\infty}_t \ = \ \Fc_t^W\vee\Fc'\vee\caln
$$
(here $\Fc'$ denotes a $\sigma$-algebra in $(\bar \Omega, \bar \calf)$, namely
$\{\Omega\times B\, :\, B\in \Fc'\}$).

\vspace{1mm}

In order to prove the inequality ${\text{\Large$\upsilon$}}_0^\Rc$ $\le$ ${\text{\Large$\upsilon$}}_0$,  we first prove two technical lemmata. In particular, in Lemma \ref{L:AWmu}  we show that the primal problem is equivalent to a new primal problem with $\F^{W,\bar\mu_\infty}$-progressive controls on the enlarged space
$(\bar\Omega,\bar\Fc)$.

\begin{Lemma}\label{L:AWmu}
We have ${\text{\Large$\upsilon$}}_0$ $=$ $\sup_{\bar\alpha\in\Ac^{W,\mu'}} \bar J(\bar\alpha)$, where
\[
\bar J(\bar\alpha) \ = \ \bar\E\left[\int_0^Tf_t(X^{\bar\alpha},\bar\alpha)\,dt+g(X^{\bar\alpha})\right],
\]
and $\Ac^{W,\mu'}$ is the set of all $\F^{W,\mu'_\infty}$-progressive processes $\bar\alpha$ with values in $A$. Moreover, $X^{\bar\alpha}$ $=$ $(X_t^{\bar\alpha})_{0\leq t\leq T}$ is the strong solution to \eqref{stateeq} $($with $\bar\alpha$ in place of $\alpha$$)$
satisfying $X_0^{\bar\alpha}=x_0$, which is unique in the class of continuous processes adapted to the filtration
$(\calf^{B}_t\vee \sigma (x_0)\vee \calf'\vee \caln)_{t\ge 0}$.
\end{Lemma}
\textbf{Proof.}
The   inequality ${\text{\Large$\upsilon$}}_0$ $\le$ $\sup_{\bar\alpha\in\Ac^{W,\mu'}} \bar J(\bar\alpha)$ is immediate, since every control
$\alpha\in\Ac^W$ also lies in $\Ac^{W,\mu'}$ and $J(\alpha)=\bar J(\alpha)$, whence $J(\alpha)\leq\sup_{\bar\alpha\in\Ac^{W,\mu'}}
\bar J(\bar\alpha)$ and so ${\text{\Large$\upsilon$}}_0$ $=$  $\sup_{\alpha\in\Ac^W} J(\alpha)\leq\sup_{\bar\alpha\in\Ac^{W,\mu'}} \bar J(\bar\alpha)$.

Let us prove the opposite inequality. Fix $\tilde\alpha\in\Ac^{W,\mu'}$ and consider the (uncompleted) filtration
$\F'':=( \Fc_t^W\vee\Fc')_{t\ge0}$. Then we can find  an  $A$-valued $\F''$-progressive process $\bar   \alpha$
such that  $\bar\alpha=\tilde\alpha$   $\bar\P(d\bar\omega)dt$-almost surely, so that  in particular $\bar J(\bar\alpha)=\bar J(\tilde\alpha)$.
It is easy to verify that, for every $\omega'\in\Omega'$,  the process $\alpha^{\omega'}$, defined by
$\alpha_t^{\omega'}(\omega):=\bar\alpha_t(\omega,\omega')$, is $\F^W$-progressive. Consider now the controlled equation on $[0,T]$
\begin{align}
\label{SDEomega'}
X_t \ &= \ x_0 + \int_0^t b_s(X, \alpha^{\omega'})\,ds + \int_0^t \sigma_s(X, \alpha^{\omega'})\,dB_s \\
&= \ x_0 + \int_0^t b_s(X, \bar\alpha(\cdot,\omega'))\,ds + \int_0^t \sigma_s(X, \bar\alpha(\cdot,\omega'))\,dB_s. \notag
\end{align}
From the first line of \eqref{SDEomega'} we see that, under Assumption {\bf (A1)}, for every $\omega'$ there exists a unique (up to indistinguishability) continuous  process $X^{\alpha^{\omega'}}=(X_t^{\alpha^{\omega'}})_{0\leq t\leq T}$ strong solution to \eqref{SDEomega'}, adapted to the filtration
$(\calf^{B}_t\vee \sigma (x_0)\vee \caln)_{t\ge 0}$.
 On the other hand, from the second line of \eqref{SDEomega'}, it follows that the process $X^{\bar\alpha}(\cdot,\omega')=(X_t^{\bar\alpha}(\cdot,\omega'))_{0\leq t\leq T}$ solves the above equation. From the pathwise uniqueness of strong solutions to equation \eqref{SDEomega'}, it follows that $X_t^{\alpha^{\omega'}}(\omega)=X_t^{\bar\alpha}(\omega,\omega')$, for all $t\in[0,T]$, $\P(d\omega)$-a.s. By the Fubini theorem
\[
\bar J(\tilde\alpha) =
\bar J(\bar\alpha)  = \int_{\Omega'} \E\bigg[\int_0^Tf_t(X^{\alpha^{\omega'}},\alpha^{\omega'})\,dt+g(X^{\alpha^{\omega'}})
\bigg] \, \P'(d\omega').
\]
Since the inner expectation equals the gain $J(\alpha^{\omega'})$,
it cannot exceed $V$ and it follows that
$\bar J(\tilde\alpha)\le {\text{\Large$\upsilon$}}_0$.  The claim follows from the  arbitrariness of $\tilde\alpha$.
\ep

\vspace{3mm}

The next result provides a decomposition of any element $\nu$ $\in$ $\Vc$, i.e.  $\calp(\F^{W,\mu})\otimes \calb(A)$-measurable and bounded.

\begin{Lemma}\label{L:nu_pred}
\textup{(i)} Let $\nu\in\Vc$, then there exists a $\bar\P$-null set $\bar N\in\Nc$ such that $\nu$ admits the following representation
\begin{align*}
\nu_t(\omega,\omega',a) \ &= \ \nu_t^{(0)}\big(\omega,a\big) \, 1_{\{0<t\leq T_1(\omega')\}} \\
&\quad \ + \sum_{n=1}^\infty \nu_t^{(n)}\big(\omega,(T_1(\omega'),A_1(\omega')),\ldots,(T_n(\omega'),A_n(\omega')),a\big) \, 1_{\{T_n(\omega')<t\leq T_{n+1}(\omega')\}}, \notag
\end{align*}
for all $(\omega,\omega',t,A)\in\bar\Omega\times\R_+\times A$, $(\omega,\omega')\notin\bar N$, for some maps $\nu^{(n)}\colon\Omega\times\R_+\times(\R_+\times A)^n\times A\rightarrow(0,\infty)$, $n\geq1$, $($resp. $\nu^{(0)}\colon\Omega\times\R_+\times A\rightarrow(0,\infty)$$)$, which are $\Pc(\F^W)\otimes\Bc((\R_+\times A)^n)\otimes\Bc(A)$-measurable $($resp. $\Pc(\F^W)\otimes\Bc(A)$-measurable$)$ and uniformly bounded with respect to $n$. Moreover, if $\nu\in\Vc_{\inf\,>\,0}$  then $\inf_{\bar \Omega\times [0,T]\times A}\nu^{(n)}>0$
as well, for every $n\geq0$.

\vspace{2mm}

\noindent\textup{(ii)} Let $\nu\in\Vc$, then there exists $\tilde N\in\Fc$, with $\P(\tilde N)=0$, such that the map $\nu^\omega=\nu_t^\omega(\omega',a): \Omega'\times \R_+\times A\to (0,\infty)$, defined by
\[
\nu_t^\omega(\omega',a) \ := \ \nu_t(\omega,\omega',a), \qquad (\omega',t,a)\in\Omega'\times \R_+\times A,
\]
belongs to $\Vc'$ whenever $\omega\notin\tilde N$. Moreover, for every $\omega\notin\tilde N$ there exists $N_\omega\in\Nc'$ such that
\begin{align}\label{representation_nu^omega}
\nu_t^\omega(\omega',a) \ &= \ \nu_t^{(0)}\big(\omega,a\big) \, 1_{\{0<t\leq T_1(\omega')\}} \\
&\quad \ + \sum_{n=1}^\infty \nu_t^{(n)}\big(\omega,(T_1(\omega'),A_1(\omega')),\ldots,(T_n(\omega'),A_n(\omega')),a\big) \, 1_{\{T_n(\omega')<t\leq T_{n+1}(\omega')\}}, \notag
\end{align}
for all $(\omega',t,A)\in\Omega'\times\R_+\times A$, $\omega'\notin N_\omega'$,
where, clearly, $\nu_\cdot^{(n)}(\omega,\cdot)$ $($resp. $\nu_\cdot^{(0)}(\omega,\cdot)$$)$ is $\Bc(\R_+)\otimes\Bc((\R_+\times A)^n)\otimes\Bc(A)$-measurable $($resp. $\Bc(\R_+)\otimes\Bc(A)$-measurable$)$.
\end{Lemma}

\noindent {\bf Proof.}
The proof is an extension of the results in
 \cite{ja} Lemma 3.3, it is based on
monotone class arguments and is left to the reader.
\ep

\vspace{3mm}

By Lemma \ref{L:nu_pred}-(ii), given $\nu\in\Vc$, consider the process $\nu^\omega$ $\in$ $\Vc'$,
with corresponding $\P$-null set $\tilde N\in\Fc$. Define the Dol\'eans exponential process $\kappa^{\nu^\omega}$ by formula \eqref{doleans} with $\nu^\omega$ in place of $\nu$. Notice that by Lemma \ref{L:nu_pred}-(ii) we have $\kappa_t^{\nu^\omega}(\omega')=\kappa_t^\nu(\omega,\omega')$, for all $(\omega',t)\in\Omega'\times\R_+$, whenever $\omega\notin\tilde N$. Moreover, for $\omega\notin\tilde N$, $(\kappa_t^{\nu^\omega})_{t\ge 0}$ is a martingale
with respect to $\P'$ and $\F^\mu$.  We claim that there exists
 a unique probability measure $\P^{\nu^\omega}$ on $(\Omega',\calf^\mu_\infty)$ such that $\P^{\nu^\omega}(d\omega')=\kappa_t^{\nu^\omega}(\omega')\P'(d\omega')$
 on each $\sigma$-algebra $\calf^\mu_t$ and,
by the Girsanov theorem,  the $\F^\mu$-compensator of $\mu$ under $\P^{\nu^\omega}$
 is given by the right-hand side of \eqref{representation_nu^omega}.

The verification  of the claim  is a standard argument: using
the boundedness of $\nu$  one first verifies that
$$
\kappa^{\nu^\omega}_{t\wedge T_n}(\omega')\le a_n\,e^{b\,T_n(\omega')}
$$
for some constants $a_n,b$,
which implies that $(\kappa^\nu_{t\wedge T_n})_{t\ge 0}$ is a
uniformly integrable martingale with respect to $\P'$ and $\F^\mu$.
Then the probabilities $\P_n^{\nu^\omega}$ defined on $\calf_{T_n}^\mu$ setting
$\P_n^{\nu^\omega}(d\omega')=\kappa_{T_n}^{\nu^\omega}(\omega')\,\P'(d\omega')$
satisfy the  compatibility condition: $\P_{n+1}^{\nu^\omega}=\P_n^{\nu^\omega}$ on $\calf_{T_n}^\mu$
for every $n$. Arguing as in Theorem  3.6 in \cite{ja}, by  the Kolmogorov extension theorem
there exists a unique probability $\P^{\nu^\omega}$ on $(\Omega',\calf^\mu_\infty)$
such that $\P^{\nu^\omega}=\P^{\nu^\omega}_n$ on each $\calf^\mu_{T_n}$,
and $\P^{\nu^\omega}$ has the required properties.


We can now state the following key  result (Lemma \ref{L:Replication_I}) from  which the required conclusion  of this subsection
follows readily (see Proposition \ref{P:Ineq_I}). Recall that $(\Omega',\Fc',\P')$ denotes the canonical probability space constructed above.

\begin{Lemma}
\label{L:Replication_I} Given $\nu\in\Vc_{\inf\,>\,0}$, there exist a  sequence $(T_n^\nu,A_n^\nu)_{n\geq1}$ on $(\bar\Omega,\bar\Fc,\bar\P)$ and a $\P$-null set $N\in\Fc$, with $\tilde N\subset N$ ($\tilde N$ is the set appearing in Lemma \ref{L:nu_pred}-(ii)), such that:
\begin{itemize}
\item[\textup{(i)}] for every $n\geq1$, $(T_n^\nu,A_n^\nu)$ takes values in $(0,\infty)\times A$ and $T_n^\nu<T_{n+1}^\nu$;
\item[\textup{(ii)}] for every $n\geq1$, $T_n^\nu$ is an $\F^{W,\mu'_\infty}$-stopping time and $A_n^\nu$ is $\Fc_{T_n^\nu}^{W,\mu'_\infty}$-measurable;
\item[\textup{(iii)}] $\lim_{n\rightarrow\infty}T_n^\nu=\infty$;
\item[\textup{(iv)}] for every $\omega\notin N$, we have
\[
\mathscr L_{\P'}\big((T_n^\nu(\omega,\cdot),A_n^\nu(\omega,\cdot))_{n\geq1}\big) \ = \ \mathscr L_{\P^{\nu^\omega}}\big((T_n,A_n)_{n\geq1}\big).
\]
\end{itemize}
Finally, let $\bar\alpha_t^\nu=a_01_{[0,T_1^\nu)}+\sum_{n=1}^\infty A_n^\nu 1_{[T_n^\nu,T_{n+1}^\nu)}(t)$ be the step process associated with $(T_n^\nu,A_n^\nu)_{n\geq1}$. Then, $\bar\alpha^\nu\in\Ac^{W,\mu'}$ and $\mathscr L_{\P'}(\bar\alpha^\nu(\omega,\cdot))=\mathscr L_{\P^{\nu^\omega}}(I)$, $\omega\notin N$.
\end{Lemma}
\textbf{Proof.}
Suppose that we have already constructed a multivariate point process $(T_n^\nu,A_n^\nu)_{n\geq1}$ satisfying points (i)-(ii)-(iii)-(iv) of the Theorem. Then, by (ii) it follows that $\bar\alpha^\nu$ is c\`adl\`ag and $\F^{W,\mu'_\infty}$-adapted, hence progressive. Moreover, by (iii), for every $(\bar\omega,t)\in\bar\Omega\times[0,T]$ the series $\sum_{n=1}^\infty A_n^\nu(\bar\omega) 1_{[T_n^\nu(\bar\omega),T_{n+1}^\nu)}(t)$ is a finite sum, and thus $\bar\alpha^\nu$
$\in$   $\Ac^{W,\mu'}$. Furthermore, by (iv) we see that $\mathscr L_{\P'}(\bar\alpha^\nu(\omega,\cdot))=\mathscr L_{\P^{\nu^\omega}}(I)$, $\omega\notin N$.

Let us now construct $(T_n^\nu,A_n^\nu)_{n\geq1}$ satisfying points (i)-(ii)-(iii)-(iv).
Fix $\nu\in\Vc_{\inf\,>\,0}$ and let  $\tilde N\in\Fc$ be as in Lemma \ref{L:nu_pred}. In particular, recall that formula \eqref{representation_nu^omega} holds for some maps $\nu^{(n)}$, $n\geq0$, satisfying $0<\inf\nu^{(n)}\leq\sup\nu^{(n)}\leq M^\nu$, for some constant $M^\nu>0$, independent of $n$.
Next recall
the construction of the map $\pi:\R\to A$ and the measure $\lambda'$ .
Accordingly, we split
the rest of the proof into two cases.
\vspace{2mm}

\noindent\textbf{Case I:
$A_c=\emptyset$, so that  $A=A_{nc}$ is uncountable.}
In this case $\pi:\R\to A$ is a Borel isomorphism, se
to shorten notation
we identify $A$ with $\R$ and use the notations $A$, $\lambda$, $A_n$, $\mu$, $\F^{W, \mu_\infty}= (\Fc_{t}^{W,\mu_\infty})_{t\ge 0}$ instead of
$\R$, $\lambda'$, $R_n$, $\bar\mu$, $\F^{W,\mu'_\infty}= (\Fc_{t}^{W,\mu'_\infty})_{t\ge 0}$.
Since we are treating the case $A_c=\emptyset$,
 we have $\lambda(\{a\})=0$ for every $a\in A$. We construct by induction on $n\geq1$ a  sequence
$(T_n^\nu,A_n^\nu)_{n\geq1}$ and a $\P$-null set $N\in\Fc$, with $\tilde N\subset N$, such that $(T_n^\nu,A_n^\nu)_{n\geq1}$ satisfies properties (i) and (ii) of the Theorem, and also the following properties:
\begin{itemize}
\item[(iii)'] for every $n\geq1$, we have $T_n^\nu\geq T_n/M^\nu$;
\item[(iv)'] for every $n\geq 1$ and $\omega\notin N$, we have
\begin{equation}\label{SameLaw}
\mathscr L_{\P'}(T_1^\nu(\omega,\cdot),A_1^\nu(\omega,\cdot),\ldots,T_n^\nu(\omega,\cdot),A_n^\nu(\omega,\cdot)) \ = \ \mathscr L_{\P^{\nu^\omega}}(T_1,A_1,\ldots,T_n,A_n).
\end{equation}
\end{itemize}
Notice that (iv)' is equivalent to (iv). Moreover, since $\lim_{n\rightarrow\infty}T_n=\infty$, we see that (iii)' implies property (iii).

\vspace{2mm}

\noindent\textbf{Step  1: the case  $n=1$}.
Define
\begin{equation}\label{vartheta_1}
\theta_t^{(1)}(\omega) \ := \ \frac{1}{\lambda(A)} \int_0^t\int_A \nu_s^{(0)}(\omega,a) \lambda(da)\,ds.
\end{equation}
Since $0<\inf\nu^{(0)}\leq\sup\nu^{(0)}\leq M^\nu$, we see that, for every $\omega\in\Omega$, the map $t\mapsto\theta_t^{(1)}(\omega)$ is continuous, strictly increasing, $\theta_0^{(1)}(\omega)=0$, $\theta_t^{(1)}(\omega)\leq M^\nu t$, and $\theta_t^{(1)}(\omega)\nearrow\infty$ as $t$ goes to infinity. Then there exists a unique $T_1^\nu\colon\bar\Omega\rightarrow\R_+$ such that
\[
\theta_{T_1^\nu(\bar\omega)}^{(1)}(\omega) \ = \ T_1(\omega').
\]
Notice that $T_1^\nu\geq T_1/M_1^\nu$. Moreover, since $T_1>0$, we also have $T_1^\nu>0$. Let $\bar E_{T_1}:=\{(\bar\omega,t)\in\bar\Omega\times\R_+\colon \theta_t^{(1)}(\omega)=T_1(\omega')\}$. Since
the process $(\bar\omega,t)\mapsto (\theta_t^{(1)}(\omega),T_1(\omega'))$
is $\F^{W,\mu_\infty}$-adapted and continuous,
 $\bar E_{T_1}$ is an $\F^{W,\mu_\infty}$-optional set (in fact, predictable). Since $T_1^\nu(\bar\omega)=\inf\{t\in\R_+\colon(\bar\omega,t)\in\bar E_{T_1}\}$
  is the d\'ebut of $\bar E_{T_1}$,
  from Theorem 1.14 of \cite{jacod_book} it follows that $T_1^\nu$ is an $\F^{W,\mu_\infty}$-stopping time. In particular, $T_1^\nu$ is $\bar\Fc$-measurable, therefore there exists a $\P$-null set $N_{T_1^\nu}\in\Fc$ such that $T_1^\nu(\omega,\cdot)$ is $\Fc'$-measurable, whenever $\omega\notin N_{T_1^\nu}$.

Now define
\[
F_b \ := \ \P'(A_1\leq b) \ = \ \frac{  \lambda((-\infty,b])}{\lambda(A)},
\qquad
F_b^{(1)}(\bar\omega) \ := \ \frac{\int_{-\infty}^b \nu_{T_1^\nu(\bar\omega)}^{(0)}(\omega,a)\,\lambda(da)}{\int_{-\infty}^{+\infty}\nu_{T_1^\nu(\bar\omega)}^{(0)}(\omega,a)\,\lambda(da)}.
\]
Since $\inf\nu^{(0)}>0$ and $\lambda(\{a\})=0$ for any $a\in A$, we see that, for every $\bar\omega\in\bar\Omega$, the map $b\mapsto F_b^{(1)}(\bar\omega)$  is continuous, strictly increasing, valued in $(0,1)$, and $\lim_{b\rightarrow-\infty}F_b^{(1)}(\bar\omega)=0$, $\lim_{b\rightarrow+\infty}F_b^{(1)}(\bar\omega)=1$. Then, there exists a unique $A_1^\nu\colon\bar\Omega\rightarrow\R$ such that
\[
F_{A_1^\nu(\bar\omega)}^{(1)}(\bar\omega) \ = \ F_{A_1(\omega')}.
\]
We note that the process
$$
(\bar\omega,t)\mapsto
\frac{
\int_{-\infty}^b \nu_{t}^{(0)}(\omega,a)\,\lambda(da)
}{
\int_{-\infty}^{+\infty}\nu_{t}^{(0)}(\omega,a)\,\lambda(da)
}
$$
is predictable with respect to $\F^W$, hence
it is also $\F^{W,\mu_\infty}$-progressive.
Substituting $t$ with $T_1^\nu(\bar\omega)$ we conclude that
$F_b^{(1)}$ is $(\Fc_{T_1^\nu}^{W,\mu_\infty})$-measurable.
Since $A_1$ is clearly $\calf'$-measurable and
$\calf'\subset \Fc_{0}^{W,\mu_\infty}\subset \Fc_{T_1^\nu}^{W,\mu_\infty}$,
$A_1$  is also $(\Fc_{T_1^\nu}^{W,\mu_\infty})$-measurable.
Recalling the continuity of $b\mapsto F_b^{(1)}(\bar\omega)$
it is easy to conclude that $A_1^\nu$ is $(\Fc_{T_1^\nu}^{W,\mu_\infty})$-measurable.
This implies that $A_1^\nu\vee a$ is also $\Fc_{T_1^\nu}^{W,\mu_\infty}$-measurable. From the arbitrariness of $a$, we deduce that $A_1^\nu$ is $\Fc_{T_1^\nu}^{W,\mu_\infty}$-measurable. In particular, $A_1^\nu$ is $\bar\Fc$-measurable, therefore there exists a $\P$-null set $N_{A_1^\nu}\in\Fc$ such that $A_1^\nu(\omega,\cdot)$ is $\Fc'$-measurable, whenever $\omega\notin N_{A_1^\nu}$.

In order to conclude the proof of the case $n=1$, let us prove that \eqref{SameLaw} holds for $n=1$, whenever $\omega\notin N_1:=\tilde N\cup N_{T_1^\nu}\cup N_{A_1^\nu}$. We begin recalling that, for every $\omega\notin\tilde N$,
the $\F^\mu$-compensator of $\mu$ under $\P^{\nu^\omega}$ is given by the
right-hand side of \eqref{representation_nu^omega}, so that in particular
we have
\[
\P^{\nu^\omega}(T_1>t) \ = \ \exp\bigg(-\int_0^t\int_A \nu_s^{(0)}(\omega,a) \lambda(da)\,ds\bigg) \ = \ \exp\big(-\lambda(A)\theta_t^{(1)}(\omega)\big).
\]
Notice that
\[
\P'\big(T_1^\nu(\omega,\cdot)>t\big) \ = \ \P'\big(\theta_{T_1^\nu(\omega,\cdot)}^{(1)}(\omega)>\theta_t^{(1)}(\omega)\big) \ = \ \P'\big(T_1>\theta_t^{(1)}(\omega)\big) \ = \ \exp\big(-\lambda(A)\theta_t^{(1)}(\omega)\big),
\]
for every $\omega\notin N_{T_1^\nu}$, where for the last equality we used the formula $\P'(T_1>t)=\exp(-\lambda(A)t)$. Therefore $\mathscr L_{\P'}(T_1^\nu(\omega,\cdot))=\mathscr L_{\P^{\nu^\omega}}(T_1)$, for every $\omega\notin\tilde N\cup N_{T_1^\nu}$. Now, recall that, for every $\omega\notin\tilde N$, we have, $\P'$-a.s.,
\[
\P^{\nu^\omega}\big(A_1\leq b\,\big|\,\sigma(T_1)\big) \ = \ \frac{\int_{-\infty}^b \nu_{T_1}^{(0)}(\omega,a)\,\lambda(da)}{\int_{-\infty}^{+\infty}\nu_{T_1}^{(0)}(\omega,a)\,\lambda(da)}.
\]
On the other hand, for every $\omega\notin N_{T_1^\nu}\cup N_{A_1^\nu}$, $\P'$-a.s.,
\begin{align*}
\P'\big(A_1^\nu(\omega,\cdot)\leq b\,\big|\,\sigma(T_1^\nu(\omega,\cdot))\big) \ &= \ \P'\big(F_{A_1^\nu(\omega,\cdot)}^{(1)}(\omega,\cdot)\leq F_b^{(1)}(\omega,\cdot)\,\big|\,\sigma(T_1^\nu(\omega,\cdot))\big) \\
&= \ \P'\big(F_{A_1}\leq F_b^{(1)}(\omega,\cdot)\,\big|\,\sigma(T_1^\nu(\omega,\cdot))\big).
\end{align*}
Since $A_1$ is independent of $T_1$ under $\P'$
and $T_1^\nu(\omega,\cdot)$ is $\sigma(T_1)$-measurable, it follows that $A_1$ is also independent of $T_1^\nu(\omega,\cdot)$. Moreover, by definition we see that $F_b^{(1)}(\omega,\cdot)$ is $\sigma(T_1^\nu(\omega,\cdot))$-measurable. Therefore, for every $\omega\notin N_{T_1^\nu}\cup N_{A_1^\nu}$, we have, $\P'$-a.s.,
\[
\P'\big(F_{A_1}\leq F_b^{(1)}(\omega,\cdot)\,\big|\,\sigma(T_1^\nu(\omega,\cdot))\big) \ = \ \P'\big(F_{A_1}\leq a\big)\big|_{a=F_b^{(1)}(\omega,\cdot)} \ = \ F_b^{(1)}(\omega,\cdot),
\]
where we used the fact that $F_{A_1}$ is uniformly distributed in $(0,1)$ under $\P'$.
As a consequence, recalling that $\mathscr L_{\P'}(T_1^\nu(\omega,\cdot))=\mathscr L_{\P^{\nu^\omega}}(T_1)$, for every $\omega\notin\tilde N\cup N_{T_1^\nu}$, we deduce that $\mathscr L_{\P'}(T_1^\nu(\omega,\cdot),A_1^\nu(\omega,\cdot))=\mathscr L_{\P^{\nu^\omega}}(T_1,A_1)$, whenever $\omega\notin N_1$. This concludes the proof of the base case $n=1$.

\vspace{2mm}

\noindent\textbf{Step  2: the inductive step}.   Fix $n\geq1$ and suppose we are given $(T_1^\nu,A_1^\nu),\ldots,(T_n^\nu,A_n^\nu)$ satisfying points (i) and (ii) of the Theorem. Suppose also that \eqref{SameLaw} holds for the fixed $n$, whenever $\omega\notin N_n$, for some $\P$-null set $N_n\in\Fc$ in place of $N$, with $\tilde N\subset N_n$.

Given $\theta^{(1)}$ as in \eqref{vartheta_1}, we define recursively, for $i=1,\ldots,n$,
\begin{align}\label{vartheta_i}
\theta_t^{(i+1)}(\bar\omega) \ &:= \ \theta_{T_i^\nu(\bar\omega)\wedge t}^{(i)}(\bar\omega) \\
&\quad \ + \frac{1}{\lambda(A)} \int_{T_i^\nu(\bar\omega)}^{T_i^\nu(\bar\omega)\vee t}\int_A \nu_s^{(i)}\big(\omega,(T_1^\nu(\bar\omega),A_1^\nu(\bar\omega)),\ldots,(T_i^\nu(\bar\omega),A_i^\nu(\bar\omega)),a\big) \lambda(da)\,ds. \notag
\end{align}
Since $0<\inf\nu^{(i)}\leq\sup\nu^{(i)}\leq M^\nu$, we see that, for every $\bar\omega\in\bar\Omega$ and $i=1,\ldots,n$, the map $t\mapsto\theta_t^{(i+1)}(\bar\omega)$ is continuous, strictly increasing, $\theta_0^{(i+1)}(\bar\omega)=0$, $\theta_t^{(i+1)}(\omega)\leq M^\nu t$, and
$\theta_t^{(i+1)}(\bar\omega)\nearrow\infty$ as $t$ goes to infinity. Then, there exists a unique $
T_{n+1}^\nu\colon\bar\Omega\rightarrow\R_+$ such that
\[
\theta_{T_{n+1}^\nu(\bar\omega)}^{(n+1)}(\bar\omega) \ = \ T_{n+1}(\omega').
\]
Notice that $T_{n+1}^\nu\geq T_{n+1}/M^\nu$. Moreover, since $T_{n+1}>T_n$, we also have $T_{n+1}^\nu>T_n^\nu$. Indeed, arguing by contradiction, suppose that $T_{n+1}^\nu(\bar\omega)\leq T_n^\nu(\bar\omega)$ for some $\bar\omega\in\bar\Omega$. Then
\[
\theta_{T_n^\nu(\bar\omega)}^{(n)}(\bar\omega) \ = \ T_n(\omega') \ < \ T_{n+1}(\omega') \ = \ \theta_{T_{n+1}^\nu(\bar\omega)}^{(n+1)}(\bar\omega) \ = \ \theta_{T_{n+1}^\nu(\bar\omega)}^{(n)}(\bar\omega),
\]
where the last equality follows from \eqref{vartheta_i}. From the monotonicity of $\theta^{(n)}$, we get
$T_n^\nu(\bar\omega)<T_{n+1}^\nu(\bar\omega)$, which yields a contradiction.

Reasoning in the same way as for $T_1^\nu$, since $T_{n+1}^\nu$ is the d\'ebut of $\bar E_{T_{n+1}}:=\{(\bar\omega,t)\in\bar\Omega\times\R_+\colon \theta_t^{(n+1)}(\bar\omega)=T_{n+1}(\omega')\}$,   it is an $\F^{W,\mu_\infty}$-stopping time. In particular, $T_{n+1}^\nu$ is $\bar\Fc$-measurable, so that there exists a $\P$-null set $N_{T_{n+1}^\nu}\in\Fc$ such that $T_{n+1}^\nu(\omega,\cdot)$ is $\Fc'$-measurable, whenever
$\omega\notin N_{T_{n+1}^\nu}$.

Now define $F_b  =   \P'(A_1\leq b)  =    \lambda((-\infty,b])/\lambda(A)$ as before and
\[
F_b^{(n+1)}(\bar\omega) \ := \ \frac{\int_{-\infty}^b \nu_{T_{n+1}^\nu(\bar\omega)}^{(n)}\big(\omega,(T_1^\nu(\bar\omega),A_1^\nu(\bar\omega)),\ldots,(T_n^\nu(\bar\omega),A_n^\nu(\bar\omega)),a\big)\,\lambda(da)}{\int_{-\infty}^{+\infty}\nu_{T_{n+1}^\nu(\bar\omega)}^{(n)}\big(\omega,(T_1^\nu(\bar\omega),A_1(\bar\omega)),\ldots,(T_n^\nu(\bar\omega),A_n^\nu(\bar\omega)),a\big)\,\lambda(da)}.
\]
Since $\inf\nu^{(n)}>0$ and $\lambda(\{a\})=0$ for any $a\in A$, we see that, for every $\bar\omega\in\bar\Omega$, the map $b\mapsto F_b^{(n+1)}(\bar\omega)$ is continuous, strictly increasing, valued in $(0,1)$, and  $\lim_{b\rightarrow-\infty}F_b^{(n+1)}(\bar\omega)=0$, $\lim_{b\rightarrow+\infty}F_b^{(n+1)}(\bar\omega)=1$. Then, proceeding along the same lines as for the construction of $A_1^\nu$, we see that there exists a unique $\Fc_{T_{n+1}^\nu}^{W,\mu_\infty}$-measurable map $A_{n+1}^\nu\colon\bar\Omega\rightarrow\R$ such that
\[
F_{A_{n+1}^\nu(\bar\omega)}^{(n+1)}(\bar\omega) \ = \ F_{A_{n+1}(\omega')}.
\]
In particular, $A_{n+1}^\nu$ is $\bar\Fc$-measurable, therefore there exists a $\P$-null set $N_{A_{n+1}^\nu}\in\Fc$ such that $A_{n+1}^\nu(\omega,\cdot)$ is $\Fc'$-measurable, whenever $\omega\notin N_{A_{n+1}^\nu}$.

In order to conclude the proof of the inductive step, let us prove that \eqref{SameLaw} holds for $n+1$, whenever $\omega\notin N_{n+1}:=N_n\cup N_{T_{n+1}^\nu}\cup N_{A_{n+1}^\nu}$.
Set  $S_{n+1}=T_{n+1}-T_{n}$ and
 recall that, for every $\omega\notin\tilde N$,
the $\F^\mu$-compensator of $\mu$ under $\P^{\nu^\omega}$ is given by the
right-hand side of \eqref{representation_nu^omega}, so that in particular
we have, $\P'$-a.s.,
\begin{align}\nonumber
&\P^{\nu_\omega}\big(S_{n+1}>t\,\big|\,\sigma(T_1,A_1,\ldots,T_n,A_n)\big) \\
&= \ \exp\bigg(-\int_{T_n}^{T_n + t}\int_A \nu_s^{(n)}\big(\omega,(T_1,A_1),\ldots,(T_n,A_n),a\big) \lambda(da)\,ds\bigg). \label{Sn1condiz}
\end{align}
Define $S_{n+1}^\nu:=T_{n+1}^\nu-T_n^\nu$ and observe that, whenever $\omega\not\in N_n\cup N_{T_{n+1}^\nu}$, $\P'$-a.s.,
\begin{align*}
&\P'\big(S_{n+1}^\nu(\omega,\cdot)>t\,\big|\,\sigma(T_1^\nu(\omega,\cdot),A_1^\nu(\omega,\cdot),\ldots,T_n^\nu(\omega,\cdot),A_n^\nu(\omega,\cdot))\big) \\
&= \ \P'\big(\theta_{T_{n+1}^\nu(\omega,\cdot)}^{(n+1)}(\omega,\cdot)>\theta_{T_n^\nu(\omega,\cdot) + t}^{(n+1)}(\omega,\cdot)\,\big|\,\sigma(T_1^\nu(\omega,\cdot),A_1^\nu(\omega,\cdot),\ldots,T_n^\nu(\omega,\cdot),A_n^\nu(\omega,\cdot))\big) \\
&= \ \P'\bigg(T_{n+1}>\frac{1}{\lambda(A)} \int_{T_n^\nu(\omega,\cdot)}^{T_n^\nu(\omega,\cdot)+t}\int_A \nu_s^{(n)}\big(\omega,T_1^\nu(\omega,\cdot),A_1^\nu(\omega,\cdot),\ldots,T_n^\nu(\omega,\cdot),A_n^\nu(\omega,\cdot),a\big) \lambda(da)\,ds \\
&\quad \ + \theta_{T_n^\nu(\omega,\cdot)}^{(n)}(\omega,\cdot)\,\bigg|\,\sigma(T_1^\nu(\omega,\cdot),A_1^\nu(\omega,\cdot),\ldots,T_n^\nu(\omega,\cdot),A_n^\nu(\omega,\cdot))\bigg) \\
&= \ \P'\bigg(S_{n+1}>\frac{1}{\lambda(A)} \int_{T_n^\nu(\omega,\cdot)}^{T_n^\nu(\omega,\cdot)+t}\int_A \nu_s^{(n)}\big(\omega,\ldots,a\big) \lambda(da)\,ds\,\bigg|\,\sigma(T_1^\nu(\omega,\cdot),\ldots,A_n^\nu(\omega,\cdot))\bigg).
\end{align*}
Recall now that $S_{n+1}$ is independent of $(T_1,A_1,\ldots,T_n,A_n)$
under $\P'$, and note that, by construction, $(T_1^\nu(\omega,\cdot),A_1^\nu(\omega,\cdot),\ldots,T_n^\nu(\omega,\cdot),A_n^\nu(\omega,\cdot))$ is $\sigma(T_1,A_1,\ldots,T_n,A_n)$-measurable. Therefore, $S_{n+1}$ is also independent of $(T_1^\nu(\omega,\cdot),A_1^\nu(\omega,\cdot),\ldots,T_n^\nu(\omega,\cdot),A_n^\nu(\omega,\cdot))$. As a consequence, for every $\omega\not\in N_n\cup N_{T_{n+1}^\nu}$, $\P'$-a.s.,
\begin{align*}
&\P'\bigg(S_{n+1}>\frac{1}{\lambda(A)} \int_{T_n^\nu(\omega,\cdot)}^{T_n^\nu(\omega,\cdot)+t}\int_A \nu_s^{(n)}\big(\omega,\ldots,a\big) \lambda(da)\,ds\,\bigg|\,\sigma(T_1^\nu(\omega,\cdot),\ldots,A_n^\nu(\omega,\cdot))\bigg) \\
&= \ \P'(S_{n+1}>r)\big|_{r=\frac{1}{\lambda(A)} \int_{T_n^\nu(\omega,\cdot)}^{T_n^\nu(\omega,\cdot)+t}\int_A \nu_s^{(n)}(\omega,\ldots,a) \lambda(da)\,ds} \\
&= \ \exp\bigg(-\int_{T_n^\nu(\omega,\cdot)}^{T_n^\nu(\omega,\cdot) + t}\int_A \nu_s^{(n)}\big(\omega,(T_1^\nu(\omega,\cdot),A_1^\nu(\omega,\cdot),\ldots,T_n^\nu(\omega,\cdot),A_n^\nu(\omega,\cdot)),a\big) \lambda(da)\,ds\bigg),
\end{align*}
where for the last equality we used the formula $\P'(S_{n+1}>r)=\exp(-\lambda(A)r)$.
Comparing with
\eqref{Sn1condiz}, we see that the conditional distribution of $S_{n+1}^\nu$ given
$T_1^\nu(\omega,\cdot),A_1^\nu(\omega,\cdot),\ldots,T_n^\nu(\omega,\cdot),A_n^\nu(\omega,\cdot)$
under $\P'$ is the same as the
conditional distribution of $S_{n+1}$ given
$T_1,A_1,\ldots,T_n,A_n$
under
$\P^{\nu^\omega}$.
Together with the inductive assumption \eqref{SameLaw}
this proves that
\begin{equation}\label{SameLaw_S_n+1}
\mathscr L_{\P'}(T_1^\nu(\omega,\cdot),\ldots,T_n^\nu(\omega,\cdot),A_n^\nu(\omega,\cdot),T_{n+1}^\nu(\omega,\cdot)) \ = \ \mathscr L_{\P^{\nu^\omega}}(T_1,\ldots,T_n,A_n,T_{n+1}),
\end{equation}
for every $\omega\notin N_n\cup N_{T_{n+1}^\nu}$. Now, recall that, for every $\omega\notin\tilde N$, $\P'$-a.s.,
\begin{equation}
\P^{\nu^\omega}\big(A_{n+1}\leq b\,\big|\,\sigma(T_1,\ldots,A_n,T_{n+1})\big)
= \ \frac{\int_{-\infty}^b \nu_{T_{n+1}}^{(n)}\big(\omega,(T_1,A_1),\ldots,(T_n,A_n),a\big)\,\lambda(da)}
{\int_{-\infty}^{+\infty}\nu_{T_{n+1}}^{(n)}\big(\omega,(T_1,A_1),\ldots,(T_n,A_n),a\big)\,\lambda(da)}.
\label{An1condiz}
\end{equation}
On the other hand, for every $\omega\not\in N_{n+1}$, we have, $\P'$-a.s.,
\begin{align*}
&\P'\big(A_{n+1}^\nu(\omega,\cdot)\leq b\,\big|\,\sigma(T_1^\nu(\omega,\cdot),\ldots,A_n^\nu(\omega,\cdot),T_{n+1}^\nu(\omega,\cdot))\big) \\
&= \ \P'\big(F_{A_{n+1}^\nu(\omega,\cdot)}^{(n+1)}(\omega,\cdot)\leq F_b^{(n+1)}(\omega,\cdot)\,\big|\,\sigma(T_1^\nu(\omega,\cdot),\ldots,A_n^\nu(\omega,\cdot),T_{n+1}^\nu(\omega,\cdot))\big) \\
&= \ \P'\big(F_{A_{n+1}}\leq F_b^{(n+1)}(\omega,\cdot)\,\big|\,\sigma(T_1^\nu(\omega,\cdot),\ldots,A_n^\nu(\omega,\cdot),T_{n+1}^\nu(\omega,\cdot))\big).
\end{align*}
Since $A_{n+1}$ is independent of $(T_1,\ldots,A_n,T_{n+1})$ under $\P'$ and $(T_1^\nu(\omega,\cdot),\ldots,A_n^\nu(\omega,\cdot),T_{n+1}^\nu(\omega,\cdot))$ is $\sigma(T_1,$ $\ldots,A_n,T_{n+1})$-measurable, it follows that $A_{n+1}$ is also independent of $(T_1^\nu(\omega,\cdot),\ldots,A_n^\nu(\omega,\cdot),$ $T_{n+1}^\nu(\omega,\cdot))$. Moreover, by definition we see that $F_b^{(n+1)}(\omega,\cdot)$ is $\sigma(T_1^\nu(\omega,\cdot),\ldots,A_n^\nu(\omega,\cdot),T_{n+1}^\nu(\omega,\cdot))$-measurable. Therefore, for every $\omega\notin N_{n+1}$, we have, $\P'$-a.s.,
\begin{align*}
&\P'\big(F_{A_{n+1}}\leq F_b^{(n+1)}(\omega,\cdot)\,\big|\,\sigma(T_1^\nu(\omega,\cdot),\ldots,A_n^\nu(\omega,\cdot),T_{n+1}^\nu(\omega,\cdot))\big) \\
&= \ \P'\big(F_{A_{n+1}}\leq a\big) \big|_{a=F_b^{(n+1)}(\omega,\cdot)} \ = \ F_b^{(n+1)}(\omega,\cdot),
\end{align*}
where we used the fact that $F_{A_{n+1}}$ is uniformly distributed in $(0,1)$ under $\P'$.
Comparing with
\eqref{An1condiz}, we see that the conditional distribution of $A_{n+1}^\nu$ given
$T_1^\nu(\omega,\cdot),A_1^\nu(\omega,\cdot),\ldots,T_n^\nu(\omega,\cdot),$
$A_n^\nu(\omega,\cdot),$
$T_{n+1}^\nu(\omega,\cdot)$
under $\P'$ is the same as the
conditional distribution of $A_{n+1}$ given
$T_1,A_1,\ldots,$ $T_n,$ $A_n, $ $T_{n+1}$
under
$\P^{\nu^\omega}$.
Therefore, by \eqref{SameLaw_S_n+1} we deduce that \eqref{SameLaw} holds for $n+1$, whenever $\omega\notin N_{n+1}$, which concludes the proof of the inductive step
and also the proof of Case I.

\vspace{2mm}

\noindent\textbf{Case II: $A_{c}\neq\emptyset$.}
Let $\pi:\R\to A$ be the
canonical projection
\eqref{defproiez}
or \eqref{defproiezdue}
according whether  $A_{nc}=\emptyset$ or $A_{nc}\neq\emptyset$.
The idea of the proof is to construct a
random sequence with values in $(0,\infty)\times \R$ using the
Case I previously addressed, and  obtain
the required sequence $(T_n^\nu,A_n^\nu)_{n\geq1}$
by projecting the second component onto $A$.  The detailed construction
and proof
is presented below in the case $A_{nc}=\emptyset$, the other
one being simpler and entirely analogous.

Given $\nu\in\Vc_{\inf\,>\,0}$, define $\bar\nu=\bar\nu_t(\omega, \omega',r): \bar\Omega\times \R_+\times\R\to (0,\infty)$ by
\[
\bar\nu_t(\omega, \omega',r) \ := \ \nu_t(\omega, \omega',\pi(r)).
\]
By a monotone class argument  we see that $\bar\nu$ is $\calp(\F^{W,\mu'})\otimes \calb(\R)$-measurable.   Then
we can perform the construction  presented in step I, with
$\R$, $\lambda'$, $\bar\nu$, $R_n$, $\mu'$,
$\F^{W,\mu'_\infty}= (\Fc_{t}^{W,\mu'_\infty})_{t\ge 0}$ instead of   $A$, $\lambda$, $\nu$, $A_n$, $\mu$,
$\F^{W, \mu_\infty}= (\Fc_{t}^{W,\mu_\infty})_{t\ge 0}$,  respectively.
This way we obtain a $\P$-null set $N\in\Fc$  and a sequence $(\bar T_n^{\bar\nu},\bar R_n^{\bar\nu})_{n\geq1}$
 with values in $(0,\infty)\times \R$ such that   $\bar T_n^{\bar\nu}
 <\bar T_{n+1}^{\bar\nu}\nearrow\infty$,
  $\bar T_n^{\bar\nu}$ is an $\F^{W, \bar\mu_\infty}$-stopping time and
  $\bar R_n^{\bar\nu}$ is $\Fc_{\bar T_n^{\bar\nu}}^{W,\bar\mu_\infty}$-measurable,
  and
\begin{equation}\label{SameLawBar}
\mathscr L_{\P'}\big((\bar T_n^{\bar\nu}(\omega,\cdot),
\bar R_n^{\bar\nu}(\omega,\cdot))_{n\geq1}\big) \ = \ \mathscr L_{\P^{\bar\nu^\omega}}\big((T_n,R_n)_{n\geq1}\big)
\end{equation}
for every $\omega\notin N$. We define
the required sequence
 $(T_n^\nu,A_n^\nu)_{n\ge 1}$ setting
$$
T_n^\nu:= \bar T_n^{\bar\nu},
\qquad
A_n^\nu:= \pi(\bar R_n^{\bar\nu}).
$$
Clearly,  conditions (i)-(ii)-(iii)  of the Theorem hold true and,  recalling the notation $A_n=\pi(R_n)$,  from \eqref{SameLawBar}
 it follows that
\begin{equation}\label{SameLawBarBis}
\mathscr L_{\P'}\big(( T_n^{ \nu}(\omega,\cdot),
 A_n^{ \nu}(\omega,\cdot))_{n\geq1}\big) \ = \ \mathscr L_{\P^{\bar\nu^\omega}}\big((T_n,A_n)_{n\geq1}\big).
\end{equation}
 Next note that, by the definition of $\bar\nu$,
 \beq \nonumber
\int_\R (1 - \bar \nu_t(r))\lambda'(dr)&=&
\int_{(0,\infty)} (1 - \bar \nu_t(r))\lambda'(dr)+
\sum_{j\in A_c}\int_{I_j} (1 - \bar \nu_t(r))\lambda'(dr)
\\\nonumber
&=&
\int_{A_{nc}} (1 -   \nu_t(s))\lambda(da)+
\sum_{j\in A_c}(1 -  \nu_t(j))\lambda(\{j\})
\\\nonumber
&=&
\int_A(1 -  \nu_s(a))\lambda(da)
\enq
and
$ \bar \nu_{T_n}(R_n)=  \nu_{T_n}(\pi(R_n))=  \nu_{T_n}(A_n)$, so that
\beq \nonumber
\kappa_t^{\bar\nu}
&=& \ \exp\left(\int_0^t\int_\R (1 - \bar \nu_s(r))\lambda'(dr)\,ds
\right)\prod_{T_n\le t}\bar\nu_{T_n}(R_n)
\\\nonumber
&=& \ \exp\left(\int_0^t\int_A (1 -  \nu_s(a))\lambda(da)\,ds
\right)\prod_{T_n\le t}\nu_{T_n}(A_n)
\\\nonumber
&=&  \kappa_t^{\nu} .
\enq
Therefore we have, for  $\omega\notin N$, on every $\sigma$-algebra
$\calf^\mu_t$,
$$
\P^{\bar\nu^\omega}(d\omega')=  \kappa_t^{\bar\nu} (\omega,\omega')\,\P'(d\omega')
= \kappa_t^{\nu} (\omega,\omega')\,\P'(d\omega') = \P^{\nu^\omega}(d\omega')
$$
and   property (iv) of the Theorem follows from \eqref{SameLawBarBis}. The
 proof is finished.
\ep

\vspace{3mm}

We can now prove the main result of this subsection and conclude the proof of the inequality ${\text{\Large$\upsilon$}}_0^\Rc$ $\le$  ${\text{\Large$\upsilon$}}_0$.

\begin{Proposition}
\label{P:Ineq_I}
For every $\nu\in\Vc_{\inf\,>\,0}$ there exists $\bar\alpha^\nu\in\Ac^{W,\mu'}$ such that
\begin{equation}\label{LawB_I=LawB_alpha}
\mathscr L_{\P^\nu}(x_0,B,I) \ = \ \mathscr L_{\bar\P}(x_0,B,\bar\alpha^\nu).
\end{equation}
In particular,  $V$ and $W$ are standard Wiener processes,  $V,W,x_0$ are all independent  under $\P^\nu$, and we have
\begin{equation}\label{LawX_I=LawX_alpha}
\mathscr L_{\P^\nu}(X,I)  =  \mathscr L_{\bar\P}(X^{\bar\alpha^\nu},\bar\alpha^\nu),
\qquad J^\Rc(\nu)=\bar J(\bar\alpha^\nu).
\end{equation}
Finally, ${\text{\Large$\upsilon$}}_0^\Rc$ $:=$ $\sup_{\nu\in\Vc}J^\Rc(\nu)$ $\leq$ $\sup_{\bar\alpha\in\Ac^{W,\mu'}}\bar J(\bar\alpha)$  $=$
$\sup_{\alpha\in\Ac^W}J(\alpha)$ $=:$ ${\text{\Large$\upsilon$}}_0$.
\end{Proposition}
\textbf{Proof.}
Suppose that \eqref{LawB_I=LawB_alpha} holds. Then, from equation \eqref{stateeq} and Assumption {\bf (A1)} it is well-known that the first
equality in \eqref{LawX_I=LawX_alpha} holds as well, and this   implies the second equality.
 From the arbitrariness of $\nu\in\Vc_{\inf\,>\,0}$, we deduce that $\sup_{\nu\in\Vc_{\inf\,>\,0}}J^\Rc(\nu)\leq\sup_{\bar\alpha\in\Ac^{W,\mu'}}\bar J(\bar\alpha)$. Since  $\sup_{\nu\in\Vc_{\inf\,>\,0}}J^\Rc(\nu)=\sup_{\nu\in\Vc}J^\Rc(\nu)$
 by \eqref{eqinfzero}, we conclude that $\sup_{\nu\in\Vc}J^\Rc(\nu)\leq\sup_{\bar\alpha\in\Ac^{W,\mu'}}\bar J(\bar\alpha) = \sup_{\alpha\in\Ac^W}J(\alpha)$, where the last equality follows from Lemma \ref{L:AWmu}.

Let us now prove \eqref{LawB_I=LawB_alpha}. Fix $\nu\in\Vc_{\inf\,>\,0}$ and consider the process $\bar\alpha^\nu$ given by Lemma \ref{L:Replication_I}. In order to prove \eqref{LawB_I=LawB_alpha}, we have to show that
\begin{equation}\label{psi_phi}
\bar\E[\kappa_T^\nu\chi(x_0)\psi(B)\phi(I)] \ = \ \bar\E[\chi(x_0)\psi(B)\phi(\bar\alpha^\nu)],
\end{equation}
for any $\chi\in B_b(\R^n)$ (the space of bounded Borel measurable real functions on $\R^n$),
for any $\psi\in B_b({\bf C}_{m+d})$ (the space of bounded Borel measurable real functions on ${\bf C}_{m+d}$, which denotes the space of continuous paths from $[0,T]$ to $\R^{m+d}$ endowed with the supremum norm) and any $\phi\in B_b({\bf D}_A)$ (the space of bounded Borel measurable real functions on ${\bf D}_A$, which denotes the space of c\`adl\`ag paths from $[0,T]$ to $A$ endowed with the supremum norm). By the Fubini theorem,  \eqref{psi_phi} can be rewritten as
\begin{align*}
&\int_\Omega\chi(x_0(\omega)) \psi(B(\omega)) \bigg(\int_{\Omega'} \kappa_T^\nu(\omega,\omega') \phi(I(\omega')) \P'(d\omega')\bigg) \P(d\omega) \\
&= \ \int_\Omega\chi(x_0(\omega)) \psi(B(\omega)) \bigg(\int_{\Omega'} \phi(\bar\alpha^\nu(\omega,\omega')) \P'(d\omega')\bigg) \P(d\omega).
\end{align*}
Let $\tilde N\in\Fc$ be as in Lemma \ref{L:nu_pred}. Then we have to prove that $\E'[\kappa_T^\nu(\omega,\cdot)\phi(I)]=\E'[\phi(\bar\alpha^\nu(\omega,\cdot))]$, whenever $\omega\notin\tilde N$, or, equivalently by definition of $\P^{\nu^\omega}$:
\[
\E^{\nu^\omega}[\phi(I)] \ = \ \E'[\phi(\bar\alpha^\nu(\omega,\cdot))], \qquad \text{whenever }\omega\notin\tilde N.
\]
This is a direct consequence of the last statement of Lemma \ref{L:Replication_I}.
\ep

\subsection{Proof of the inequality ${\text{\Large$\upsilon$}}_0$ $\le$ ${\text{\Large$\upsilon$}}_0^\Rc$}
\label{secVleqVR}

In this proof we borrow some constructions from
 \cite{FP15} Proposition 4.1,  but we
 need to  obtain improved results and we
 simplify  considerably some arguments.

Suppose we are given a setting $(\Omega,\calf,\P, \F, V,W, x_0)$ for the optimal control problem
with partial observation, satisfying the conditions in Section \ref{Primal}, and  consider the controlled equation \eqref{stateeq} and the gain
\eqref{gaineq}. We fix an $\F^W$-progressive process $\alpha$ with values in $A$.  We will show how to construct a sequence of settings
$(\hat \Omega, \hat \calf,\hat \P_k, \hat V,\hat W, \hat \mu_k, \hat x_0)_k$
for the randomized control problem of Section \ref{randomizedformulation},
and a sequence $(\hat\nu^k)_k$ of corresponding admissible controls
(both sequences depending on $\alpha$),
  such that for the corresponding gains, defined by \eqref{defJrandomized}, we have:
\beq \label{desiredineq}
J^\Rc(\hat\nu^k) & \rightarrow &  J(\alpha), \;\;\; \mbox{ as } \; k \rightarrow \infty.
\enq
Admitting for a moment that this has been done, the proof
is easily concluded by the following arguments. By \reff{desiredineq}, we can find,
 for any $\eps$ $>$ $0$, some $k$ such that $J^\Rc(\hat\nu^k)$ $>$ $J(\alpha)-\eps$.
Since $J^\Rc(\hat\nu^k)$ is a gain of a randomized control problem, it can not exceed the value ${\text{\Large$\upsilon$}}_0^\Rc$ defined in \eqref{dualvalue} which,
by Proposition \ref{indepofthesetting}, does not depend on $\epsilon$ nor on $\alpha$.
It follows that
\beqs
   {\text{\Large$\upsilon$}}_0^\Rc &>&  J(\alpha)-\epsilon
\enqs
and by the arbitrariness of $\epsilon $ and $\alpha$, we obtain the required inequality ${\text{\Large$\upsilon$}}_0^\Rc$ $\ge$ ${\text{\Large$\upsilon$}}_0$.

In order to construct the sequences $(\hat \Omega, \hat \calf,\hat \P_k, \hat V,\hat W, \hat \mu_k, \hat x_0)_k$ and   $(\hat \nu^k)_k$
satisfying \eqref{desiredineq}, we apply Proposition \ref{extensionapproximation}, in the Appendix
below,  to the
probability  space $(\Omega,\calf,\P)$ of the partially observed control problem and to the filtration $\G:=\F^W$.
In that Proposition a suitable probability space $(\Omega',\calf',\P')$ is introduced and the product space $(\hat\Omega,
\hat \calf,\Q)$ is constructed:
\beqs
\hat\Omega \; = \; \Omega\times \Omega',
\qquad
\hat \calf \; = \; \calf\otimes \calf',
\qquad
\Q \; = \; \P\otimes \P'.
\enqs
Then  the random variable $x_0$ and
the processes $\alpha$ and $B=(V, W)$ are extended to $\hat\Omega$
in a natural way. We denote $\hat x_0$ and $\hat\alpha$ the extensions
of $x_0$ and $\alpha$.
The extension of $B$, denoted $\hat B=(\hat V,\hat W)$, remains a Wiener process
under $\Q$. The filtration $\F^W$ can also be canonically extended
to a filtration in  $(\hat\Omega,
\hat \calf)$, which coincides with the filtration $\F^{\hat W}$ generated
by $\hat W$.

Following \cite{80Krylov}, for any pair  $\alpha^1,\alpha^2:\hat\Omega\times [0,T]\to {\bf A}$ of   measurable
processes in $(\hat\Omega,\hat \calf,\Q)$  we define a distance $\tilde\rho(\alpha^1,\alpha^2)$ setting
\beq\label{MetricKrylov}
\tilde \rho(\alpha^1,\alpha^2) &=&
\E^\Q\Big[\int_0^T
\rho(\alpha_t^1,\alpha_t^2) dt\Big],
\enq
where  $\E^\Q$ denotes the expectation under  $\Q$, and $\rho$ is a metric in $A$ satisfying $\rho$ $<$ $1$.
By Proposition
\ref{extensionapproximation},
  for any integer $k\ge1$
 there exists a marked point process $(\hat S_n^k,\hat \eta^k_n)_{n\ge 1}$ defined in
$(\hat\Omega,
\hat \calf,\Q)$
satisfying the following conditions.
\begin{enumerate}
\item
Setting
$
\hat S_0^k=0$, $\hat \eta_0^k=a_0$, $
\hat I_t^k=\sum_{n\ge 0}\hat  \eta^k_{n}1_{ [\hat S^k_n,\hat S^k_{n+1})}(t)$,
we have
$    \tilde\rho (\hat I^k, \hat \alpha)<1/k$.
\item Denote
$\hat \mu_k=\sum_{n\ge1}\delta_{(\hat S_n^k,\hat \eta^k_n)}$ the random measure
associated to $(\hat S_n^k,\hat \eta_n^k)_{n\ge 1}$ and
$\F^{\hat\mu_k}=(\calf_t^{\hat\mu_k})_{t\geq 0}$
 the natural filtration of
$\hat\mu_k$;
then
the  compensator of    $\hat \mu_k$ under $\Q$
with respect to the filtration $ (\calf^{\hat W}_t\vee\calf_t^{\hat\mu_k})_{t\geq 0}$
is absolutely continuous
 with respect to $\lambda(da)\,dt$ and
it can be written in the form
$$
 \hat\nu_t^k(\hat\omega,a)\, \lambda(da)\,dt
$$
 for  some  nonnegative $\calp(\F^{\hat W,\hat \mu})
 \otimes \calb(A)$-measurable  function $\hat\nu^k$ satisfying
$
 \inf_{\hat\Omega\times [0,T]\times A}\hat\nu^k>0$ and
$
  \sup_{\hat\Omega\times [0,T]\times A}\hat\nu^k<\infty$.
\end{enumerate}
We note that $\hat\mu_k$ (and so also $\hat I^k$ and $\hat\nu^k$)
depend on $\alpha$ as well, but we do not make it explicit in the notation.
Let us now consider the completion of the probability
space $(\hat\Omega,\hat\calf,\Q)$, that will be denoted
by the same symbol for simplicity of notation, and let
$\caln$ denote the family of $\Q$-null sets of the completion.
Then the filtration
$(\calf^{\hat W}_t\vee\calf_t^{\hat\mu_k}\vee \caln)_{t\geq 0}$
coincides with  the filtration previously denoted by
$\F^{\hat W,\hat \mu_k}=(\calf^{\hat W,\hat \mu_k}_t)_{t\ge 0}$
(compare with formula  \eqref{expandedfiltration} in
  section \ref{randomizedformulation}). It is easy to see that
$ \hat\nu_t^k(\hat\omega,a)\, \lambda(da)\,dt$ is the compensator
of $\hat\mu_k$  with respect to $\F^{\hat W,\hat \mu_k}$ and the extended
probability $\Q$ as well.

Using the Girsanov theorem for point processes (see e.g. \cite{ja}) we next
construct an equivalent probability under which $\hat\mu_k$ becomes
a Poisson random measure with intensity $\lambda$.
Since $\hat\nu^k$ is a strictly positive $\calp(\F^{\hat W,\hat\mu_k})\otimes\calb(A)$-measurable
random field with bounded inverse, the Dol\'eans exponential process
\begin{equation}\label{defgirsanovdensity}
 M_t^k \; := \; \exp\Big(\int_0^t\int_A(1-\hat\nu^k_s(a)^{-1})\,
 \hat\nu_t^k(a) \lambda(da)\,ds\Big) \prod_{\hat S_n^k\le t} \hat\nu^k_{\hat S^k_n}
 (\hat\eta^k_n)^{-1},\qquad t\in [0,T],
\end{equation}
is a strictly positive martingale (with respect to $\F^{\hat W,\hat\mu_k}$ and $\Q$),
and we can
define an equivalent probability $\hat\P_k$ on the space $(\hat\Omega,\hat\calf)$ setting
 $\hat\P_k(d\hat\omega)$ $=$ $M_T^k(\hat\omega)\Q(d\hat\omega)$. The expectation under
 $\hat\P_k$ will be denoted  $\hat\E_k$. By the Girsanov theorem,
 the restriction of
 $\hat\mu_k$      to $(0,T]\times A$
  has $(\hat\P_k,\F^{\hat W,\hat\mu_k})$-compensator $\lambda(da)\,dt$,
so that  in particular    it is a Poisson random measure.
It can also be proved by standard arguments (see e.g. \cite{FP15},
page 2155,
for   detailed verifications in a similar framework)  that
$\hat B$ is a $(\hat\P_k,\F^{\hat W,\hat\mu_k})$-Wiener process and
that $\hat B$ and $\hat \mu_k$ are independent under $\hat\P_k$.
We have thus constructed a setting
$(\hat \Omega, \hat \calf,\hat \P_k, \hat V,\hat W, \hat \mu_k, \hat x_0)$ for a randomized
control problem as in Section \ref{randomizedformulation}.
Since $\hat\nu^k$ is a bounded,  strictly positive and
$\calp(\F^{W,\hat\mu})\otimes\calb(A)$-measurable
random field
  it  belongs to the class $\hat \calv^k$ of admissible controls for the
  randomized control problem and we now proceed
  to evaluating its gain $J^\Rc(\hat\nu^k)$ and to comparing it with $J(\alpha)$.
Our  aim is to show that,
as a consequence of the fact that
$    \tilde\rho (\hat I^k, \hat \alpha)<1/k$,
we have $J^\Rc(\hat\nu^k)   \to J(\alpha)$ as $k\to\infty$.

 We introduce the
 Dol\'eans exponential process $\kappa^{\hat \nu^k}$
 corresponding to $\hat\nu^k$
  (compare formula
\eqref{doleans}):
\begin{equation}\label{doleansbis}
\kappa_t^{\hat\nu^k} =
 \exp\left(\int_0^t\int_A (1 - \hat \nu^k_s(a))\lambda(da)\,ds
\right)\prod_{\hat S_n^k\le t}\nu^k_{\hat S^k_n}(\hat \eta^k_n),\qquad t\in [0,T],
\end{equation}
we define
 the probability
$d\hat\P_k^{\hat\nu^k} =\kappa_T^{\hat\nu^k} d\hat\P_k $ and we obtain
the  gain
$$
J^\Rc(\hat\nu) =  \hat\E^{\hat\nu^k}
\left[\int_0^Tf_t(\hat X^k,\hat I^{k})\,dt+g(\hat  X^k)\right],
$$
where      $\hat X^k$ is the
solution to
the equation
\begin{equation}\label{controlledhatkappa}
d\hat X_t ^k=  b_t( \hat X^k, \hat I^{k})\,dt
+ \sigma_t(\hat X^k, \hat I^{k})\,d\hat B_t,
\qquad
\hat X_0^k =\hat x_0.
\end{equation}
  However comparing \eqref{defgirsanovdensity}  and
 \eqref{doleansbis}  shows that
  $\kappa^{\hat \nu^k}_T\,M_T^k\equiv 1$, so that
 the  Girsanov transformation $\hat\P_k\mapsto \hat\P_k^{\hat\nu^k}$
is the inverse to the transformation $\Q\mapsto \hat\P_k$ made above,
and   changes back the probability $\hat\P_k$ into
$\Q$ considered above. Therefore we have
$\hat\P_k^{\hat\nu^k}=\Q$ and also
\beq \label{JR}
J^\Rc(\hat\nu^k) &=&  \E^\Q \left[ \int_0^Tf_t(\hat X^k,\hat I^{k})\,dt+g(\hat X^k)\right].
\enq
On the other hand, the gain $J(\alpha)$ of the initial control problem with partial
observation was defined in \eqref{gaineq} in terms
of the solution $X^\alpha$ to the controlled equation \eqref{stateeq}.
Denoting $\hat X^\alpha$ the extension
of $X^\alpha$
to $\hat\Omega$, it is easy to verify
that it is the solution to
\begin{equation}\label{controlledhat}
        d\hat X_t^\alpha \ = \ b_t(\hat  X^\alpha,\hat  \alpha)\,dt +
\sigma_t(\hat  X^\alpha, \hat \alpha)\,d\hat B_t, \qquad \hat X_0^\alpha=\hat x_0,
\end{equation}
and that
\beq \label{Ja}
J(\alpha) & = & \E^\Q\left[\int_0^Tf_t(\hat X^\alpha,\hat \alpha)\,dt+
g(\hat X^\alpha)\right].
\enq
Equations \eqref{controlledhat} and \eqref{controlledhatkappa} are
 considered in the same probability space $(\hat\Omega,
\hat \calf,\Q)$.  In \eqref{controlledhatkappa}  we find a
solution adapted to the filtration
$\G^k:=\F^{\hat x_0,\hat B,\hat \mu_k}$
(defined as in
\eqref{expandedfiltration}) and in
\eqref{controlledhat}
 we find a
solution adapted to the filtration
$\G^0:=\F^{\hat x_0,\hat B}$ generated by $\hat x_0$ and $\hat B$
(since $\alpha$ was $\F^W$-progressive and so
$\hat\alpha$ is progressive with respect to $\F^{\hat W}\subset
\F^{\hat x_0,\hat B}$).

\vspace{2mm}

In order to conclude, we need the following stability lemma,
where the continuity condition {\bf (A1)}-(iii) plays a role.

\begin{Lemma}
\label{contrhotilde}
Given a probability space $(\hat\Omega,
\hat \calf,\Q)$ with filtrations
$\G^k=(\calg_t^k)_{t\ge 0}$ ($k\ge 0$)
consider
the equations
$$
d  Y_t ^k=  b_t( Y^k, \gamma^{k})\,dt
+ \sigma_t(Y^k, \gamma^{k})\,d\beta_t,
\qquad
\hat Y_0^k =y_0,
$$
where $\beta$  is a Wiener process with respect to each $\G^k$,
$\E^\Q|y_0|^p<\infty$,
 $y_0$ is $\calg_t^k$-measurable
and $\gamma^k$ is $\G^k$-progressive for every $k$.
If $\tilde\rho (\gamma^k, \gamma^0)\to 0$
as $k\to\infty$,
then
$$
\E^\Q\sup_{t\in [0,T]}|Y^k_t-Y^0_t|^p\to 0,
\quad
\E^\Q \Big[ \int_0^Tf_t(Y^k,\gamma^{k})\,dt+g(Y^k)\Big]
\to
\E^\Q \Big[ \int_0^Tf_t(Y^0,\gamma^{0})\,dt+g(Y^0)\Big].
$$
\end{Lemma}
\noindent {\bf Proof.}
This stability result
for control problems was first proved in \cite{80Krylov} in the standard diffusion case.
 The extension to the non-Markovian case presented
 in \cite{FP15}, Lemma 4.1 and Remark 4.1, also holds in our case
 with the same proof, using the continuity assumption
 {\bf (A1)}-(iii) and the Lipschitz and growth conditions
 \eqref{lipbsig}-\eqref{PolGrowth_f_g}.
\ep

\vspace{3mm}

Applying the Lemma to $\beta=\hat B$, $Y^k=\hat X^k$,
$\gamma^k=\hat I^k$
 (for $k\ge 1$) and $Y^0=\hat X^\alpha$, $\gamma^0=\hat \alpha$,
and recalling that $\tilde\rho (\hat I^k, \hat \alpha)<1/k\to 0$ we conclude by \reff{JR}, \reff{Ja} that   $J^\Rc(\hat\nu^k)   \to J(\alpha)$ as $k\to\infty$.
Therefore relation  \eqref{desiredineq} is satisfied for this choice of the sequence $(\hat\nu^k)_k$ and for the corresponding  settings $(\hat \Omega, \hat \calf,\hat \P_k, \hat V,\hat W, \hat \mu_k, \hat x_0)_k$.
This ends the proof of the inequality ${\text{\Large$\upsilon$}}_0$ $\le$ ${\text{\Large$\upsilon$}}_0^\Rc$.

\subsection{Some extensions}
\label{Sec:locallylip}

The aim of this paragraph  is to generalize Theorem \ref{MainThm} by weakening the required assumptions. We introduce the following hypothesis on the data $A$, $b,\sigma,f,g$.

\vspace{3mm}

\noindent {\bf (A3)}

\medskip

Points (i)-(ii)-(iii)-(v) of  Assumption {\bf (A1)} are supposed to hold. Moreover,

\begin{itemize}
\item [(iv)']
 There exist nonnegative constants  $M$, $K$, $r$ and, for
 any integer $N\ge 1$, there exist constants $L_N\ge 0$
such that
\beq
(x)^*_t \le N,\; (x')^*_t \le N & \Rightarrow&
|b_t(x,a) - b_t(x',a)| + |\sigma_t(x,a)-\sigma_t(x',a)|
 \leq  L_N  (x-x')^*_t; \label{loclipbsig} \\
|b_t(x,a)| + |\sigma_t(x,a)| & \leq & M  (1 +  (x)^*_t  );  \label{lineargrowthbsig}
\\
|f_t(x,a)| + |g(x)| & \leq & K \big(1 + \|x\|_{_\infty}^r \big); \label{PolGrowth_f_g_bis}
\enq
for all $(t,x,x',a)$ $\in$ $[0,T]\times\bfC_n \times\bfC_n \times \bfM_A$.

\end{itemize}

\vspace{2mm}

Thus, the global Lipschitz condition
\eqref{lipbsig} is replaced by the local Lipschitz condition \eqref{loclipbsig}  and by
the linear growth property
\eqref{lineargrowthbsig},
whereas
\eqref{PolGrowth_f_g_bis} is the same as \eqref{PolGrowth_f_g}.
Again, instead of the measurability condition {\bf (A1)}-(ii) we may assume
condition (ii)'
in Remark  \ref{A1_misurabilita}.

 Besides the desire of greater generality,
we need   this extension in order  to include the motivating examples
 presented in paragraphs  \ref{SubS:Finance} and \ref{SubS:ClassicalPr} under
 the scope of our results.
 It is an easy task, left to the reader, to formulate assumptions on
the data  $\bar b$, $\bar\sigma^1$, $\bar\sigma^2$, $k$, $h$, $\bar\beta$, $\gamma^1$, $\gamma^2$, $\bar f$ and $\bar g$ of the optimal control problem
formulated in paragraph  \ref{SubS:Finance} in such a way that the
state equations
\eqref{statouno}-\eqref{statodue}-\eqref{statotre}-\eqref{dynMP}
for the
four-component process  $X^\alpha$ $=$ $(\bar X^\alpha,M,O,Z)$
and the gain functional \eqref{gainrecast}, after reformulation in the
form  \eqref{stateeq}-\eqref{gaineq},
  satisfy
the requirements in  {\bf (A3)}. A similar consideration applies to
the data  $\bar b$, $\bar\sigma^1$, $\bar\sigma^2$, $k$, $h$,
$\bar f$ and $\bar g$ of the optimal control problem
formulated in paragraph  \ref{SubS:ClassicalPr}, so that both examples
can be treated by our results.
However, the global Lipschitz condition in {\bf (A1)} is not satisfied
except in trivial cases, due to the occurrence of the linear
(hence unbounded) term  $Z_t$ in equation
\eqref{statodue} (respectively,  the term $Z_t^\alpha$
in equation $\eqref{statoduebis}$).

From now on we will assume that assumptions {\bf (A2)} and {\bf (A3)}
are satisfied. We will formulate the partially observed control
problem and the randomized control problem, showing that
the corresponding values are well defined. This is done by a
classical truncation method and a priori inequalities
where however the dependence on the control needs to be carefully
studied.

In order to formulate
the partially observed control problem we fix a setting
$(\Omega,\calf,\P, \F, V,W, x_0)$ as in Subsection \ref{Primal}
and we introduce the same state equation \eqref{stateeq} and gain functional
\eqref{gaineq}.

Next, we fix a continuously differentiable function $\eta:\R\to [0,1]$ such that
 $\eta(r)=1$ for $r\le 0$,  $\eta(r)=0$ for $r\ge 1$ and we introduce
 the truncated coefficients
 $$
 b_t^N(x,a)= b_t(x,a)\;\eta \big((x)^*_t - N\big),
 \qquad   \sigma_t^N(x,a)=\sigma_t(x,a)\;\eta \big((x)^*_t - N\big).
 $$
 We note that $b^N$, $\sigma^N$ satisfy assumption {\bf (A1)}
and we introduce  approximating optimal control problems
given by the state equation
\begin{equation}\label{stateeq_N}
    dX_t^{\alpha,N} \ = \ b_t^N( X^{\alpha,N}, \alpha)\,dt +
\sigma_t^N( X^{\alpha,N}, \alpha)\,dB_t
\end{equation}
on the interval $[0,T]$ with initial condition $X_0^{\alpha,N}=x_0$,
the gain functional
\begin{equation}\label{gaineq_N}
J^N(\alpha) \ = \ \E\left[\int_0^Tf_t(X^{\alpha,N},\alpha)\,dt+g(X^{\alpha,N})\right]
\end{equation}
and the corresponding value
\begin{equation}\label{primalvalue_N}
{\text{\Large$\upsilon$}}^N_0 = \inf_{\alpha\in\Ac^W} J^N(\alpha)
\end{equation}
where we maximize over the same set $\Ac^W$ of admissible controls.

By a standard use of the Burkholder-Davis-Gundy inequalities and the Gronwall
lemma, one proves that  for every $\alpha\in\Ac^W$
there exists a unique $\F$-adapted strong solution  $(X_t^{\alpha,N})_{0\leq t \leq T}$   to \eqref{stateeq_N}
with continuous trajectories and satisfying the estimate
\begin{equation}\label{stimaunif_N}
    \E\,\Big[\sup_{t\in [0,T]}|X_t^{\alpha,N}|^p\Big]   \le
    (1+\E\,|x_0|^p)\, C(p,T,M)
\end{equation}
where $C(p,T,M)$ denotes a constant depending only on the indicated parameters
but independent of $N$ and $\alpha$. Moreover setting
$$
\tau_\alpha^N=\inf\{t\in [0,T]\,:\, |X_t^{\alpha,N}|\ge N\},
$$
with the convention that $\inf \emptyset =T$, one easily checks that
$X^{\alpha,N+1}=X^{\alpha,N}$ on $[0,\tau_\alpha^N]$. We also clearly have
$$
    \E\,\Big[\sup_{t\in [0,\tau^N_\alpha]}|X_t^{\alpha,N}|^p\Big]   \le
    (1+\E\,|x_0|^p)\, C(p,T,M)
$$
and since $\sup_{t\in [0,\tau^N_\alpha]}|X_t^{\alpha,N}|
=\sup_{t\in [0,\tau^N_\alpha]}|X_t^{\alpha,N+1}|
\le \sup_{t\in [0,\tau^{N+1}_\alpha]}|X_t^{\alpha,N+1}|$ we have,
by monotone convergence,
\begin{equation}\label{stimauniftau_N}
        \E\,\Big[\sup_N\sup_{t\in [0,\tau^N_\alpha]}|X_t^{\alpha,N}|^p\Big]
    =\lim_{N\to\infty}\E\,\Big[\sup_{t\in [0,\tau^N_\alpha]}|X_t^{\alpha,N}|^p\Big]
    \le
    (1+\E\,|x_0|^p)\, C(p,T,M).
\end{equation}
Since $\tau_\alpha^N<T$ implies $|X_{\tau^N_\alpha}^{\alpha,N}|\ge N$, we have
$
\{\tau_\alpha^N<T \text{ for infinitely many } N\}\subset
\{ \sup_N|X_{\tau^N_\alpha}^{\alpha,N} |=\infty\}
$
and the latter is a null set by \eqref{stimauniftau_N}. It follows that
$\P$-a.s. we have $\tau_\alpha^N=T$ for all but a finite number of integers $N$.
We can then define, up to indistinguishability, a unique process
$(X_t^{\alpha})_{0\leq t \leq T}$   such that $X^{\alpha}=X^{\alpha,N}$
on $[0,\tau_\alpha^N]$ for every $N$.  It is easily checked that $X^{\alpha}$
is the unique strong solution to the state equation \eqref{stateeq}, it
satisfies
\begin{equation}\label{stimaunifinalpha}
\E\,\Big[\sup_{t\in [0,T]}|X_t^{\alpha}|^p\Big]
    \le
    (1+\E\,|x_0|^p)\, C(p,T,M)
\end{equation}
for every $\alpha\in\Ac^W$
 and the gain functional
\eqref{gaineq} and the value  \eqref{primalvalue}
are well defined and finite.

Next we proceed to formulate
the randomized control problem by fixing a setting
$(\hat \Omega, \hat \calf$, $\hat \P$, $\hat V$, $\hat W$, $\hat \mu, \hat x_0)$
as in Section
\ref{randomizedformulation}, defining  the
 filtrations
$\F^{\hat W,\hat \mu}$,
$\F^{\hat x_0,\hat B,\hat \mu}$
and
the process
$\hat I$ by formula
\eqref{I} as before, and considering the
equation \eqref{dynXrandom} for the process $\hat X$.
We  introduce the same set $\hat \calv$ of admissible controls
and, for every $\hat\nu\in \hat \calv$, the corresponding probability
$\hat\P^{\hat\nu}$.

We also consider the approximating problems with
state equation
$$
d\hat X_t^N =  b_t^N( \hat X^N,\hat I)\,dt + \sigma^N_t(\hat X^N,\hat I)\,dB_t,
$$
 and initial condition $X^N_0=x_0$, corresponding to the truncated coefficients
$b^N$, $\sigma^N$, and with the gain functional
\beq \label{defJrandomized_N}
J^{\Rc,N}(\hat\nu) &=&  \hat\E^{\hat\nu}
\left[\int_0^Tf_t(\hat X^N,\hat I)\,dt+g(\hat X^N)\right]
\enq
and the corresponding value
\begin{equation}\label{dualvalue_N}
{\text{\Large$\upsilon$}}_0^{\Rc,N} \;=\;   \sup_{\hat \nu\in\hat \calv} J^{\Rc,N}(\hat \nu).
\end{equation}

Since
 $\hat B$ remains a Brownian motion under each probability $\hat\P^{\hat\nu}$
  with respect to $\F^{\hat x_0,\hat B,\hat \mu}$,
the same arguments that led to \eqref{stimaunif_N} also yield the inequality
\begin{equation}\label{stimaunif_N_hat}
\hat\E^{\hat\nu} \,\Big[\sup_{t\in [0,T]}|\hat X_t^{N}|^p\Big]   \le
    (1+\hat \E\,|x_0|^p)\, C(p,T,M)
\end{equation}
with the same constant
  $C(p,T,M)$, noting also that $\hat \E^{\hat\nu}\,|x_0|^p=\hat \E\,|x_0|^p$
  since $x_0$ has the same law under $\hat\P^{\hat\nu} $ and $\hat\P$.
Setting
$$
\tau^N=\inf\{t\in [0,T]\,:\, |\hat X_t^{N}|\ge N\},
$$
with the convention that $\inf \emptyset =T$, and arguing as before,
we conclude that
$\hat\P$-a.s. (or equivalently $\hat\P^{\hat\nu} $-a.s.)
we have $\tau^N=T$ for all but a finite number of integers $N$, that
we can   define  a unique process
$(\hat X_t)_{0\leq t \leq T}$   such that $\hat X =\hat X^{N}$
on $[0,\tau^N]$ for every $N$, that   $\hat X$
is the unique strong solution to the state equation \eqref{dynXrandom}, and that it
satisfies
\begin{equation}\label{stimaunifrandomizz}
\hat\E^{\hat\nu} \,\Big[\sup_{t\in [0,T]}|\hat X_t |^p\Big]
    \le
    (1+\hat \E\,|x_0|^p)\, C(p,T,M)
\end{equation}
for every $\hat\nu\in \hat \calv$
 and the gain functional
\eqref{defJrandomized} and the value  \eqref{dualvalue}
are well defined and finite.

We are ready to state and prove the required generalization of
Theorem
\ref{MainThm}. The reader will notice that,
instead of the inequality $p\ge \max(2,2r)$  in {\bf (A1)}-(v), here
 it is enough to suppose that  $p\ge 2$ and $p>r$.

\begin{Theorem}
\label{MainThm_bis}
Under {\bf (A2)} and {\bf (A3)} we have
${\text{\Large$\upsilon$}}_0 = {\text{\Large$\upsilon$}}_0^\Rc$
and this common value only depends on the objects
$A,b,\sigma,f,g,\rho_0$ appearing in assumption {\bf (A3)}.
\end{Theorem}
{\bf Proof.} Since $b^N,\sigma^N$ satisfy assumption {\bf (A1)}, by
Theorem
\ref{MainThm} we have
\begin{equation}\label{uguaglicon_N}
    \inf_{\alpha\in\Ac^W} J^N(\alpha)={\text{\Large$\upsilon$}}^N_0=
{\text{\Large$\upsilon$}}_0^{\Rc,N} \;=\;   \sup_{\hat \nu\in\hat \calv} J^{\Rc,N}(\hat \nu).
\end{equation}

We claim that $  J^N(\alpha)\to  J(\alpha)$ uniformly with respect to $\alpha\in\Ac^W$
and that $J^{\Rc,N}(\hat \nu)\to J^{\Rc}(\hat \nu)$  uniformly with respect to
$\hat \nu\in\hat \calv$. This allows to pass to the limit in \eqref{uguaglicon_N}
as $N\to\infty$ and to obtain the required equality
${\text{\Large$\upsilon$}}_0=
{\text{\Large$\upsilon$}}_0^{\Rc}$.

In order to prove the claim we note that
$$
 J^N(\alpha)-J(\alpha)  =  \E\bigg[
\int_0^T\big\{f_t(X^{\alpha,N},\alpha)-f_t(X^{\alpha},\alpha)\big\}\,
1_{\tau_\alpha^N<T}\,dt
+\big\{g(X^{\alpha,N})- g(X^{\alpha,N})\big\}  \, 1_{\tau_\alpha^N<T}\bigg]
$$
so that by \eqref{PolGrowth_f_g_bis}
and by the H\"older inequality with conjugate exponents $p/r$ and $p/(p-r)$,
\begin{eqnarray*}
|J^N(\alpha)-J(\alpha)|
&\le & \E\left[
K\big(2+ \sup_{t\in [0,T]}|X^{\alpha,N}_t|^r + \sup_{t\in [0,T]}|X^{\alpha}_t|^r \big)
(1+T)\, 1_{\tau_\alpha^N<T}\right]
\\
&\le& c\,
\left\{
1+ \E\, \big[\sup_{t\in [0,T]}|X^{\alpha,N}_t|^p\big]
+\E\, \big[ \sup_{t\in [0,T]}|X^{\alpha}_t|^p \big]
\right\}^\frac{r}{p}
 \, \left\{P(\tau_\alpha^N<T)\right\}^\frac{p-r}{p},
\end{eqnarray*}
for some constant $c$ only depending on $K,T,p,r$. By \eqref{stimaunif_N}
and \eqref{stimaunifinalpha} the first term in curly brackets is bounded by a
constant independent of $\alpha$ and $N$. By the Markov inequality
and \eqref{stimaunif_N} we have
$$
P(\tau_\alpha^N<T)\le
P(\sup_{t\in [0,T]}|X_t^{\alpha,N}| \ge N)\le N^{-p}\,
    \E\,\Big[\sup_{t\in [0,T]}|X_t^{\alpha,N}|^p\Big]   \le  N^{-p}\,
    (1+\E\,|x_0|^p)\, C(p,T,M)
$$
which tends to $0$ as $N\to\infty$, uniformly with respect to $\alpha\in\Ac^W$.
This proves that $  J^N(\alpha)\to  J(\alpha)$ uniformly with respect to $\alpha\in\Ac^W$,
and the proof that
 $J^{\Rc,N}(\hat \nu)\to J^{\Rc}(\hat \nu)$  uniformly with respect to
$\hat \nu\in\hat \calv$ is obtained by similar arguments.
\qed

\section{The randomized equation}
\label{Sec:separandom}

In this section the assumptions {\bf (A2)} and {\bf (A3)}
are assumed to hold. We will show how the randomized formulation of the control problem leads to a randomized equation in terms of a backward SDE. We choose
a setting for the randomized control problem \reff{dualvalue} as in Remark \ref{remchoice}, i.e. a product extension $(\bar \Omega, \bar \calf,\bar \P,   V,  W,   \mu, x_0)$
of the  setting $(\Omega,\Fc,\P,V,W,x_0)$ for the initial control problem \reff{primalvalue}.
In view of Proposition \ref{indepofthesetting},  entirely analogous results hold true
in any setting $(\hat \Omega, \hat \calf,\hat \P, \hat V,\hat W, \hat \mu, \hat x_0)$
for the randomized control problem as described in section
\ref{randomizedformulation}.

\subsection{The separated randomized control problem}

We first consider the (path-dependent) filtering of the randomized process $X$ solution to \reff{dynXrandom}, which consists in the process of conditional distributions $\rho_t$ of $X_{\cdot\wedge t}$ given ${\cal F}_t^{W,\mu}$. More precisely, let $\mathscr{P}(\bfC_n)$ be the space of probability measures on $\bfC_n$ and let $B_b(\bfC_n)$ denote the space of bounded Borel measurable real functions on $\bfC_n$. We define $\rho=(\rho_t)_{0\leq t\leq T}$ as an $\F^{W,\mu}$-optional process valued in $\mathscr{P}(\bfC_n)$ satisfying, for every $\varphi\in B_b(\bfC_n)$, (we use the notation $\rho_t(\varphi)=\int_{\bfC_n}\varphi(x)\,\rho_t(dx)$)
\beq \label{defrho}
\rho_t (\varphi) &=& \bar \E\big[ \varphi(X_{\cdot\wedge t}) \mid \Fc_t^{W,\mu} \big], \qquad t\in [0,T],\,\bar\P\text{-a.s.}
\enq
The process
$t\mapsto\bar\E[\varphi(X_{\cdot\wedge t})\mid \calf^{W,\mu}_t]$
is understood as an optional projection. The existence of such a process $\rho$ follows for example from Theorem 2.24 in \cite{BaCri}. While \reff{defrho} is defined for bounded $\varphi$,  since $\rho_t$ is constructed as a $\mathscr{P}(\bfC_n)$-valued process, relation \reff{defrho}
holds for unbounded $\varphi$ once the conditional expectation is well-defined, i.e. $\rho_t (|\varphi|)$ $<$ $\infty$ for all $t\in[0,T]$, $\bar\P$-a.s. (see e.g. Remark 2.27 in \cite{BaCri}).

We can now express the randomized gain functional in terms of the $\F^{W,\mu}$-optional processes $\rho$ and $I$.

\begin{Lemma}
\label{L:Jrho}
For any $\nu$ $\in$ $\Vc$, we have
\beqs
J^{\Rc}(\nu) &=& \E^\nu \Big[ \int_0^T \rho_t( f_t(\cdot,I))\,dt + \rho_T(g) \Big]
\enqs
and, more generally,
\beq
\label{Separation}
\E^\nu\Big[\int_t^Tf_s(X,I)\,ds+g(X)\big|\Fc_t^{W,\mu}\Big] &=& \E^\nu \Big[ \int_t^T \rho_s( f_s(\cdot,I))\,ds + \rho_T(g) \big| \Fc_t^{W,\mu} \Big],
\enq
for all $0\leq t\leq T$.
\end{Lemma}
{\bf Proof.} The result
follows
from the Bayes formula and the $(\bar\P,\F^{W,\mu})$-martingale property of the density process $\kappa^\nu$.
\ep

\vspace{3mm}

The above Lemma \ref{L:Jrho} together with Theorem \ref{MainThm} proves that the randomized control pro\-blem, and thus the primal control problem under partial observation, can be written in a \emph{separated} form involving $\F^{W,\mu}$-optional state processes:
\beq \label{separ}
{\text{\Large$\upsilon$}}_0 &=& \sup_{\nu\in\Vc} \E^\nu \Big[ \int_0^T \rho_t( f_t(\cdot,I)) dt + \rho_T(g) \Big].
\enq

\subsection{BSDE representation}

The purpose of this paragraph is to show that the separated randomized control problem, described by the right-hand side of \reff{separ}, admits a dual representation in terms of a constrained backward SDE,  which then characterizes both the primal control problem and the randomized control problem (as well as the separated randomized control problem). We shall refer to it as the \emph{randomized equation}.

On the space $(\bar \Omega, \bar \calf,\bar \P)$ equipped with the filtration $\F^{W,\mu}$, let us consider the following constrained BSDE
on the time interval $[0,T]$:
\begin{equation}\label{BSDEconstrained}
\begin{cases}
\vspace{2mm} \dis Y_t \ = \ \rho_T(g)  + \int_t^T\rho_s(f_s(\cdot,I)) ds + K_T - K_t - \int_t^TZ_s\,dW_s - \int_t^T\!\int_A U_s(a)\,\mu(ds\,da), \\
\dis U_t(a) \ \le \ 0.
\end{cases}
\end{equation}

We look for a (minimal) solution to \eqref{BSDEconstrained}  in the sense of the following definition.

\begin{Definition}\label{BSDEdef}
A quadruple $(Y_t,Z_t,U_t(a),K_t)$ $($$t\in [0,T]$, $a\in A$$)$
is called a solution to the BSDE  \eqref{BSDEconstrained} if
\begin{enumerate}
  \item $Y$ $\in$  $\Sc^2(\F^{W,\mu})$, the set of real-valued c\`adl\`ag  $\F^{W,\mu}$-adapted processes  satisfying  $\|Y\|_{\Sc^2}^2$ $:=$ $\bar\E[\sup_{0\leq t\leq T}|Y_t|^2]$ $<$ $\infty$;
  \item $Z$ $\in$ $L_W^2(\F^{W,\mu})$, the set of     $\F^{W,\mu}$-predictable  processes
  with values in $\R^d$
  satisfying $\|Z\|_{L_W^2}^2$ $:=$
      $\bar\E\big[\int_0^T|Z_t|^2dt\big]<\infty$;
  \item $U$ $\in$ $L_{\tilde\mu}^2(\F^{W,\mu})$, the set of real-valued $\calp(\F^{W,\mu})\otimes \calb(A)$-measurable processes  satisfying $\|U\|^2_{L_{\tilde\mu}^2}$ $:=$
  $\bar\E\big[\int_0^T\int_A|U_t(a)|^2\lambda(da)dt\big]$ $<$ $\infty$;
  \item $K$ $\in$ $\Kc^2(\F^{W,\mu})$, the subset of $\Sc^2(\F^{W,\mu})$ consisting of  $\F^{W,\mu}$-predictable nondecreasing process with $K_0=0$;
  \item $\bar\P$-a.s. the equality in \eqref{BSDEconstrained}
holds for every $t\in[0,T]$
  and the constraint $U_t(a)\le 0$ is understood to hold
  $\bar\P(d\bar\omega)\lambda(da)dt$-almost everywhere.
\end{enumerate}
A minimal solution $(Y,Z,U,K)$ is a solution to \eqref{BSDEconstrained} such that for any other solution $(Y',Z'$, $U',K')$, we have
$\bar\P$-a.s., $Y_t\le Y'_t$ for all $t\in [0,T]$.
\end{Definition}

We now state the main result of this section.

\begin{Theorem}
\label{Thm:RandomizedFormula}
There exists a unique minimal solution $(Y,Z,U,K)$ $\in$ $\Sc^2(\F^{W,\mu})\times L_W^2(\F^{W,\mu})\times L_{\tilde\mu}^2(\F^{W,\mu})\times\Kc^2(\F^{W,\mu})$ to the randomized equation \eqref{BSDEconstrained}. Moreover, we have $Y_0=\sup_{\nu\in\calv} J^\Rc(\nu)$, and, more generally,
\beq
\label{RandomizedFormula}
Y_t &=& \esssup_{\nu\in\calv} \E^\nu\bigg[ \int_t^T\rho_s (f_s(\cdot,I))\,ds + \rho_T(g) \,\bigg|\, \calf_t^{W,\mu}\bigg].
\enq
\end{Theorem}
\begin{Remark}
{\rm
Combining Theorems \ref{MainThm} and \ref{Thm:RandomizedFormula} we deduce the BSDE representation for the primal problem
\[
Y_0 \ = \ \sup_{\alpha\in\Ac^W} J(\alpha).
\]
We refer sometimes to $Y_0=\sup_{\nu\in\calv} J^\Rc(\nu)$ as \emph{duality relation}, since $Y_0$ coincides with the infimum $\inf\{Y_0'\colon$$(Y',Z',U',K')$ $\in$ $\Sc^2(\F^{W,\mu})\times L_W^2(\F^{W,\mu})\times L_{\tilde\mu}(\F^{W,\mu})\times\Kc^2(\F^{W,\mu})$ solution to \eqref{BSDEconstrained}$\}$.
\ep
}
\end{Remark}
\textbf{Proof (of Theorem \ref{Thm:RandomizedFormula})}
Let us introduce for every $n\in\N$ the following penalized BSDE on $[0,T]$:
\begin{equation}\label{BSDEpenalized}
Y_t^n \ = \ \rho_T(g)  + \int_t^T\rho_s(f_s(\cdot,I))\,ds + K_T^n - K_t^n - \int_t^TZ_s^n\,dW_s - \int_t^T\int_A U_s^n(a)\,\mu(ds\,da),
\end{equation}
where
\[
K_t^n \ = \ n \int_0^t \int_A \big(U_s^n(a)\big)^+ \,\lambda(da)\,ds.
\]
Set $\xi:=\rho_T(g)$ and $F_t:=\rho_t(f_t(\cdot,I))$. By \eqref{PolGrowth_f_g} and
\eqref{stimaunifrandomizz}
 we see that
\[
\bar\E|\xi|^2 \ < \ \infty, \qquad\qquad \bar\E\bigg[\int_0^T |F_t|^2\,dt\bigg] \ < \ \infty
\]
(here we use the assumption that $p\ge 2r$  in {\bf (A1)}-(v)).
Then, from Lemma 2.4 in \cite{tang_li} it follows that, for every $n\in\N$, there exists a unique solution $(Y^n,Z^n,U^n)$ $\in$ $\Sc^2(\F^{W,\mu})\times L_W^2(\F^{W,\mu})\times L_{\tilde\mu}^2(\F^{W,\mu})$ to the above penalized BSDE.

Now, proceeding along the same lines as in the proof of Lemma 4.8 in \cite{FP15}, we obtain the formula
\beq
\label{RandomizedFormula_n}
Y_t^n &=& \esssup_{\nu\in\calv_n} \E^\nu\bigg[ \int_t^T\rho_s (f_s(\cdot,I))\,ds + \rho_T(g) \,\bigg|\, \calf_t^{W,\mu}\bigg],
\enq
where $\Vc_n=\{\nu\in\Vc\colon\nu\text{ takes values in }(0,n]\}$. By \eqref{Separation}, together with estimates \eqref{PolGrowth_f_g} and \eqref{EstimateX_nu}, we deduce that
\begin{equation}
\label{YnUpperBound}
\sup_n Y_t^n \ < \ \infty, \qquad \text{for all }0\leq t\leq T.
\end{equation}
Notice that equation \eqref{BSDEconstrained} can be written as follows:
\begin{equation}\label{BSDEconstrained2}
\begin{cases}
\vspace{2mm} \dis Y_t \ = \ \rho_T(g)  + \int_t^T \bigg(\rho_s(f_s(\cdot,I)) - \int_A U_s(a)\,\lambda(da)\bigg) ds + K_T - K_t \\
\vspace{2mm} \dis \qquad\;\;\; - \int_t^TZ_s\,dW_s - \int_t^T\!\int_A U_s(a)\,\tilde\mu(ds\,da), \\
\dis U_t(a) \ \le \ 0.
\end{cases}
\end{equation}
Then, we see that the above equation is a particular case of a backward stochastic differential equation studied in a general non-Markovian framework in \cite{KP12}. In particular, existence and uniqueness of the minimal solution $(Y,Z,U,K)$ to equation \eqref{BSDEconstrained2} (or, equivalently, to equation \eqref{BSDEconstrained}) follow from Theorem 2.1 in \cite{KP12}. Indeed, Assumption {\bf (H0)} in \cite{KP12} is clearly satisfied. Concerning Assumption {\bf (H1)}, this is only used in Lemma 2.2 of \cite{KP12} to prove that the sequence $(Y^n)_n$ satisfies \eqref{YnUpperBound}, a property that in our setting follows directly from \eqref{PolGrowth_f_g} and \eqref{EstimateX_nu}.

Finally, from Theorem 2.1 in \cite{KP12} we also have that $Y_t^n(\bar\omega)$ converges increasingly to $Y_t(\bar\omega)$ as $n\rightarrow\infty$, $\bar\P(d\bar\omega)$-a.s. Since $\Vc=\cup_n\Vc_n$, letting $n\rightarrow\infty$ in \eqref{RandomizedFormula_n} we obtain \eqref{RandomizedFormula}.
\qed

\vspace{3mm}

We end this section proving the following generalization of formula \eqref{RandomizedFormula}.

\begin{Theorem}[Randomized dynamic programming principle]\label{T:DPPRandomizedProblem}
\quad\\
For all $0\leq t\leq T$, we have
\beq\label{DPPRandomizedProblem}
Y_t &=& \esssup_{\nu\in\calv} \esssup_{\tau\in\Tc} \E^\nu\bigg[ \int_t^\tau\rho_r (f_r(\cdot,I))\,dr + Y_\tau \,\bigg|\, \calf_t^{W,\mu}\bigg] \notag \\
&=& \esssup_{\nu\in\calv} \essinf_{\tau\in\Tc} \E^\nu\bigg[ \int_t^\tau\rho_r (f_r(\cdot,I))\,dr + Y_\tau \,\bigg|\, \calf_t^{W,\mu}\bigg],
\enq
where $\mathcal T$ denotes the class of $[0,T]$-valued $\F^{W,\mu}$-stopping times.
\end{Theorem}
\textbf{Proof.}
For every $n$, proceeding along the same lines as in the proof of Lemma 4.8 in \cite{FP15}, we obtain
\begin{align}
Y_t^n \ &= \ \E^\nu\bigg[Y_\tau^n + \int_t^\tau\rho_r(f_r(\cdot,I))\,dr + \int_t^\tau\int_A \big[n\big(U_r^n(a)\big)^+ - \nu_r(a) U_r^n(a)\big]\lambda(da)dr\,\bigg|\, \calf_t^{W,\mu}\bigg] \label{DPPRandomizedProblem4} \\
&\geq \ \E^\nu\bigg[Y_\tau^n + \int_t^\tau\rho_r(f_r(\cdot,I))\,dr\,\bigg|\, \calf_t^{W,\mu}\bigg], \qquad \text{for all }\nu\in\Vc_n,\,\tau\in\Tc. \notag
\end{align}
Recalling that $Y\geq Y^n$, we find
\[
Y_t \ \geq \ \E^\nu\bigg[Y_\tau^n + \int_t^\tau\rho_r(f_r(\cdot,I))\,dr\,\bigg|\, \calf_t^{W,\mu}\bigg], \qquad \text{for all }\nu\in\Vc_n,\,\tau\in\Tc.
\]
Letting $n\rightarrow\infty$, and recalling that $\Vc_n\subset\Vc$, we end up with
\[
Y_t \ \geq \ \E^\nu\bigg[Y_\tau + \int_t^\tau\rho_r(f_r(\cdot,I))\,dr\,\bigg|\, \calf_t^{W,\mu}\bigg], \qquad \text{for all }\nu\in\Vc,\,\tau\in\Tc.
\]
The above inequality yields
\begin{align*}
Y_t \ \geq \ \esssup_{\nu\in\calv} \esssup_{\tau\in\Tc} \E^\nu\bigg[ \int_t^\tau\rho_r (f_r(\cdot,I))\,dr + Y_\tau \,\bigg|\, \calf_t^{W,\mu}\bigg], \\
Y_t \ \geq \ \esssup_{\nu\in\calv} \essinf_{\tau\in\Tc} \E^\nu\bigg[ \int_t^\tau\rho_r (f_r(\cdot,I))\,dr + Y_\tau \,\bigg|\, \calf_t^{W,\mu}\bigg].
\end{align*}
It remains to prove the reverse inequalities. As in the proof of Lemma 4.8 in \cite{FP15}, for every $n$ and $\eps\in(0,1)$, we define $\nu_r^{\eps,n}(a)=n\,1_{\{U_r^n(a)\geq0\}} + \eps\,1_{\{-1<U_r^n(a)<0\}} - \eps\,(U_r^n(a))^{-1}\,1_{\{U_r^n(a)\leq1\}}$. Then, $\nu^{\eps,n}\in\Vc_n$ and
\[
n\big(U_r^n(a)\big)^+ - \nu_r^{\eps,n}(a) U_r^n(a) \ \leq \ \eps, \qquad \text{for all }r\in[0,T].
\]
Therefore, from equality \eqref{DPPRandomizedProblem4}, we find
\begin{align*}
Y_t^n \ &= \ \E^{\nu^{\eps,n}}\bigg[Y_\tau^n + \int_t^\tau\rho_r(f_r(\cdot,I))\,dr + \int_t^\tau\int_A \big[n\big(U_r^n(a)\big)^+ - \nu_r(a) U_r^n(a)\big]\lambda(da)dr\,\bigg|\, \calf_t^{W,\mu}\bigg] \\
&\leq \ \E^{\nu^{\eps,n}}\bigg[Y_\tau^n + \int_t^\tau\rho_r(f_r(\cdot,I))\,dr\,\bigg|\, \calf_t^{W,\mu}\bigg] + \eps\,\lambda(A)\,T, \qquad \text{for all }\tau\in\Tc.
\end{align*}
Then, we obtain the two following inequalities:
\begin{align*}
Y_t^n \ &\leq \ \esssup_{\tau\in\Tc}\E^{\nu^{\eps,n}}\bigg[Y_\tau^n + \int_t^\tau\rho_r(f_r(\cdot,I))\,dr\,\bigg|\, \calf_t^{W,\mu}\bigg] + \eps\,\lambda(A)\,T, \\
Y_t^n \ &\leq \ \essinf_{\tau\in\Tc}\E^{\nu^{\eps,n}}\bigg[Y_\tau^n + \int_t^\tau\rho_r(f_r(\cdot,I))\,dr\,\bigg|\, \calf_t^{W,\mu}\bigg] + \eps\,\lambda(A)\,T.
\end{align*}
As a consequence, we get (we continue the proof with ``$\text{ess\,inf}$'' over $\tau\in\Tc$, since the proof with ``$\text{ess\,sup}$'' can be done proceeding along the same lines)
\[
Y_t^n \ \leq \ \esssup_{\nu\in\calv_n}\essinf_{\tau\in\Tc}\E^\nu\bigg[Y_\tau^n + \int_t^\tau\rho_r(f_r(\cdot,I))\,dr\,\bigg|\, \calf_t^{W,\mu}\bigg] + \eps\,\lambda(A)\,T.
\]
Using the arbitrariness of $\eps$, and recalling that $\Vc_n\subset\Vc$ and $Y^n\leq Y$, we obtain
\[
Y_t^n \ \leq \ \esssup_{\nu\in\calv}\essinf_{\tau\in\Tc}\E^\nu\bigg[Y_\tau + \int_t^\tau\rho_r(f_r(\cdot,I))\,dr\,\bigg|\, \calf_t^{W,\mu}\bigg].
\]
The claim follows letting $n\rightarrow\infty$.
\ep

\section{Numerical implications}
\label{Sec:numerics}

We discuss briefly in this section how the BSDE representation \reff{BSDEconstrained} can provide a new probabilistic numerical scheme for solving
partial observation control problem.  We postpone the detailed study for a future research work.
We shall  consider the case without path-dependence in the state and control, i.e. $b(t,x,a)$, $\sigma(t,x,a)$, $f(t,x,a)$ are functions
on $[0,T]\times\R^n\times A$, and $g(x)$ is a function defined on $\R^n$,  and we  focus our attention on the partial observation feature in this  Markov setting.
In this case, the forward  component of the BSDE \reff{BSDEconstrained} is given by the so-called randomized filter process
\beqs
\rho_t(\varphi) & =& \bar \E \big[ \varphi(X_t) | \Fc_t^{W,\mu} \big],  \;\;\; t \in [0,T],
\enqs
for  any $\varphi$ $\in$ $B_b(\R^n)$, the set of bounded Borel measurable functions on $\R^n$, and where $(X,I)$ is the randomized regime switching diffusion process in \reff{dynXrandom}, hence given by
\beqs
X_t &=& X_0 + \int_0^t b(s,X_s,I_s) ds + \int_0^t \sigma(s,X_s,I_s) dB_s \\
I_t &=& I_0 +  \int_{(0,t]} \int_A (a-I_{s^-}) \mu(ds,da).
\enqs

The first step is to get an approximation for the randomized filter process: the exact filter is in general intractable in infinite dimension and numerical
approximation techniques are needed. Several numerical techniques have been developed for reducing the approximation of the filter to a random discrete measure with $M$ points, by particle  or quantization methods, and  we refer to \cite{BaCri} for an overview of these approximations.
Anyway, this  leads to an approximation of the filter process with $\bar\rho^\pi$ $=$ $\{\bar\rho_{t_k}^\pi, k =0,\ldots,N\}$ valued in $K^M$, the simplex of
$\R^M$, for a partition $\pi$ $=$ $\{t_0=0< \ldots < t_k < \ldots t_N=T\}$ of the time interval $[0,T]$.
Next, the numerical issue is to deal with the nonpositive jump constraint in the BSDE representation \reff{BSDEconstrained}, and following \cite{KLP15} (originally  designed  for
full observation control problem),  we propose a discrete time approximation scheme of the form:
\begin{equation}\label{schemeBSDE}
\left\{ \begin{array}{rcl}
\bar Y_T^\pi \; = \; \bar \Yc_T^\pi &=& \bar\rho_T^\pi(g)  \\
\bar\Yc_{t_k}^\pi &=& \E\Big[ \bar Y_{t_{k+1}}^\pi  \big| \Fc_{t_k}^{W,\mu} \Big] +  (t_{k+1}-t_k) \; \bar\rho_{t_k}^\pi\big( f(t_k,.,I_{t_k})\big) \\
\bar Y_{t_k}^\pi &=& \esssup_{a \in A} \E \Big[ \bar\Yc_{t_k}^\pi \big| \Fc_{t_k}^{W,\mu},  I_{t_k} = a \Big], \;\;\; k = 0, \ldots, N-1,
\end{array}
\right.
\end{equation}
The interpretation of this scheme is the following. The  first two lines in \reff{schemeBSDE} correspond to the standard scheme $\bar \Yc^\pi$ for a discretization of a BSDE with jumps, where we omit here the computation of the Brownian and jump component since they do not appear in the generator of the BSDE.  The last line in \reff{schemeBSDE} for computing the approximation $\bar Y^\pi$ of the minimal solution $Y$  corresponds precisely to the minimality condition for the nonpositive jump constraint and should be understood as follows. By the Markov property of the forward process $(\rho_t,I_t)$, the solution $\Yc$  to the BSDE with jumps (without constraint)  is in the form $\Yc_t$ $=$
$\vartheta(t,\rho_t,I_t)$ for some deterministic function $\vartheta$.  Assuming that $\vartheta$ is a continuous function, the jump component of the BSDE, which is induced by a jump of the forward component $I$,
is equal to  $\Uc_t(a)$ $=$ $\vartheta(t,\rho_t,a)-\vartheta(t,X_t,I_{t^-})$.  Therefore, the nonpositive jump constraint  means that:
$\vartheta(t,\rho_t,I_{t^-})$ $\geq$ $\esssup_{a\in A} \vartheta(t,\rho_t,a)$. The  minimality condition is thus written as:
\beqs
Y_t \; = \;  v(t,\rho_t) \; = \;  \esssup_{a \in A} \vartheta(t,\rho_t,a) &=& \esssup_{a \in A} \E [ \Yc_t | \rho_t, I_{t} = a ],
\enqs
whose discrete time version is the last line in scheme \reff{schemeBSDE}.  The practical implementation of the above discrete time scheme
requires the estimation and computation of the conditional expectations together with the supremum. This can be achieved  for example with
regression methods on  basis functions defined on $K^M\times A$, based on simulation of the approximate randomized filter  $\bar\rho^\pi$ together with the pure jump process $I$.

\appendix

\section{Appendix}

\setcounter{Theorem}{0}
\setcounter{Definition}{0}
\setcounter{Proposition}{0}
\setcounter{Assumption}{0}
\setcounter{Lemma}{0}
\setcounter{Corollary}{0}
\setcounter{Remark}{0}
\setcounter{Example}{0}
\setcounter{equation}{0}

This section is devoted to the proof of Proposition \ref{extensionapproximation} below, which was used in the proof of
Theorem \ref{MainThm}. We assume that $A$ is a Borel space,
and that $\lambda$ and $a_0$
 are given and satisfy the assumption {\bf (A2)}. Our starting point is also a probability  space $(\Omega,\calf,\P)$, with a filtration
 $\G=(\calg_t)_{t\ge 0}$.

Following \cite{80Krylov},
for any pair  $\alpha^1,\alpha^2:\Omega\times [0,T]\to A$ of   $\G$-progressive
  processes we define a distance $\tilde\rho(\alpha^1,\alpha^2)$ setting
\beqs
\tilde \rho(\alpha^1,\alpha^2) &=& \E \Big[\int_0^T\rho(\alpha^1_t,\alpha^2_t)\,dt \Big].
\enqs
where $\rho$ is an arbitrary metric in $A$ satisfying $\rho<1$.

Below we will use an auxiliary
probability space denoted
 $(\Omega',\calf',\P')$. This can be taken as an arbitrary  probability space
 where appropriate random objects  are defined.
For integers $m,n,k\ge1$,  we assume that real random variables  $ U_n^m $, $V^m_n$
and random measures
 $\pi^k$ are defined on  $(\Omega',\calf',\P')$ and satisfy the following conditions:
\begin{enumerate}
\item every $ U_n^m $ is uniformly  distributed on $(0,1)$;
\item  each   $ V_n^m$ has exponential distribution with parameter $\lambda_{nm}$
and $\sum_{n=1}^\infty\lambda_{nm}^{-1}=1/m$ for every $m\ge 1$;
\item  every $\pi^k$ is a Poisson random measure on $(0,\infty)\times A$, admitting compensator $k^{-1}\lambda(da)\,dt$ with respect to its natural filtration;
 \item the random elements  $U^m_n,V^h_j$, $\pi^k$    are all independent.
 \end{enumerate}

The role of these random elements   will become clear in the constructions that
follow. Notice that for the construction of the space
$(\Omega',\calf',\P')$ only the knowledge of the measure $\lambda$ is required.
Morevoer  by a classical result,
 see \cite{Zabczyk96} Theorem 2.3.1, we may take $\Omega'=[0,1]$,
 $\calf'$ the corresponding Borel sets and $\P'$ the Lebesgue measure.
Next we define
\beqs
\hat\Omega \; = \; \Omega\times \Omega',
\qquad
\hat \calf \; = \; \calf\otimes \calf',
\qquad
\Q \; = \; \P\otimes \P'
\enqs
and note that
the filtration $\G$ can be canonically extended to a filtration $\hat\G=(\hat\calg_t)_{t\geq 0}$
in $(\hat\Omega,\hat\calf)$ setting $\hat\calg_t=\{A\times \Omega'\,:\, A\in\calg_t\}$.
Similarly, any process $\alpha$ in $(\Omega,\calf)$ admits an extension $\hat\alpha$
to $(\hat\Omega,\hat\calf)$ given by $\hat\alpha_t(\hat\omega)=\alpha_t(\omega)$,
where $\hat\omega=(\omega,\omega')$.
The metric $\tilde\rho$ can also be extended to
  any pair  $\beta^1,\beta^2:\hat\Omega\times [0,T]\to A$ of   $\hat\G$-progressive
  processes  setting
\beqs
\tilde \rho(\beta^1,\beta^2) &=& \E^\Q \Big[\int_0^T\rho(\beta^1_t,\beta^2_t)\,dt \Big].
\enqs
We use the same symbol $\tilde\rho$ to denote  the extended metric as well.

\vspace{1mm}

Our aim in this section is to prove the following result.

\begin{Proposition}
\label{extensionapproximation}
Let  $A$ be a Borel space, and let $\lambda$ and $a_0$
 satisfy   {\bf (A2)}.
Let  $(\Omega,\calf,\P)$ be  any probability
 space with a filtration
 $\G=(\calg_t)_{t\ge 0}$ and let
$(\hat\Omega,
\hat \calf,\Q)$ be the product space defined above. Then
 for any $\G$-progressive  $A$-valued
 process $\alpha$, and for any $\delta >0$,
 there exists a marked point process $(\hat S_n,\hat \eta_n)_{n\ge 1}$ defined in
$(\hat\Omega,
\hat \calf,\Q)$
satisfying the following conditions:
\begin{enumerate}
\item
setting
$$
\hat S_0=0,\qquad\hat \eta_0=a_0,
\qquad
\hat I_t=\sum_{n\ge 0}\hat  \eta_{n}1_{ [\hat S_n,\hat S_{n+1})}(t),
$$
 the process
$\hat I$
satisfies
\begin{equation}\label{distkrylovdelta}
    \tilde\rho (\hat I, \hat \alpha)=
\E_\Q\left[\int_0^T  \rho (\hat I_t,\hat \alpha_t)\,dt\right]<\delta;
\end{equation}
\item denoting
$\hat \mu=\sum_{n\ge1}\delta_{(\hat S_n,\hat \eta_n)}$ the random measure
associated to $(\hat S_n,\hat \eta_n)_{n\ge 1}$,
$\F^{\hat\mu}=(\calf_t^{\hat\mu})_{t\geq 0}$
 the natural filtration of
$\hat\mu$ and $\hat\G\vee \F^{\hat\mu}=(\hat\calg_t\vee\calf_t^{\hat\mu})_{t\geq 0}$,
then
the $\hat\G\vee \F^{\hat\mu}$-compensator of    $\hat \mu$ under $\Q$ is absolutely continuous
 with respect to $\lambda(da)\,dt$ and
it can be written in the form
\begin{equation}\label{compensatoreBmu}
 \hat\nu_t(\hat\omega,a)\, \lambda(da)\,dt
\end{equation}
 for  some  nonnegative $\calp(\hat\G\vee \F^{\hat\mu})
 \otimes \calb(A)$-measurable  function $\hat\nu$ satisfying
 \begin{equation}\label{compensatoreBmubounds}
 \inf_{\hat\Omega\times [0,T]\times A}\hat\nu>0,
 \qquad
  \sup_{\hat\Omega\times [0,T]\times A}\hat\nu<\infty.
  \end{equation}
\end{enumerate}

\end{Proposition}
\textbf{Proof.}
Fix $\alpha$ and $\delta$ as in the statement of the Proposition.
It can be proved that there exists an $A$-valued process  $\bar\alpha$     such that
$\tilde\rho (\alpha,\bar\alpha)<\delta$ and $\bar\alpha$
has the form $\bar\alpha_t=\sum_{n= 0}^{N-1} \alpha_{n}1_{ [t_n,t_{n+1})}(t)$,
where  $0=t_0<t_1<\ldots t_N=T$ is a deterministic subdivision of $[0,T]$,
$\alpha_0,\ldots,\alpha_{N-1}$ are $A$-valued
random variables that take only a finite number of values, and each $\alpha_n$ is $\calg_{t_n}$-measurable:
this is an immediate consequence of
 Lemma 3.2.6 in \cite{80Krylov}, where  it is proved that the set of admissible
controls $\bar\alpha$ having the form specified in the lemma are dense in
 the set of all  $\G$-progressive $A$-valued processes with respect to the metric $\tilde \rho$.

We can  (and will)  choose $\bar\alpha$ satisfying $\alpha_0=a_0$
($a_0$ is the same as in {\bf (A2)}). Indeed
this additional requirement can be fulfilled by adding, if necessary, another point $t'$
close to $0$ to the subdivision
and modifying $\bar\alpha$ setting $\bar\alpha_t=a_0$  for   $t\in [0,t')$.
This modification is as close as we wish to the original process with respect
to the metric $\tilde \rho$, provided $t'$ is chosen sufficiently small.

Finally, we further extend   $\bar\alpha$ to a function
defined on $\Omega\times [0,\infty)$ in a trivial way setting
$\bar\alpha_t=\sum_{n= 0}^{\infty} \alpha_{n}1_{ [t_n,t_{n+1})}(t)$ where
$\alpha_n =  \alpha_{N-1}$ for $n\ge N$ and $t_n=t+n-N$ for
$n>N$.
This way $\bar\alpha$ is associated to the marked point process $(t_n,\alpha_n)_{n\ge 1}$ and $\bar \alpha_0=a_0$.

\vspace{1mm}

Next recall  the spaces
 $(\Omega',\calf',\P')$ and
 $(\hat\Omega,\hat\calf,\Q)$ and the filtration $\hat\G$
 introduced before the statement
 of Proposition
\ref{extensionapproximation}.
We  extend the processes $\alpha$ and $\bar\alpha$ to $\hat\Omega\times
[0,\infty)$ and denote $\hat\alpha$ and $\hat{\bar{\alpha}}$ the corresponding
extensions. We note that clearly
\beq\label{rhotildeuno}
\tilde \rho(\hat\alpha,\hat{\bar{\alpha}})
 &=& \tilde \rho(\alpha,{\bar{\alpha}})<\delta/3.
\enq

The next step of the proof consists in constructing a sequence of random
measures $\kappa^m$ whose associated piecewise constant
trajectories, denoted $\hat\alpha_t^m$,  approximate $\hat   \alpha$
in the sense of the metric $\tilde\rho$.
The construction will be carried out in such a way that $\kappa^m$
admits a compensator absolutely continuous with respect to the measure
$\lambda(da)\,dt$.

For every   $m\ge 1$, let $\bfB(b,1/m)$ denote the open ball  of radius $1/m$, with respect to the metric $\rho$, centered at $b\in A$.
Since $\lambda(da)$ has full support, we have $\lambda(\bfB(b,1/m))>0$ and we can  define a transition kernel $q^m(b,da)$ in $A$ setting
\beqs
q^m(b,da)&=& \frac{1}{\lambda(\bfB(b,1/m))}\, 1_{\bfB(b,1/m)}(a) \lambda(da).
\enqs
We recall that we require $A$ to be a Borel space,
and we denote by $\calb(A)$ its Borel $\sigma$-algebra.
There exists a Borel measurable function $q^m:A\times [0,1]\to A$
such that for every $b\in A$ the measure $B\mapsto q^m(b,B)$ ($B\in \calb(A)$) is the image of the Lebesgue measure  on $[0,1]$ under the mapping
$u\mapsto q^m(b,u)$; equivalently,
$$
\int_A k(a)\,q^m(b,da) \; = \; \int_0^1k( q^m(b,u))\,du,
$$
for every nonnegative measurable function $k$ on $A$.
Thus, if  $U$ is a random variable defined on some probability space and having uniform law on $[0,1]$ then, for fixed $b\in A$, the $A$-valued random variable
$q^m(b,U)$ has law $q^m(b,da)$. The use of the same symbol $q^m$ should not generate confusion.
The existence of the function $q^m$ (even for a general transition kernel on $A$)
 is well known when $A$ is a
 separable complete metric space,
in particular, when $A$ is the unit interval $[0,1]$,
(see e.g. \cite{Zabczyk96}, Theorem 3.1.1) and the general case reduces to this one,
since it is known that any Borel  space is either finite or countable (with the discrete topology)
or isomorphic, as a measurable space, to the interval  $[0,1]$: see e.g.
\cite{BertsekasShreve78}, Corollary 7.16.1.

\vspace{1mm}

For fixed   $m\ge 1$,
define    $V^m_0=R^m_0=S^m_0=0$ and
 $$
R_n^m=t_n+V^m_1+\ldots +V^m_n,
\quad
S_n^m=t_n+V^m_1+\ldots +V^m_{n-1},
\quad \beta_n^m=q^m(\alpha_n,U^m_n),
\qquad \quad  n\ge 1.
$$

Since we assume $t_n<t_{n+1}$ and since $V^m_n>0$
we see that $(R^m_n,\beta^m_n)_{n\ge1}$
is a marked point process in $A$. Also note that
$R^m_{n-1}<S^m_n<R^m_n$ for $n\geq 1$.
Let
$$\kappa^m=\sum_{n\ge1}\delta_{(R^m_n,\beta^m_n)},
\qquad
\hat\alpha_t^m \; = \; \sum_{n\ge 0} \beta^m_{n}1_{ [R^m_n,R^m_{n+1})}(t),
$$
(with the convention $\beta^m_0=a_0$)
denote the corresponding random measure and the associated
trajectory. We claim that
\begin{equation}\label{compensatorofmodifiedppmbis}
\tilde\rho (\hat{\bar{\alpha}}, \hat\alpha^m)\to 0
\end{equation}
as $m\to\infty$.
Indeed, since $0=t_0<t_1<\ldots t_N=T$ we have
\begin{equation}\label{decompdist}
\tilde\rho (\hat{\bar{\alpha}}, \hat\alpha^m)=
\sum_{n=0}^{N-1}\E_\Q \int_{t_n}^{t_{n+1}} \rho(\hat{\bar{\alpha}}_t, \hat\alpha^m_t)\,dt.
\end{equation}
Note that $t_n<R_n^m$, and whenever $ R_n^m\le t<t_{n+1}<R_{n+1}^m$
we have $\hat{\bar{\alpha}}_t=\alpha_n$, $\hat\alpha^m_t=\beta_n^m$ and so
$\rho(\hat{\bar{\alpha}}_t, \hat\alpha^m_t)=\rho(\alpha_n,\beta^m_n)<1/m$
since, for every $b$ $\in$ $A$ , $q^m(b,da)$ is supported in $\bfB(b,1/m)$.
If $R_n^m<t_{n+1}$
then, recalling that $\rho <1$,
\beqs
  \int_{t_n}^{t_{n+1}} \rho(\hat{\bar{\alpha}}_t, \hat\alpha^m_t)\,dt
  &=&
\int_{t_n}^{R_n^m}   \rho(\hat{\bar{\alpha}}_t, \hat\alpha^m_t)\,dt
+\int_{R_n^m}^{t_{n+1}} \rho(\hat{\bar{\alpha}}_t, \hat\alpha^m_t)\,dt
\\
&\le&
(R_n^m-t_n) + \frac1m (t_{n+1}- R_n^m)
\\
&\le&
V_1^m+\ldots +V_n^m + \frac1m (t_{n+1}- t_n).
\enqs
If $R_n^m\ge t_{n+1}$
then the same inequality still holds since we even have
$$
\int_{t_n}^{t_{n+1}} \rho(\hat{\bar{\alpha}}_t, \hat\alpha^m_t)\,dt
\le
t_{n+1}- t_n \le R_n^m - t_n = V_1^m+\ldots +V_n^m.
$$
Substituting in \eqref{decompdist} and
computing the expectation of the exponential random variables $V_n^m$
we arrive at
$$
\tilde\rho (\hat{\bar{\alpha}}, \hat\alpha^m)=
\sum_{n=0}^{N-1}\left(
\lambda_{1m}^{-1}+\ldots +\lambda_{nm}^{-1} + \frac1m (t_{n+1}- t_n)\right)
\le
\sum_{n=1}^\infty\lambda_{nm}^{-1} +\frac{T}{m}\le \frac1m +\frac{T}{m}
$$
 which proves the claim \eqref{compensatorofmodifiedppmbis}.
From now on we fix  a value of $m$ so large that
\beq\label{rhotildedue}
\tilde\rho (\hat{\bar{\alpha}}, \hat\alpha^m)&<&\delta/3.
\enq

Let
$\F^{\kappa^m}=(\calf^{\kappa^m}_t)$ denote the natural filtration of $\kappa^m$ and
 set
\[
\H^m = (\calh_t^m)_{t\ge 0}=
(\hat\calg_t\vee \calf^{\kappa^m}_t)_{t\ge 0}.
\]
We have the following technical result that describes the compensator
$\tilde\kappa^m$ of $\kappa^m$ with respect to the filtration $\H^m$.

 \begin{Lemma}\label{MPPperturbed}
With the previous assumptions and notations, the compensator of the random measure
$\kappa^m$ with respect to $\H^m$ and $\Q$ is given by the formula
 \beqs
 \tilde\kappa^m(dt,da) &=& \sum_{n\ge 1}  1_{(S_n^m,R^m_n]}(t)\,
 q^m(\alpha_n,da) \lambda_{nm}.
 \enqs
\end{Lemma}
{\bf Proof of Lemma \ref{MPPperturbed}.} To shorten notation, we drop all the sub- and superscripts
$m$ and write $\tilde\kappa$, $S_n$, $R_n$,
$ q$, $ \lambda_{n}$,
$\F^{\kappa}$,
$\H = (\calh_t)=
(\hat\calg_t\vee \calf^{\kappa}_t)$
instead of
$\tilde\kappa^m$, $S_n^m$, $R^m_n$,
$ q^m$, $ \lambda_{nm}$, $\F^{\kappa^m}$, $\H^m $ $=$ $ (\calh_t^m)$ $=$
$(\hat\calg_t\vee \calf^{\kappa^m}_t)$.

Let us first check that $\tilde\kappa(dt,da)$, defined by the  formula above, is an $\H$-predictable random measure. The variables $R_n$ are clearly
$\F^\kappa$-stopping times and hence
$\H$-stopping times and therefore
 $S_n= R_{n-1}+t_n-t_{n-1}$  are also $\H$-stopping times.
 Since
  $\alpha_n$    are $\calf_{t_n}$-measurable and
  $\calf_{t_n}\subset  \calh_{t_n}\subset \calh_{S_n}$, $\alpha_n$
  are also
$\calh_{S_{n}}$-measurable. It follows that for every $C\in\calb(A)$ the process
$
1_{(  S_{n},R_n]}(t)\, q(\alpha_n,C) \lambda_n
$
is $\H$-predictable and finally that  $\tilde\kappa(dt,da)$ is an $\H$-predictable random measure.

To finish the proof  we need now to verify that for every  positive $\calp(\H)\otimes \calb(A)$-measurable
 random field $H_t(\omega,a)$ we have
 $$
 \E \Big[\int_0^\infty\int_AH_t(a)\,\kappa(dt\,da) \Big]  \;= \;
 \E \Big[ \int_0^\infty\int_AH_t(a)\,\tilde\kappa(dt\,da)\Big].
 $$
 Since $\calh_t=\calf_t\vee \calf^\kappa_t$, by a monotone class argument
 it is enough to consider $H$ of the form
 $$H_t(\omega,a)=H_t^1(\omega)H_t^2(\omega)k(a),
 $$
 where $H^1$ is a positive $\hat\G$-predictable random process,
$H^2$ is a positive $\F^\kappa$-predictable random process
 and $k$ is a positive $\calb(A)$-measurable function.
 Since $\F^\kappa$ is the natural filtration of $\kappa$,
 by a known result (see e.g. \cite{ja} Lemma (3.3))
 $H^2$ has the following form:
 \begin{eqnarray*}
    H^2_t &=& b_1(t)1_{(0,R_1]}(t)+
 b_2(\beta_1,R_1,t)1_{( R_1,R_2]}(t)
 \\
 &&+
 b_3(\beta_1,\beta_2 ,R_1,R_2,t)1_{(R_2,R_3]}(t)+ \ldots
 \\
     && +
 b_n(\beta_1,\ldots,\beta_{n-1},R_1,\ldots, R_{n-1},t)1_{(R_{n-1},R_n]}(t)+
 \ldots,
 \end{eqnarray*}
where each $b_n$ is a positive measurable deterministic function of
$2n-1$ real variables.
Since
$$
\E \Big[ \int_0^\infty\int_AH_t(a)\,\kappa(dt\,da) \Big] \; = \;
\E \Big[  \sum_{n\ge 1} H_{R_n}(\beta_n) \Big]
 $$
to prove the thesis  it is enough to check that for every $n\ge1$ we have the equality
$$
\E \big[ H_{R_n}(\beta_n) \big]  =
 \E \Big[ \int_0^\infty \int_A H_t(a)\,  q(\alpha_n,da) \lambda_n\,
 1_{  S_{n}<t\le R_n }\,dt \Big]
 $$
which can also be written
$$
\begin{array}{l}\dis
 \E\big[H^1_{R_n}  b_n(\beta_1,\ldots,\beta_{n-1},R_1,\ldots, R_{n-1},R_n)    k(\beta_n) \big]
 =
 \\\dis
   \E \Big[ \int_0^\infty \int_A H_t^1  b_n(\beta_1,\ldots,\beta_{n-1},R_1,\ldots, R_{n-1},t)k(a)
 q(\alpha_n,da)
\lambda_n\,  1_{  S_{n}<t\le R_n }\,dt \Big].
\end{array}
$$
We use the notation
$$
K_n(t) \; = \; H_t^1\,
 b_n(\beta_1,\ldots,\beta_{n-1},R_1,\ldots, R_{n-1},t)
 $$
 to reduce the last equality to
\begin{equation}\label{thesisrewritten}
    \E\, [K_n(R_n)
    k(\beta_n)]
=
 \E\Big[\int_0^\infty \int_A K_n (t)\,k(a)\,
 q(\alpha_n,da) \lambda_n\,
 1_{  S_{n}<t\le R_n }\,dt\Big].
\end{equation}
By the definition of $R_n$ and $\beta_n$, we have
$\E[K_n(R_n)k(\beta_n)]$ $=$  $\E[K_n(S_n+V_n)   k(q(\alpha_n,U_n))] $.
As noted above,   since $U_n$ has uniform law on $(0,1)$,  the random variable
$q(b,U_n)$ has law $q(b,da)$ on $A$,   for any fixed $b\in A$.
We note that $R_1,\ldots, R_{n-1}, S_n$ are measurable with respect to
$\sigma (V_1,\ldots, V_{n-1})$, that $\beta_1,\ldots, \beta_{n-1}$ are
measurable with respect to
$\hat\calg_\infty\vee\sigma (U_1,\ldots, U_{n-1})$ and therefore that
 the random elements $U_n $,
 $S_n$  and $(\hat\calg_\infty, \beta_1,\ldots,\beta_{n-1},R_1,\ldots, R_{n-1},S_n)$ are all  independent.
 Recalling that $V_n$ is exponentially distributed with parameter $\lambda_n$
 we obtain
\begin{equation}\label{S_nindep}
        \E\, [K_n(R_n)k(\beta_n)]
    \; = \;
    \E \Big[ \int_0^\infty\int_A\, K_n(S_n+s)\,    k(a)\,q(\alpha_n,da)\,\lambda_n
    e^{-\lambda_ns}\,ds\Big].
\end{equation}
 Using  again the independence of $V_n$  and
$(\hat\calg_\infty, \beta_1,\ldots,\beta_{n-1},R_1,\ldots, R_{n-1},S_n)$ we also have
\beqs
& &   \E \Big[ \int_0^\infty\int_A\, K_n(S_n+s)\,    k(a)\,q(\alpha_n,da)\,
    \lambda_n\,1_{V_n\ge s}\, ds \Big] \\
 &=&     \E \Big[ \int_0^\infty\int_A\, K_n(S_n+s)\,    k(a)\,q(\alpha_n,da)\,  \lambda_n\,
     \P(V_n\ge s)\, ds\Big]
 \enqs
and since  $\P(V_n\ge s)=e^{-\lambda_ns}$, this coincides with the
right-hand side of \eqref{S_nindep}. By a change of variable we arrive at
equality  \eqref{thesisrewritten}:
\[
\begin{array}{lll}\dis
        \E\, [K_n(R_n)k(\beta_n)]
    & = &\dis
    \E \Big[ \int_{S_n}^\infty\int_A\, K_n(t)\,    k(a)\,q(\alpha_n,da)\,
\lambda_n\,    1_{V_n\ge t-S_n}\, dt \Big]
    \\
& =&\dis
    \E \Big[ \int_{0}^\infty\int_A\, K_n(t)\,    k(a)\,q(\alpha_n,da)\,
  \lambda_n\,   1_{S_n<t\le R_n}\,dt\Big].
    \end{array}
\]
This concludes the proof of Lemma \ref{MPPperturbed}.    \qed

\vspace{3mm}

It follows from this lemma that
the $\H^m$-compensator of   $\kappa^m$ under $\Q$
is absolutely continuous with respect to $\lambda(da)\,dt$ and it can be written in the form
  $$
 \tilde\kappa^m(dt,da) \; = \;  \phi^m_t(a)\, \lambda(da)\,dt
 $$
for a suitable nonnegative $\calp(\H^m)\otimes \calb(A)$-measurable  function $\phi^m$
which is bounded on $\hat\Omega\times [0,T]\times A$. Indeed,
from the choice of the kernel $q^m(b,da)$ we obtain
\beqs
\phi^m_t(a) &=& \sum_{n\ge 1}  1_{(S^m_{n},R^m_n]}(t)\,
\frac{1}{\lambda(\bfB(\alpha_n,1/m))}\, 1_{\bfB(\alpha_n,1/m)}(a) \lambda_{nm}.
\enqs
which is bounded on $\hat\Omega\times [0,T]\times A$
since each $\alpha_n$ takes only a finite number of values and
$S^m_{N}> t_N=T$, so that the values of $\phi_t^m$ on $[0,T]$
only depend on the first $N-1$ summands.

In the final step of the proof we will modify the random
measure $\kappa^m$ by adding an independent Poisson process $\pi^k$
with ``small'' intensity. This will not affect too much the
$\tilde\rho$-distance between the corresponding trajectories and will
produce a random measure whose compensator
remains absolutely continuous with respect to the measure
$\lambda(da)\,dt$ and has a bounded density which,  in addition, is bounded
away from zero.

Recall that on the space
$(\Omega',\calf',\P')$ we assumed that for every integer $k\ge 1$ there
existed    a Poisson random measure $\pi^k$  on $(0,\infty)\times A$, admitting compensator $k^{-1}\lambda(da)\,dt$ with respect to its natural filtration.
We will consider $\pi^k$  as defined in $(\hat\Omega,\hat\calf)$.
Each
$\pi^k$ has the form
\beqs
\pi^k &=& \sum_{n\ge 1}\delta_{(T_n^k , {\xi_n^k })},
\enqs
for a marked point process $(T_n^k , {\xi_n^k })_{n\ge1}$
on $(0,\infty)\times A$,
and we denote  $\F^{\pi^k}$ $=$ $(\calf^{\pi^k}_t)$ its natural filtration.
Let us define another random measure setting
$$
\mu^{km} =\kappa^m+\pi^k.
$$
Note that the jumps times $(R_n^m)_{n\ge1}$  are independent
of the jump times $(T_n^k )_{n\ge1}$, and the latter have
absolutely continuous laws. It follows that, except possibly
on a set of $\Q$ probability zero, their graphs  are
disjoint, i.e. $\kappa^m$ and $\pi^k $ have no common jumps. Therefore,
  the random measure $\mu^{km} $
and its associated pure jump process (denoted $ I^{km}$)
    admit a representation
$$
\mu^{km} =\sum_{n\ge 1}\delta_{(S_n^{km} , {\eta^{km} _n})},
\qquad
 I^{km} _t =
 \sum_{n\ge 0} \eta^{km} _{n}1_{ [S_n^{km} ,S^{km} _{n+1})}(t), \;\;\; t \geq 0,
$$
where $\eta^{km} _{0}=a_0$,
$(S_n^{km} ,\eta^{km} _n)_{n\ge 1}$ is a marked point process,  each $S_n^{km} $ coincides with one of the times $R_n^m$ or one of the times
$T_n^k $, and each $\eta_n^{km} $ coincides with one of the random variables $\xi_n^k $ or one of the random variables $\beta_{n}^m$.
We claim that, for large $k$, $I^{km}$ is close to $\hat\alpha^n$ with respect to the metric $\tilde \rho$,
namely that
\begin{equation}\label{approxvalue}
\tilde \rho(I^{km},\hat\alpha^m)\to 0
\end{equation}
as  $k\to\infty$.
To prove this claim it suffices
to prove that $I^{km}\to \hat\alpha^m$ in $dt\otimes d\Q$-measure.
Recall that the jump times of $\pi^k $ are denoted $T_n^k $.
Since  $T_1^k $ has exponential law with parameter
$ \lambda(A)/k$ the event $B_k =\{T_1^k  >T\}$ has probability
$e^{-\lambda(A)T/k}$, so that $\Q(B_k)\to 1$ as $k\to \infty$.
Noting that  on the set $B_k$, we have $\hat\alpha^m_t =
\alpha_0=a_0=\eta_0^{km}=
I^{km} _t$ for all $t\in [0,T]$,
 the claim
 \eqref{approxvalue} follows immediately. We will fix from now on an integer
$k$ so large that
\beq\label{rhotildetre}
\tilde\rho (\hat\alpha^m,I^{km})<\delta/3.
\enq

Having fixed both $m$ and $k$ we now define, for $n\ge 0$,
$$
\hat S_n= S_n^{km},\quad
\hat \eta_n =\hat \eta_n^{km},
\quad
\hat \mu=\sum_{n\ge1}\delta_{(\hat S_n,\hat \eta_n)},
\quad
\hat I_t=\sum_{n\ge 0}\hat  \eta_{n}1_{ [\hat S_n,\hat S_{n+1})}(t),
$$
so that the random measure $\hat\mu$ and the associated process
$\hat I$ coincide with
$\mu^{km}$ and
$I^{km}$ respectively.
The inequalities \eqref{rhotildeuno}, \eqref{rhotildedue}, \eqref{rhotildetre}
imply that $\tilde\rho (\hat\alpha,\hat I)<\delta$, which gives
\eqref{distkrylovdelta}.

To finish the proof
it remains to prove \eqref{compensatoreBmu}-\eqref{compensatoreBmubounds}.
We first note that,
since $\kappa^m$ and $\pi^k $ are independent, it is easy to prove
 that
 $\hat\mu=\mu^{km} $ has compensator  $(\phi^m_t(a)+k^{-1})\,\lambda(da)\,dt$ with respect
 to the filtration $\H^m\vee\F^{\pi^k}$ $:=$ $(\calh^m_t\vee\calf^{\pi^k }_t)_{t\ge0}
 =(\hat\calg_t\vee \calf^{\kappa^m}_t\vee \calf^{\pi^k }_t)_{t\ge0}$.
Let $\F^{\hat\mu}=(\calf_t^{\hat\mu})_{t\geq 0}$ denote
 the natural filtration of
$\hat\mu$ and let $\hat\G\vee \F^{\hat\mu}$ be the filtration $(\hat\calg_t\vee\calf_t^{\hat\mu})_{t\geq 0}$, which is smaller
than $\H^m\vee\F^{\pi^k}$.
We wish to compute the compensator of $\hat\mu$ with respect to $\hat\G\vee \F^{\hat\mu}$
under $\Q$.
To this end, consider the measure space
$([0,\infty)\times\Omega\times A, \calb([0,\infty))\otimes\calf\otimes \calb(A),dt\otimes\Q(d\omega)\otimes\lambda(da))$.
Although this is not a probability space, one can define in a standard way the conditional expectation
 of any positive measurable
function, given an arbitrary sub-$\sigma$-algebra.
Let us denote by  $\hat\nu_t(\hat\omega,a)$ the conditional expectation
of the  random field $ \phi_t^m(\hat \omega,a)+k^{-1}$
 with respect to the $\sigma$-algebra $\calp(\hat\G\vee \F^{\hat\mu})
 \otimes \calb(A)$. It is then easy to verify that
the compensator of $\hat\mu $ with respect to
$\hat\G\vee \F^{\hat\mu}$ coincides with $\hat\nu$.
Moreover, since  $ \phi_t^m(\hat \omega,a)$ is nonnegative and bounded on
$\hat\Omega\times [0,T]\times A$,
we can take a version of $\hat\nu$ satisfying
$$
 k^{-1}\le \inf_{\hat\Omega\times [0,T]\times A}\hat\nu\le
  \sup_{\hat\Omega\times [0,T]\times A}\hat\nu<\infty.
  $$
Now \eqref{compensatoreBmu}-\eqref{compensatoreBmubounds}
are proved and the proof of Proposition
\ref{extensionapproximation}
is finished.
\qed

\vspace{9mm}

\small

\bibliographystyle{plain}

\bibliography{biblio6}

\end{document}